\documentclass[11pt]{article}

\pdfoutput=1

\usepackage{color}
\usepackage{epsfig}
\usepackage{afterpage}
\usepackage{amsmath,amssymb}

\newcommand{\proof}{\noindent{\bf Proof:} \quad}
\newcommand{\qed}{\hfill\mbox{\raggedright\rule{0.07in}{0.1in}}\vspace{0.1in}}

\newtheorem{theorem}{Theorem}[section]

\newtheorem{definition}[theorem]{Definition}

\newtheorem{lemma}[theorem]{Lemma}

\begin{document}

\begin{center}
\section*{Regions of Stability for a Linear Differential Equation with
   Two Rationally Dependent Delays}
\end{center}

\begin{center}
  Joseph M. Mahaffy \\
  \begin{small}
    Nonlinear Dynamical Systems Group \\
    Department of Mathematics \\
    San Diego State University \\ San Diego, CA 92182, USA \\
  \end{small}
  \vspace{0.15in}
  Timothy C. Busken \\
  \begin{small}
    Department of Mathematics \\
    Grossmont College \\ El Cajon, CA 92020, USA \\
  \end{small}
\end{center}

\begin{abstract}
Stability analysis is performed for a linear differential equation with two delays. Geometric arguments show that when the two delays are rationally dependent, then the region of stability increases. When the ratio has the form $1/n$, this study finds the asymptotic shape and size of the stability region. For example, a delay ration of $1/3$ asymptotically produces a stability region 44.3\% larger than any nearby delay ratios, showing extreme sensitivity in the delays. The study provides a systematic and geometric approach to finding the eigenvalues on the boundary of stability for this delay differential equation. A nonlinear model with two delays illustrates how our methods can be applied.
\end{abstract}

\noindent
{\bf Keywords:} Delay differential equation; bifurcation; stability analysis; exponential polynomial; eigenvalue
\vspace{0.1in}

\begin{small}
\noindent
{\bf Submitted:} 7/24/2013
\end{small}

\section{Introduction}

Delay differential equations (DDEs) are used in a variety of applications,
and understanding their stability properties is a complex and important
problem. The addition of a second delay significantly increases the
difficulty of the stability analysis. E.\ F.\ Infante \cite{Inf} stated that
an economic model with two delays, which are rationally related, has
a region of stability that is larger than one with delays nearby that are
irrationally related. This meta-theorem inspires much of the work below,
where we examine the linear two-delay differential equation:
\begin{equation}
 \dot{y}(t) + A\,y(t) + B\,y(t-1) + C\,y(t-R) = 0, \label{DDE2}
\end{equation}
as the parameters $A$, $B$, $C$, and $R \in (0,1)$ vary. (Note that
time has been scaled to make one delay unit time). Our efforts concentrate
on the stability region near delays of the form $R = \frac{1}{n}$ with $n$
a small integer. The stability analysis of Eqn.~(\ref{DDE2}) for the cases
$R = \frac{1}{3}$ and $\frac{1}{4}$ were studied in some detail in
Mahaffy {\it et al.} \cite{MZJa, MZJ} and Busken \cite{Busk}, and this work
extends those ideas.

Discrete time delays have been used in the mathematical modeling of many
scientific applications to account for intrinsic lags in time in the
physical or biological system. Often there are numerous stages in the
process, such as maturation or transport, which utilize multiple discrete
time delays. Some biological examples include physiological control
\cite{BEL, BeCa, CamBel}, hematopoietic systems \cite{BEM, BMM, McDc, McDd},
neural networks \cite{BCvdD,GopMo,GCC}, epidemiology \cite{CoY},
and population models \cite{BRV,MNB}. Control loops in optics \cite{MiI} and
robotics \cite{HaS} have been modeled with multiple delays. Economic
models \cite{BEMb, HOR, MACe} include production and distribution
time lags. Bifurcation analysis of these models is often quite complex.

\begin{figure}[htb]
  \centerline{
    \includegraphics[width=2.5in]{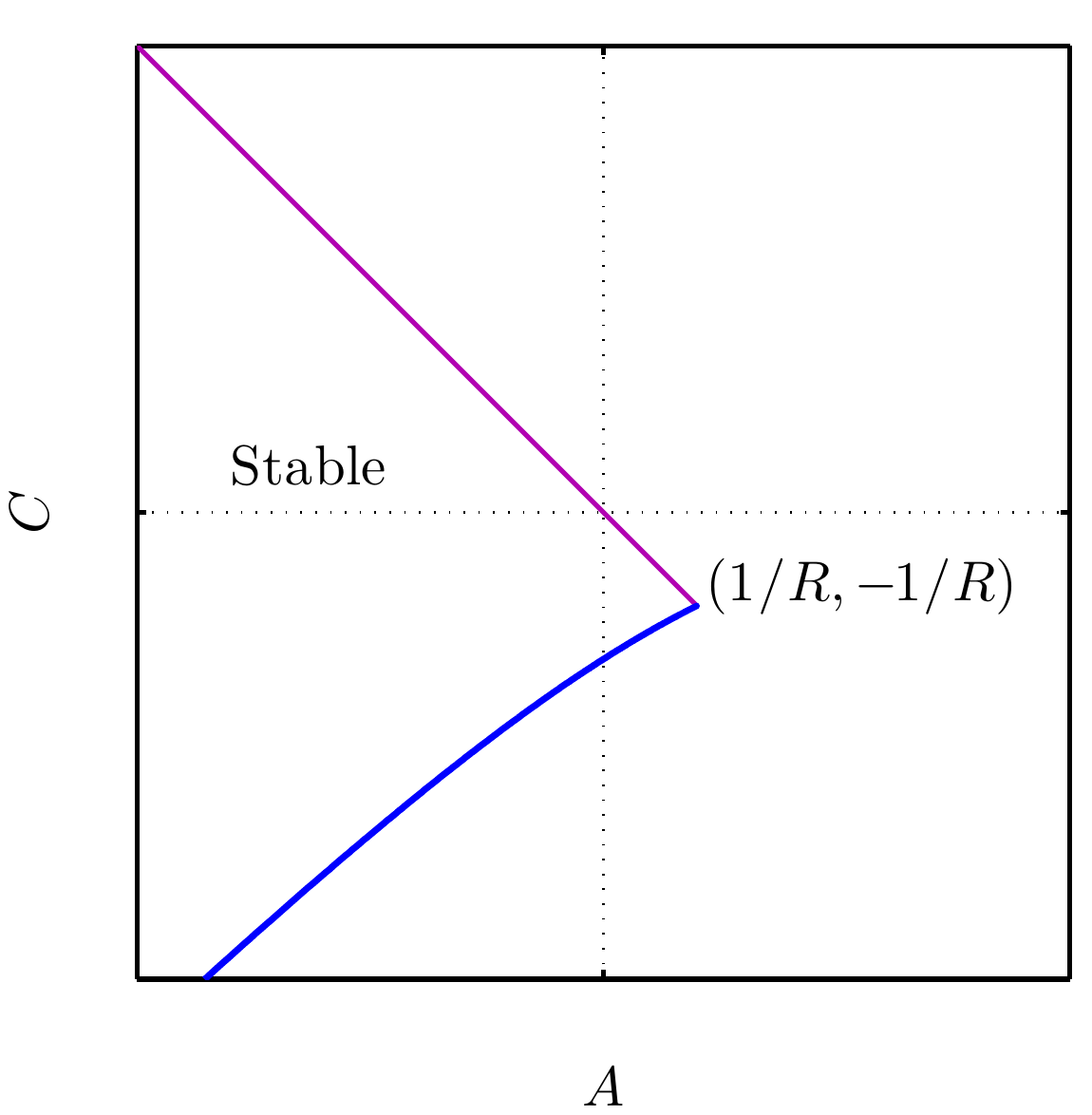}
  }
  \caption{\small{Stability region for one-delay, $\dot{y}(t) = A\,y(t)
      + C\,y(t-R)$. Violet line gives the real root crossing, while the
      blue line gives the imaginary root crossing.}}
  \label{fig:one_delay}
\end{figure}

The bifurcation analysis of the one-delay version of (\ref{DDE2}) ($B = 0$)
began with the work of Hayes \cite{HAY}. The complete stability region in the
$AC$-plane has been characterized by several authors \cite{BelC,bohay,Boese94,ELS,HIT}
with the stability boundary easily parameterized by the delay $R$ (which
can be scaled out). Fig.~\ref{fig:one_delay} shows this region of stability.
The boundary of the stability region for the two-delay equation~(\ref{DDE2})
has been studied by many researchers \cite{BEL,elsken,HAHw,HT,LRW,NUS,MZJ}.
Several authors study the special case where $A = 0$
\cite{Hal1,LEV,LRW,NUS,PIOT,RM,RUC,saka,YTJ}. Hale and Huang \cite{HAHw}
performed a stability analysis of the two-delay problem,
\begin{equation}
 \dot{y}(t) + a\,y(t) + b\,y(t-r_1) + c\,y(t-r_2) = 0, \label{dde2h}
\end{equation}
where they fixed the parameters, $a$, $b$, and $c$, then constructed the
boundary of stability in the $r_1r_2$ delay space. Braddock and van Driessche
\cite{BRV} completely determined the stability of (\ref{dde2h}) when $b = c$,
and partially extended the results outside that special case. Most of these
analyses have studied the 2D stability structure of either (\ref{DDE2}) or
(\ref{dde2h}) with one parameter equal to zero or fixing some of the parameters.
Often the 2D analyses result in observing disconnected stability regions
for (\ref{dde2h}). Elsken \cite{elsken} has proved that the stability region
of (\ref{dde2h}) is connected in the $abc$-parameter space with fixed $r_1$
and $r_2$. Recently, Bortz \cite{bortz} developed an asymptotic expansion using
Lambert W functions to efficiently compute roots of the characteristic equation
for (\ref{dde2h}) with some restrictions.

Mahaffy {\it et al.} \cite{MZJa, MZJ} studied (\ref{DDE2}) for specific
values of $R$, examining the 2D cross-sections for fixed $A$ and developing
3D bifurcation surfaces in the $ABC$-parameter space. The work below shows why
delays of the form $R = \frac{1}{n}$ have enlarged regions of stability.
\vspace{0.25in}

\setcounter{equation}{0}
\setcounter{theorem}{0}
\setcounter{figure}{0}
\setcounter{table}{0}
\section{Motivating Example}

\vspace{.15in}

In 1987, B{\'e}lair and Mackey \cite{BEM} developed a two delay model for platelet
production. The time delays resulted from a delay of maturation and another delay
representing the finite life-span of platelets. The resulting numerical simulation
for certain parameters produced fairly complex dynamics. Here we examine a slight
modification of their model and demonstrate the extreme sensitivity of the model
behavior near rationally dependent delays.

\begin{figure}[htb]
\begin{center}
\begin{tabular}{cc}
\includegraphics[width=0.45\textwidth]{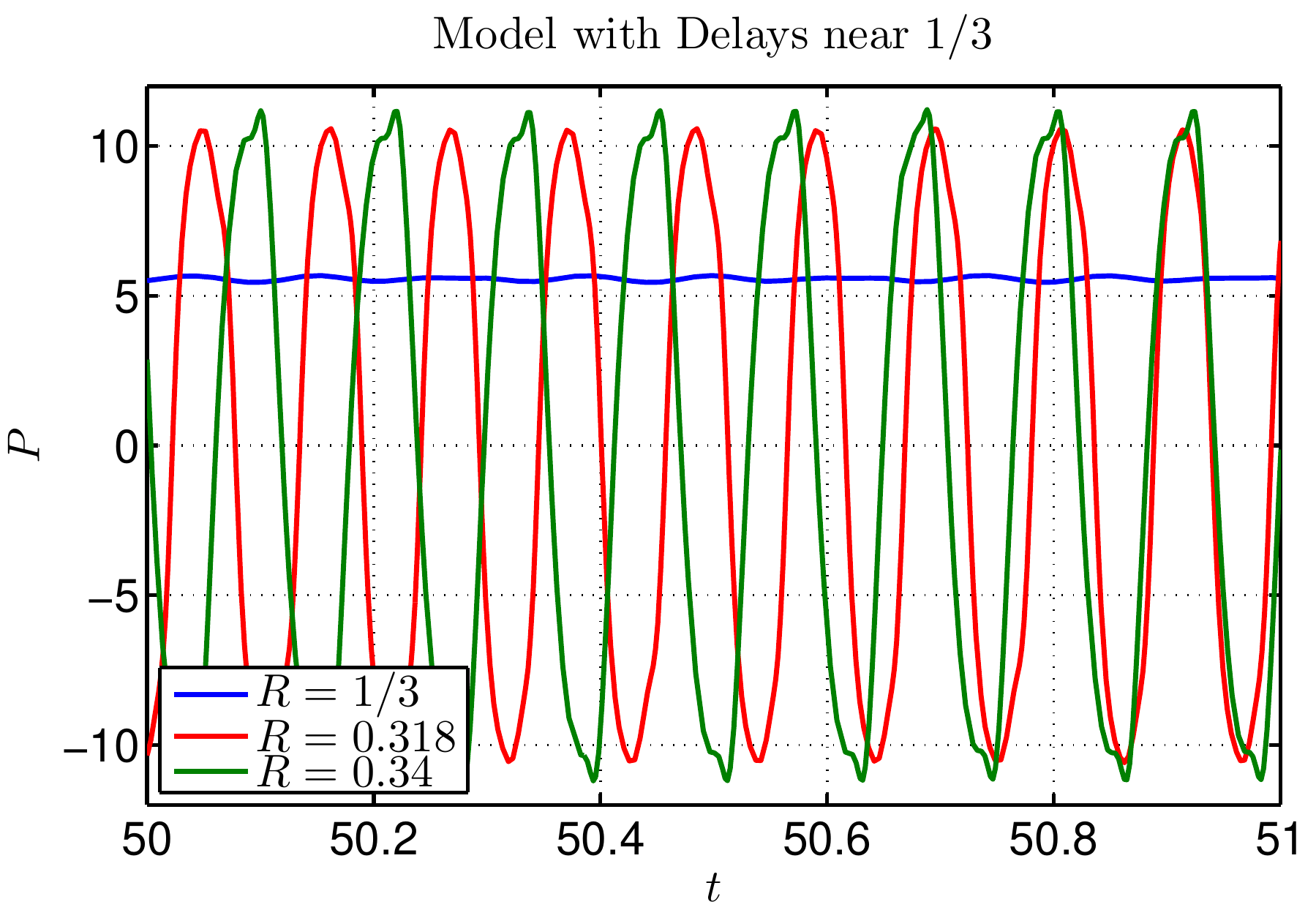}
&
\includegraphics[width=0.45\textwidth]{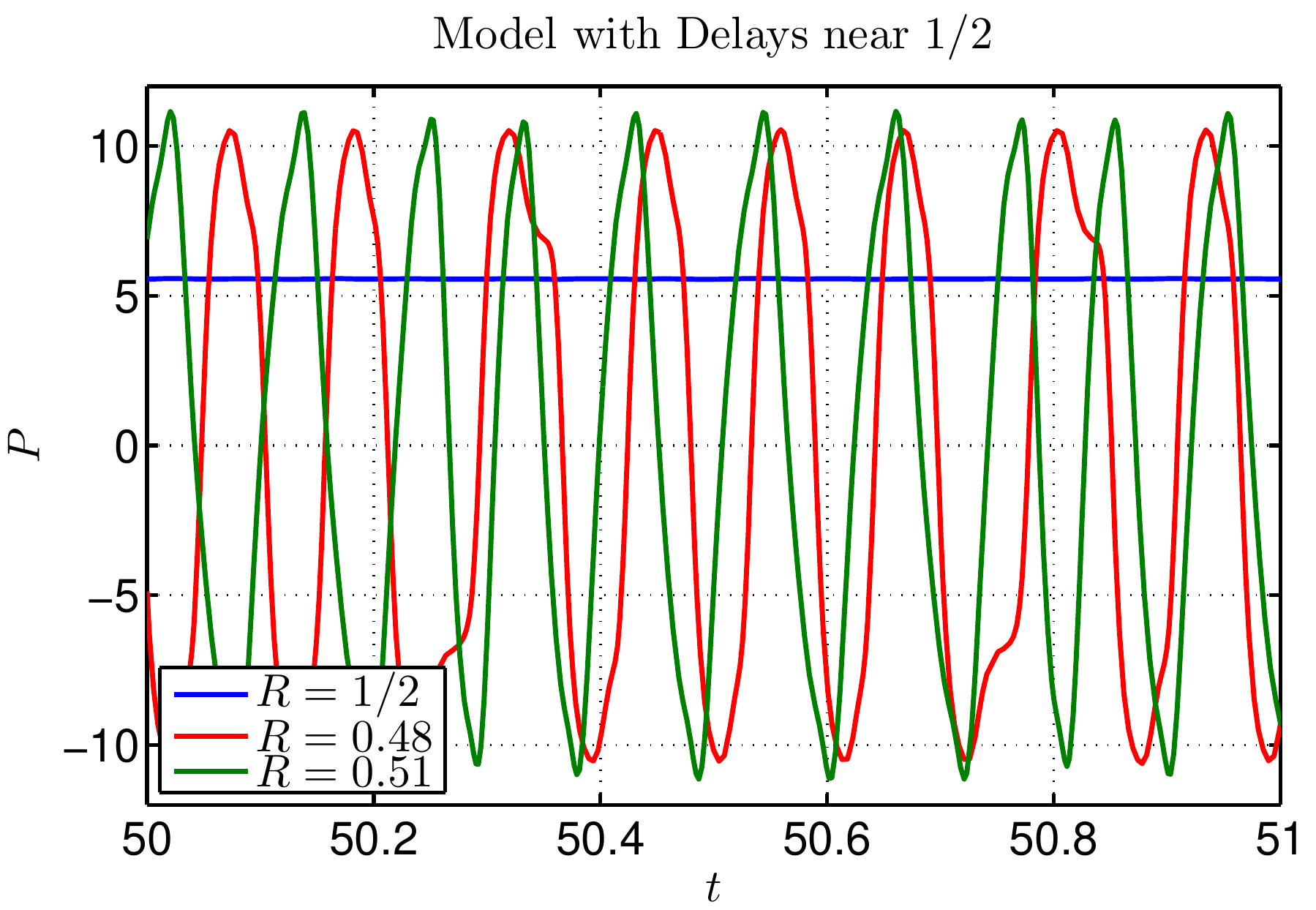}

\end{tabular}
\caption{Modified platelet model \cite{BEM}. Simulations show the model sensitivity to the delays $R = \frac{1}{3}$ and $\frac{1}{2}$ with nearby delays showing complex oscillations.}
\label{fig:platelet}
\end{center}
\end{figure}

The modified model that we consider is given by:
\[
 \frac{dP}{dt} = -\gamma P(t) + \beta(P(t-R)) - f\cdot\beta(P(t-1)),
\]
where $\beta(P) = \frac{\beta_0\theta^nP}{\theta^n + P^n}$. This model has a standard
linear decay term and nonlinear delayed production and destruction terms. The total
lifespan is normalized to one, while the maturation time is $R$. This model differs
from the platelet model by choosing an arbitrary fractional multiplier, $f$, instead
of having a delay dependent fraction. For our simulations we fixed $\gamma = 100$,
$\beta_0 = 168.6$, $n = 4$, $\theta = 10$, and $f = 0.35$. This gives the equilibria
$P_e = 0$ and $P_e \approx 5.565$ with $\beta'(5.565) \approx 100$.

Fig.~\ref{fig:platelet} shows six simulations near $R = \frac{1}{2}$ and $R = \frac{1}{3}$,
where the model is asymptotically stable. However, fairly small perturbations of the delay
away from these values result in unstable oscillatory solutions, as is readily seen in the
figure. The oscillating solutions are visibly complex. This paper will explain some of the
results shown in Fig.~\ref{fig:platelet}.

\setcounter{equation}{0}
\setcounter{theorem}{0}
\setcounter{figure}{0}
\setcounter{table}{0}
\section{Background}

\vspace{.15in}

\subsection{Definitions and Theorems}

There are a number of key definitions and theorems that are needed to
build the background for our study. Our analysis centers around finding
the stability of (\ref{DDE2}). Stability analysis of a linear DDE begins
with the characteristic equation, which is found in a manner similar to
ordinary differential equations by seeking solutions of the form
$y(t) = c\,e^{\lambda t}$. The characteristic equation for (\ref{DDE2})
is given by:
\begin{equation}
  \lambda + A + B\,e^{-\lambda} + C\,e^{-\lambda\,R} = 0. \label{chareqn}
\end{equation}
This is an exponential polynomial, which has infinitely many solutions, as
one would expect because a DDE is infinite dimensional. Stability occurs if
all of the eigenvalues satisfy
\[
 {\textnormal{Re}} (\lambda) \le 0, \qquad {\rm for\ all}\ \lambda.
\]
One can readily see from (\ref{chareqn}) that the $A + B + C = 0$ plane, $\Lambda_0$,
provides one boundary where a real eigenvalue $\lambda$ crosses between
positive and negative, so creates a bifurcation surface.

To understand what is meant by rationally dependent delays resulting in
larger regions of asymptotic stability, we need an important theorem about
the minimum region of stability for (\ref{DDE2}).

\begin{theorem} \textbf{Minimum Region of Stability (MRS)}\label{thm1}
For $A>|B|+|C|$, \,all solutions $\lambda$ to Eqn.~(\ref{chareqn}) have
${\textnormal{Re}}(\lambda) < 0$, \,which implies that Eqn.~(\ref{DDE2}) is asymptotically
stable inside the pyramidal-shaped region centered about the positive $A$-axis,
independent of $R$.
\end{theorem}

The proof of the MRS Theorem can be found in both Zaron \cite{ZAR} and
Boese  \cite{BOE}. Note that one face of this MRS is formed by the plane
$A + B + C = 0$, $\Lambda_0$, where the zero root crossing occurs. The other way that
(\ref{DDE2}) can lose stability is by roots passing through the imaginary
axis or $\lambda = i\omega$. This is substituted into (\ref{chareqn}).
Since the real and imaginary parts are zero, we obtain a parametric
representation of the bifurcation curves for $B(\omega)$ and $C(\omega)$.
These are given by the expressions:
\begin{eqnarray}
  B(\omega) & = & \frac{A\sin(\omega\,R)+\omega\cos(\omega\,R)}
                  {\sin(\omega(1-R))},  \label{bifB} \\
  C(\omega) & = & -\frac{A\sin(\omega)+\omega\cos(\omega)}
                  {\sin(\omega(1-R))}, \label{bifC}
\end{eqnarray}
where $\frac{(j-1)\pi}{1-R} < \omega < \frac{j\pi}{1-R}$, and
$j\in\mathbb{Z^+}$. Clearly, there are singularities for $B(\omega)$ and
$C(\omega)$ at $\omega = \frac{j\pi}{1-R}$. This leads to the following
definition for bifurcation surfaces.

\begin{definition}
When a value of $R$ in the interval $(0,1)$ is chosen,
\textit{\textbf{Bifurcation Surface j}}, $\Lambda_j$, is determined by
Eqns.~(\ref{bifB}) and (\ref{bifC}), and is defined parametrically
for $\frac{(j-1)\pi}{1-R}<\omega<\frac{j\pi}{1-R}$ and $A\in\mathbb{R}$.
This creates a separate parameterized surface representing solutions
of the characteristic equation, (\ref{chareqn}), $\lambda = i\omega$,
which can be sketched in the $ABC$ coefficient-parameter space of
(\ref{DDE2}), for each positive integer, $j$.
\end{definition}

Because the MRS is centered on the $A$-axis, we often choose to fix $A$
and view the cross-section of the bifurcation surfaces. Thus, we have the
related definition:

\begin{definition}
\textit{\textbf{Bifurcation Curve j}}, $\Gamma_j$, is determined by Eqns.~(\ref{bifB})
and (\ref{bifC}) and is defined parametrically for
$\frac{(j-1)\pi}{1-R} < \omega < \frac{j\pi}{1-R}$ with the values of $R$ and $A$ fixed.
This creates a parametric curve, which can be drawn in the $BC$-plane for each $j$.
\end{definition}

\begin{figure}[htb]
\begin{center}
\begin{tabular}{ccc}
\includegraphics[width=0.3\textwidth]{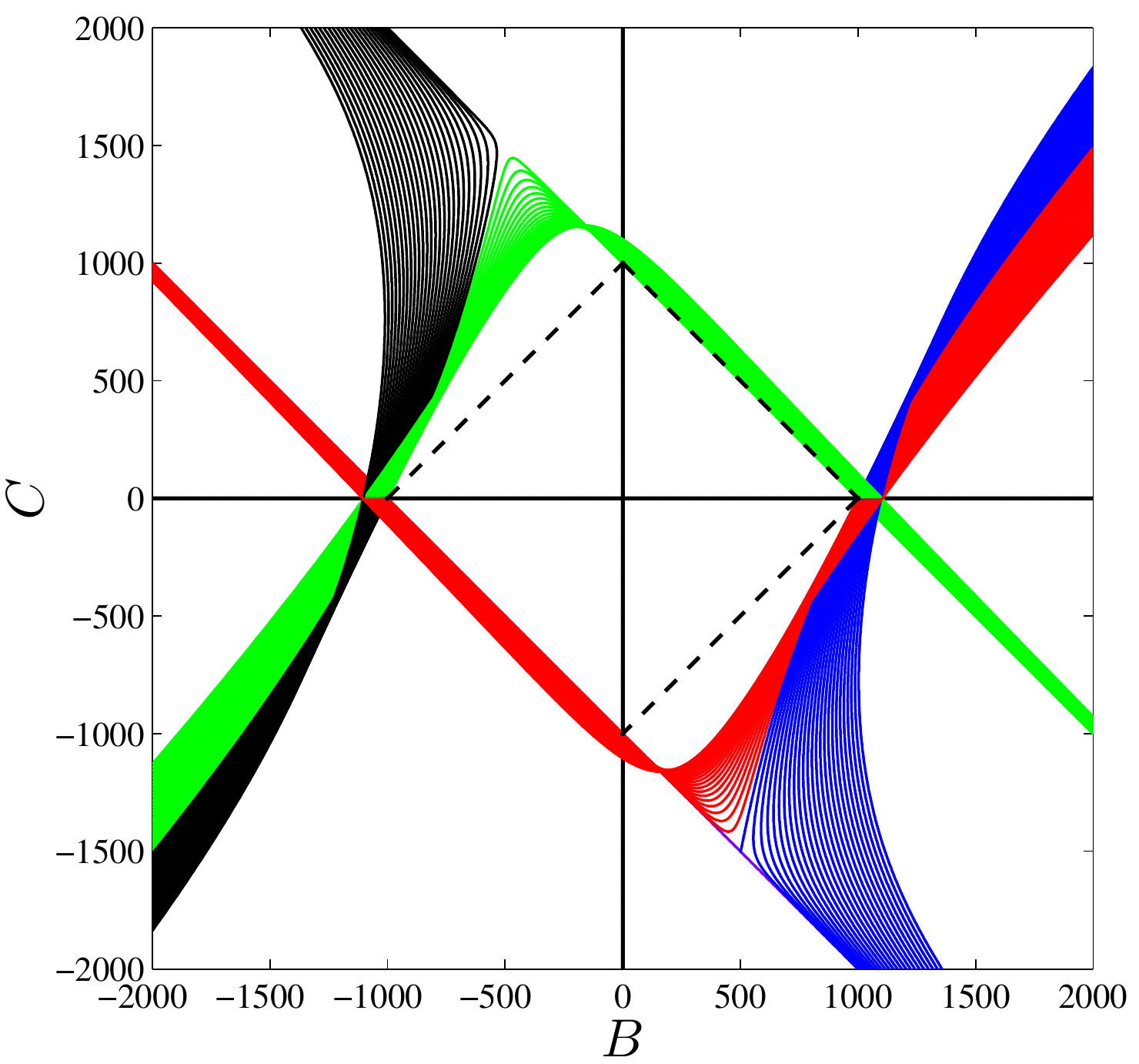}
&
\includegraphics[width=0.3\textwidth]{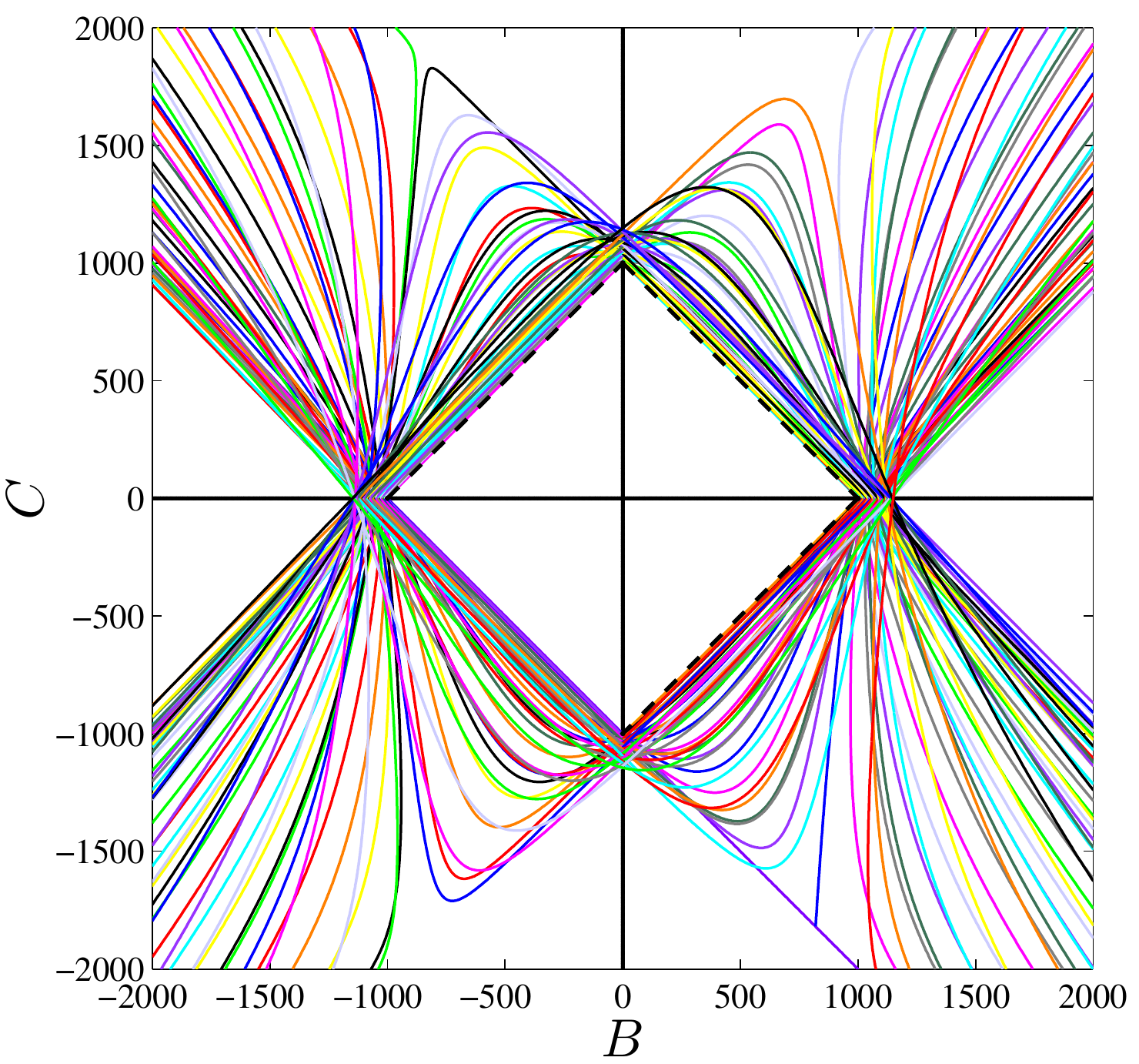}
&
\includegraphics[width=0.3\textwidth]{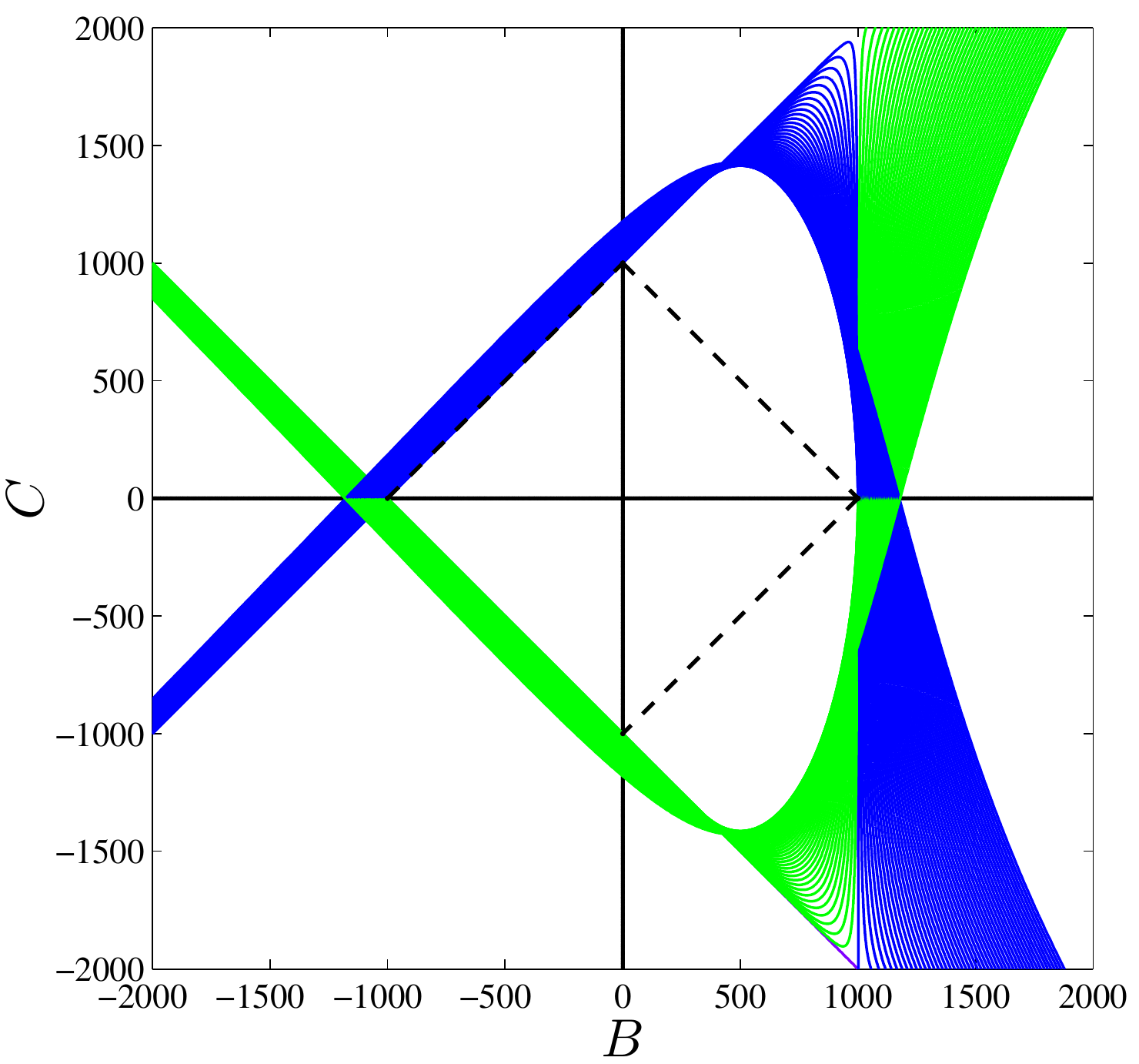} \\
 $R = \frac{1}{3}$ &  $R = 0.45$  &  $R = \frac{1}{2}$
\end{tabular}
\caption{Bifurcation curves: Shows first 100 parametric curves in the
$BC$-parameter space for $A = 1000$ and various delays. The dashed curve shows the
boundary of the MRS.}
\label{fig:bif_curves}
\end{center}
\end{figure}

For most values of A, the bifurcation curves in the $BC$ plane generated by
(\ref{bifB}) and (\ref{bifC}) tend to infinity parallel to the lines $B+C=0$
or $B-C=0$. When $A$ and $R$ are fixed, one can show using the partitioning method of d'El'sgol'Ts \cite{ELS} that a finite number of bifurcation curves will
intersect in the $BC$ parameter space to form the remainder of the boundary
of the stability region not given by part of the real root-crossing plane. It is
along this boundary where eigenvalues of (\ref{chareqn}) cross the imaginary
axis in the complex plane. Fig.~\ref{fig:bif_curves} gives three examples showing
the first 100 bifurcation curves for $A = 1000$ and delays of $R = \frac{1}{3}$,
$R = 0.45$, and $R = \frac{1}{2}$. Assuming that 100 bifurcation curves give a
good representation of the stability region, this figure shows how different
the regions of stability are for the different delays. The stability region
for $R = \frac{1}{2}$ is significantly greater than the others, and the stability
region of $R = 0.45$ is very close to the MRS.

As seen in Fig.~\ref{fig:bif_curves}, the bifurcation curves can intersect often
along the boundary of the region of stability, which creates challenges in
describing the evolution of the complete stability surface for (\ref{DDE2})
in the $ABC$-parameter space. We need to discuss how we construct the 3D
bifurcation surface using a few more defined quantities.

Mahaffy {\it et al.} \cite{MZJ} proved that if $R > R_0$ for $R_0 \approx
0.012$, then the stability surface comes to a point with a smallest value, $A_0$.

\begin{theorem}[Starting Point]  \label{Ao}
If $R > R_{0}$, then the stability surface for Eqn.~(\ref{DDE2}) comes to a
point at
$(A_{0},B_{0},C_{0}) = \left(-\frac{R+1}{R}_{,} \  \frac{R}{R-1}_{,} \
\frac{1}{R(1-R)}\right)$, and Eqn.~(\ref{DDE2}) is unstable for $A < A_{0}$.
\end{theorem}

For some range of $A$ values with $A > A_0$, the stability surface is exclusively
composed of $\Lambda_1$ and $\Lambda_0$. As $\omega \to 0^+$, $\Lambda_1$
intersects $\Lambda_0$.
The surface $\Lambda_1$ bends back and intersects $\Lambda_0$ again, enclosing
the stability region. As $A$ increases, $\Lambda_2$ approaches $\Lambda_1$,
and at least for a range of $R$, $\Lambda_2$
self-intersects. In the 2D $BC$-parameter plane, this creates a disconnected
stability region, which later joins the main bifurcation surface emanating from
the starting point. The $A$ value, where this self-intersecting bifurcation
surface joins, is the $1^{st}$ transition. Transitions are one of the most
important occurrences that affect the shape of the bifurcation surface.

\begin{figure}[htb]
\begin{center}
\begin{tabular}{ccc}
\includegraphics[width=0.3\textwidth]{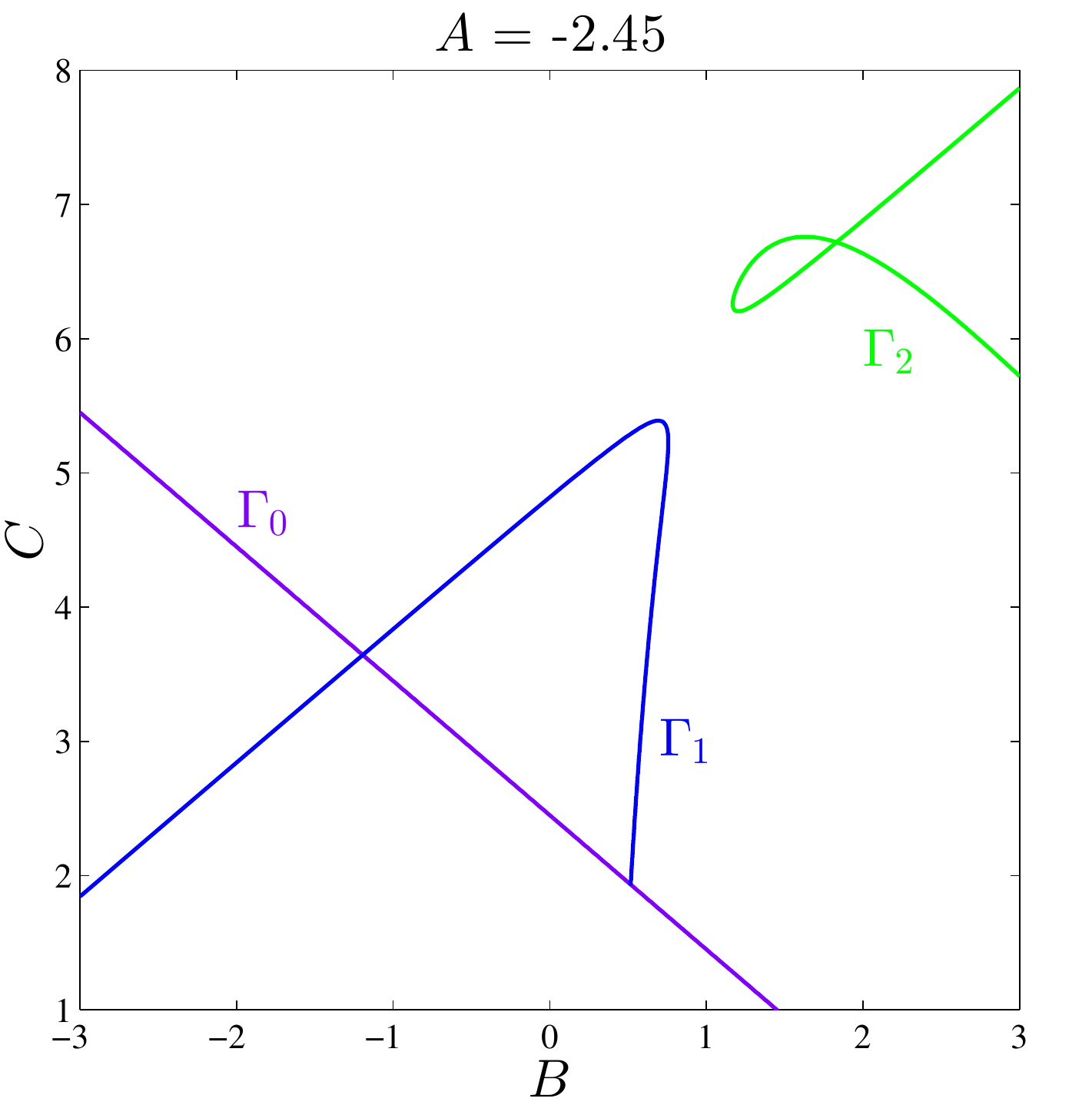}
&
\includegraphics[width=0.3\textwidth]{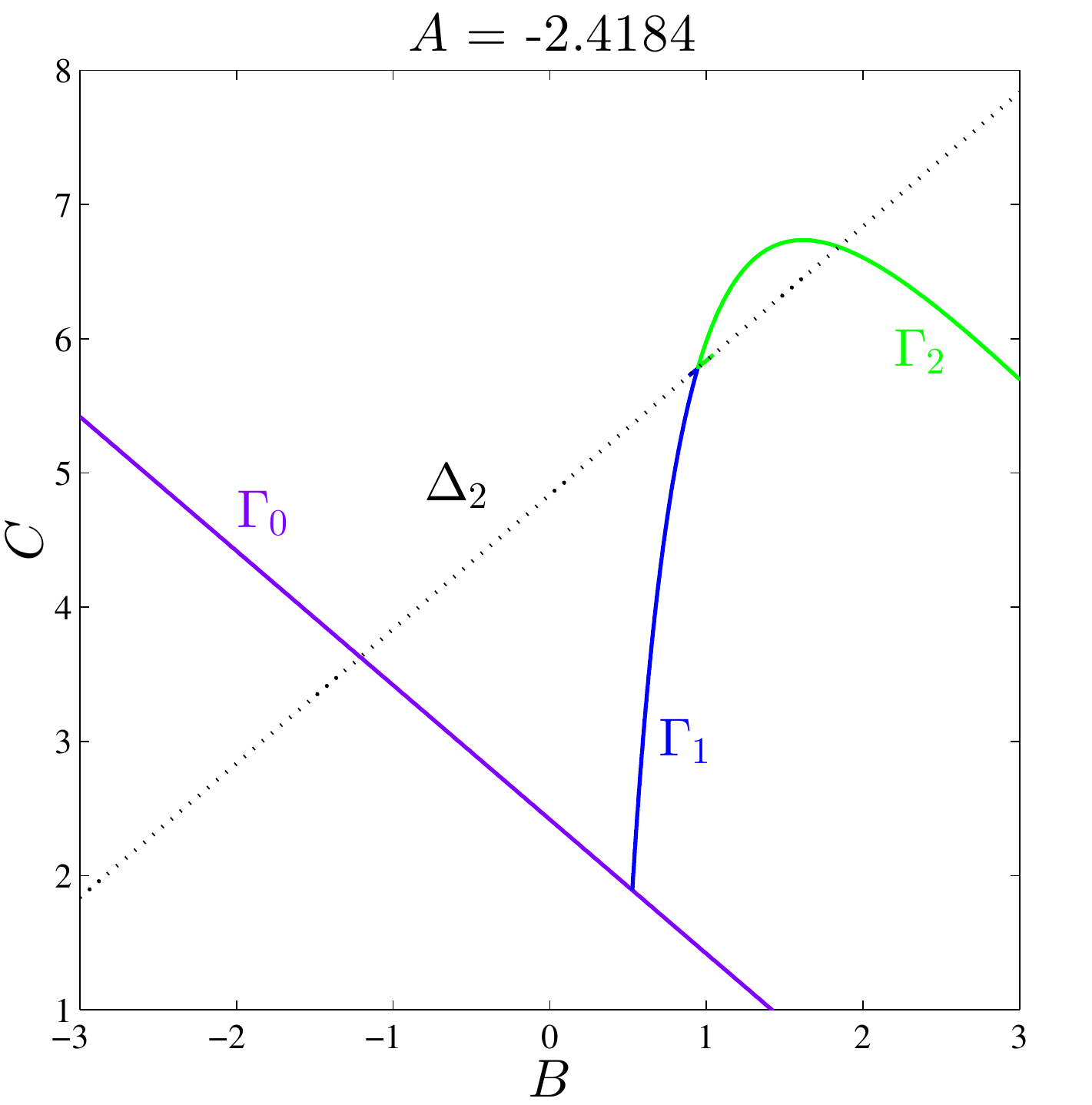}
&
\includegraphics[width=0.3\textwidth]{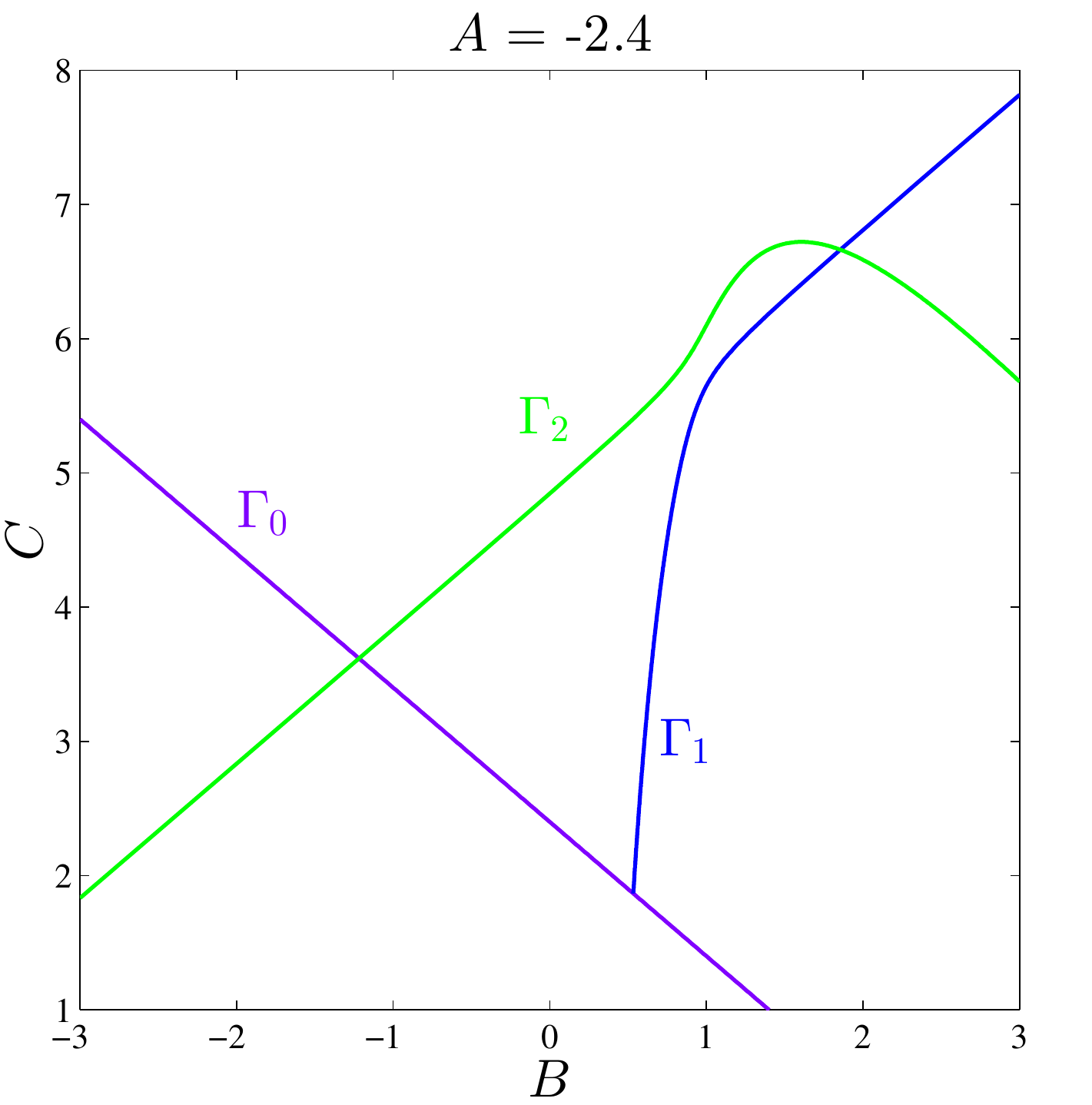}
\end{tabular}
\caption{Transition: Enclosed regions in the figures above
are the regions of stability for $R = 0.25$ as $A$ moves through a transition and
the stability spur joins the main stable surface. }
\label{fig:spur}
\end{center}
\end{figure}

\begin{definition}[Transition and Degeneracy Line]
There are critical values of $A$ corresponding to where Eqns.~(\ref{bifB}) and
(\ref{bifC}) become indeterminate at $\omega=\frac{j\pi}{1-R}$. These
\textbf{\textit{transitional}} values of $A$ are denoted by $A^{*}_{j}$, where
\begin{equation}  \label{atrans}
 A^{*}_{j}  = \textstyle{-\left(\frac{j\pi}{1-R}\right)
            \cot\left(\frac{jR\pi}{1-R}\right)}, \ j = 1,2,... \ .
\end{equation}
At a \textbf{\textit{transition}},  Curves $j$ and $(j+1)$  coincide at the
specific point $(B^{*}_{j}, C^{*}_{j})$, where
\begin{equation*}
  B^{*}_{j}  = (-1)^{j}\frac{(1-R)\cos\left(\frac{jR\pi}{1-R}\right)
            - jR\pi\csc\left(\frac{jR\pi}{1-R}\right)}{(1-R)^{2}} \\[1.0ex]
\end{equation*}
\begin{equation}  \label{bctrans}
  C^{*}_{j}  = -(-1)^{j}\frac{(1-R)\cos\left(\frac{j\pi}{1-R}\right)
               -j\pi\csc\left(\frac{j\pi}{1-R}\right)}{(1-R)^{2}}
             = \frac{j\pi\csc\left(\frac{jR\pi}{1-R}\right)-
               (1-R)\cos\left(\frac{jR\pi}{1-R}\right)}{(1-R)^{2}}\\[1.0ex]
\end{equation}
All along the \textbf{\textit{degeneracy line}}, $\Delta_j$,
\begin{equation}\label{deg}
  (B-B^{*}_{j})  +(-1)^{j}(C-C^{*}_{j})=0, \hspace{0.2in} A^{*}_{j},
\end{equation}
there are two roots of (\ref{chareqn}) on the imaginary axis with
$\lambda=\frac{j\pi}{1-R}i$.
\end{definition}

If $\Lambda_j$ is on the boundary of the stability region for $A$
slightly less than $A^{*}_{j} $, then $\Delta_j$ becomes part of the
stability region's boundary, at Transition $j$. Subsequently,
$\Lambda_{j+1}$ enters the boundary of the stability region. These transitions
create the greatest distortion to the stability surface and attach stability
spurs. It is important to note that many transitions occur outside the stability
surface and only affect the organization of the bifurcation curves.
Fig.~\ref{fig:spur} shows cross-sections in the $BC$-plane as $A$ goes through a
transition.

\begin{definition}[Stability Spur]
If $\Lambda_{j+1}$ self-intersects and encloses a region of stability for (\ref{DDE2})
as $A$ increases with $A_j^p$ being
the cusp or \textbf{the Starting Point of Spur $j$}, then this quasi-cone-shaped
\textbf{stability spur} has its cross-sectional area monotonically increase with
$A$ until $A$ reaches the transitional value,  $A^*_j$. For $A = A^*_j$, the
Stability Spur $j$, $Sp_j(R)$, connects with the larger stability surface, via the degeneracy
line, $\Delta_j$.
\end{definition}

There are a couple of other ways for bifurcation surfaces to enter (or leave)
the boundary of the main stability region as $A$ increases. We define these
means of altering the boundary as \textit{transferrals} and \textit{tangencies},
which relate to higher frequency eigenvalues becoming part of the boundary (or
being lost) as $A$ increases.

\begin{definition}[Transferral and Reverse Transferral]
The \textbf{transferral} value of $A = A_{i,j}^z$ \,is the value of $A$
corresponding to the intersection of $\Lambda_j$ (or $\Gamma_j$) with
$\Lambda_i$ (or $\Gamma_i$) at $\Lambda_0$.
$\Lambda_j$ (or $\Gamma_j$) enters the boundary of the stability region for $A > A_{i,j}^z$.
For some values of $R$ the stability surface can undergo a \textbf{reverse transferral}, \,$\tilde{A}_{j,i}^z$, \,which is a transferral characterized by $\Lambda_j$ (or $\Gamma_j$) leaving the boundary, \,or  a transferring \textit{back over} to $\Lambda_i$ (or $\Gamma_i$) the portion of the boundary originally taken by $\Lambda_j$ (or $\Gamma_j$) at $A_{i,j}^z (<\tilde{A}^z_{j,i})$.
\end{definition}

\begin{figure}[htb]
\begin{center}
\begin{tabular}{ccc}
\includegraphics[width=0.3\textwidth]{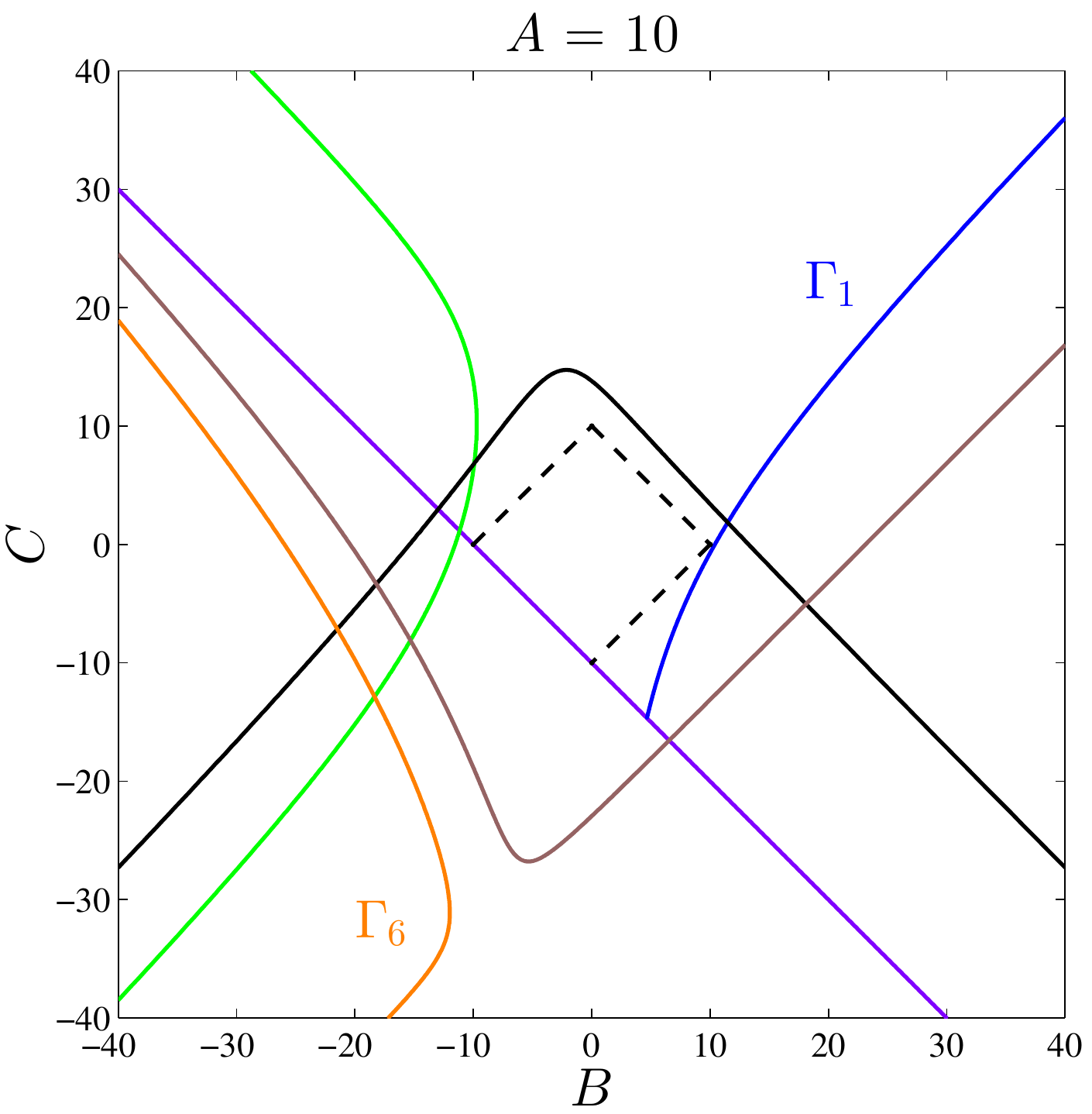}
&
\includegraphics[width=0.3\textwidth]{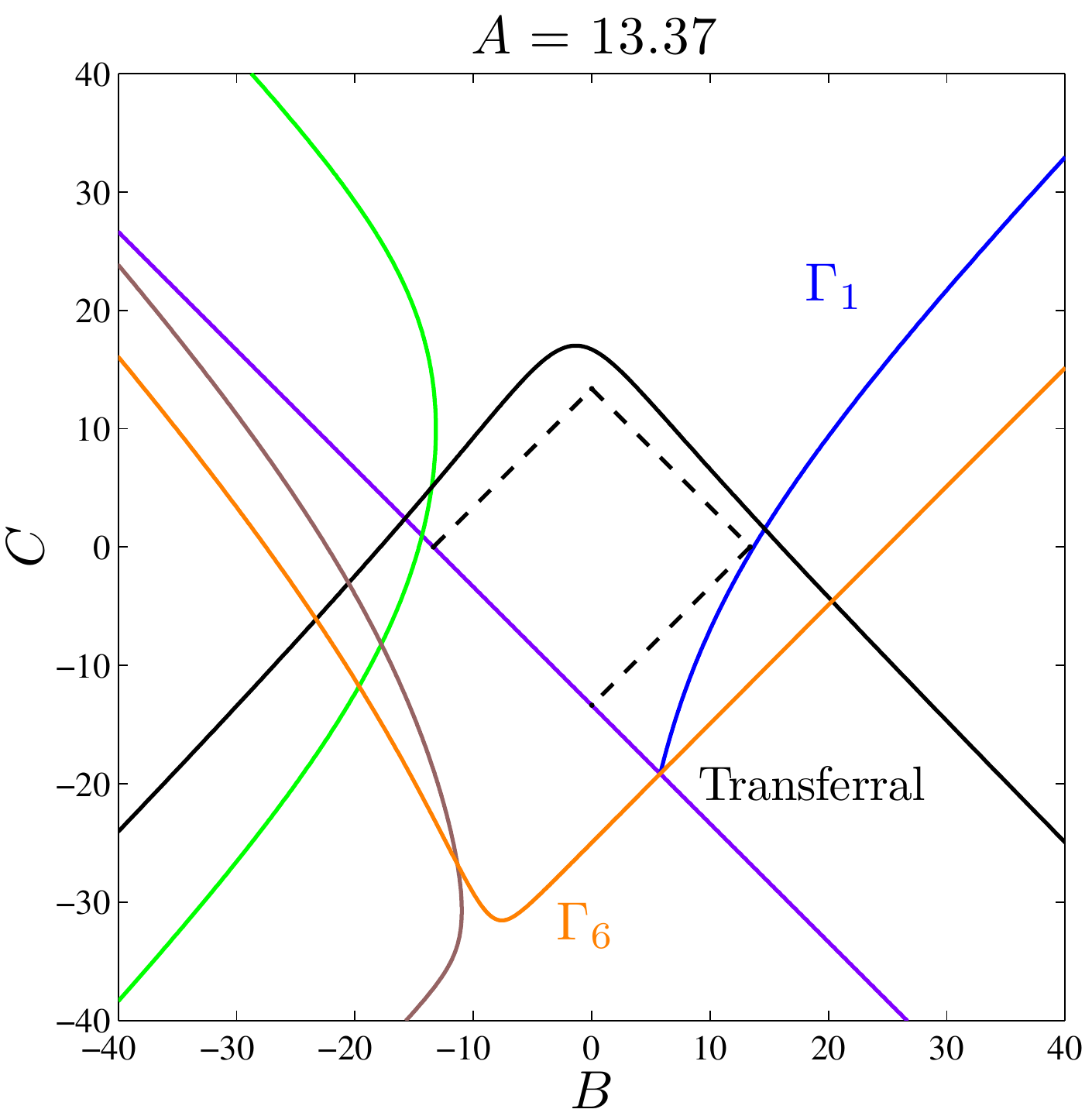}
&
\includegraphics[width=0.3\textwidth]{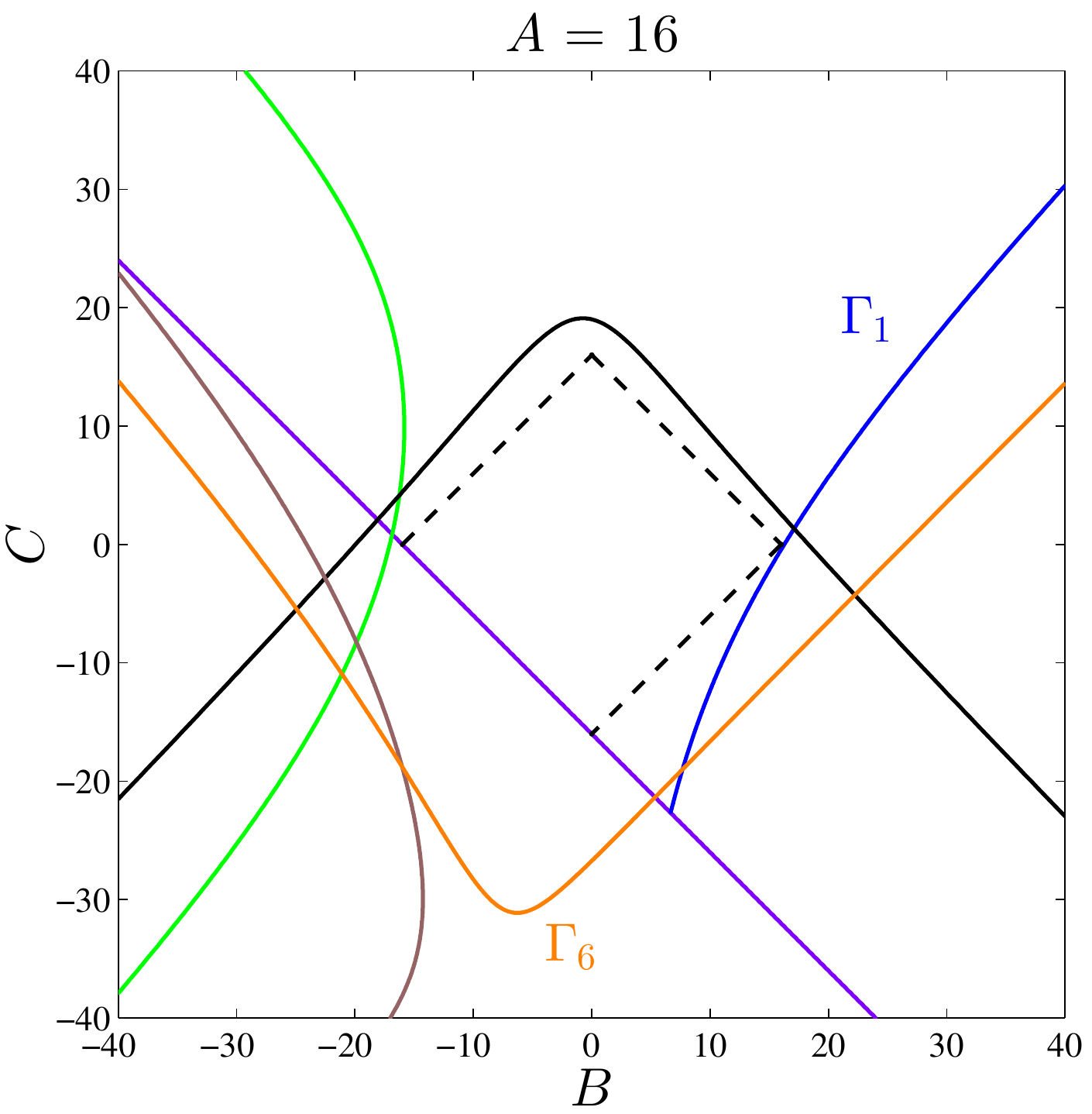}
\end{tabular}
\caption{Transferral 1, $A^z_{1,6}$. The stability boundary for $R=0.25$ located in the $BC$ plane before, during and after the first transferral is given by the closed region enveloping the MRS (bold dashed line). That closed region is formed by $\Lambda_0$ (violet), $\Gamma_1$ (blue), $\Gamma_2$ (green), $\Gamma_3$ (black), $\Gamma_5$ (dark green), and $\Gamma_6$ (orange). }
\label{fig:transfer}
\end{center}
\end{figure}

\begin{definition}[Tangency and Reverse Tangency]
The value of $A$ corresponding to the \textbf{tangency} of two surfaces $i$ and $j$ is
denoted $A_{i,j}^t$. $\Lambda_j$ (or $\Gamma_j$) becomes tangent to
$\Lambda_i$ (or $\Gamma_i$), where $\Lambda_i$ (or $\Gamma_i$) is a part of the
stability boundary prior to $A = A_{i,j}^t$. As $A$ increases from $A_{i,j}^t$,
$\Lambda_j$ (or $\Gamma_j$) becomes part of the boundary of the stability region,
separating segments of the bifurcation surface to which it was tangent.
However, many times as $A$ is increased $\Lambda_j$ (or $\Gamma_j$), the same surface (curve) which entered the boundary through tangency $A_{i,j}^t$, can be seen leaving the stability boundary via a \textbf{reverse tangency}, denoted $\tilde{A}_{j,i}^t$.
\end{definition}

\begin{figure}[htb]
\begin{center}
\begin{tabular}{ccc}
\includegraphics[width=0.3\textwidth]{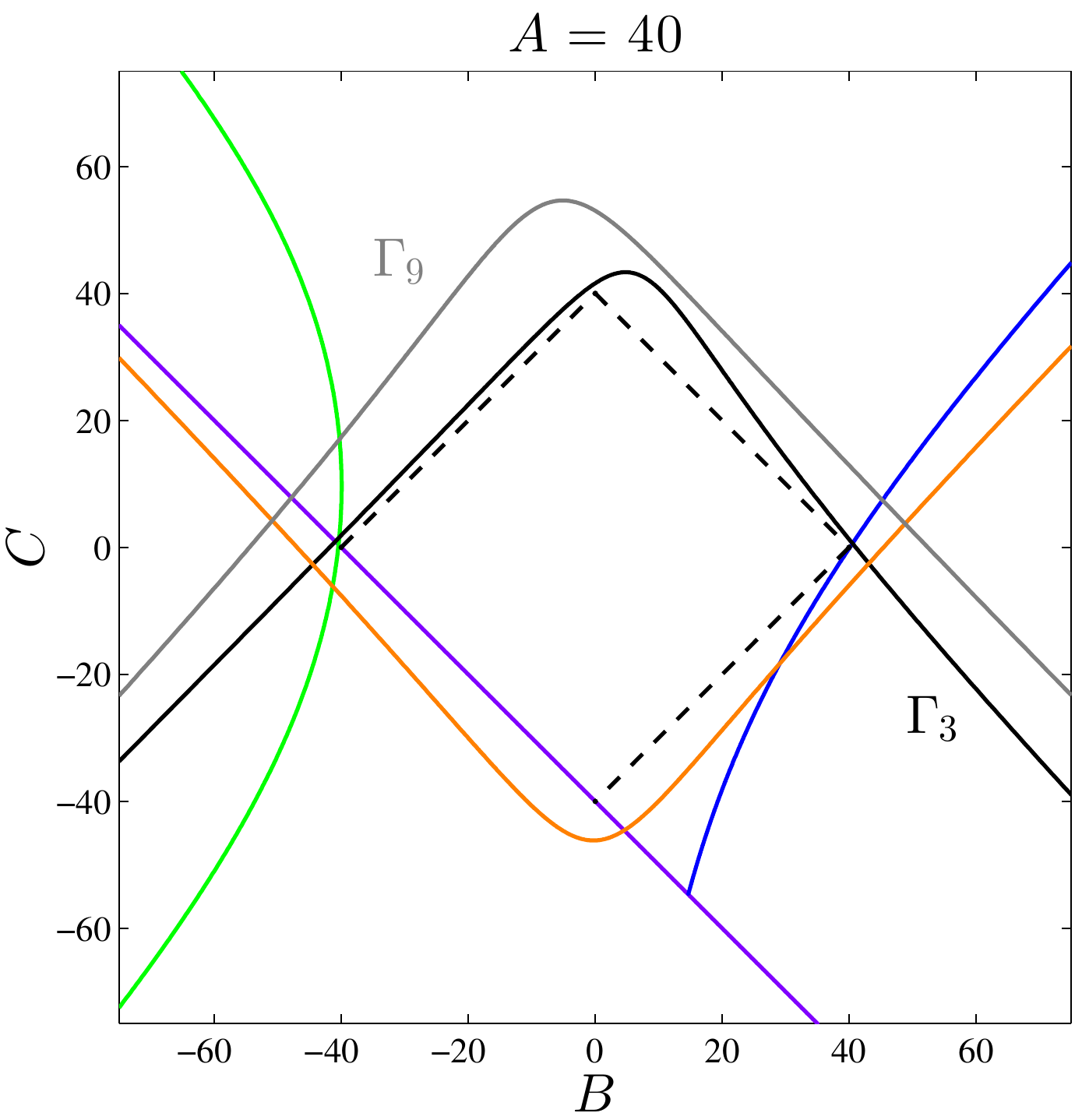}
&
\includegraphics[width=0.3\textwidth]{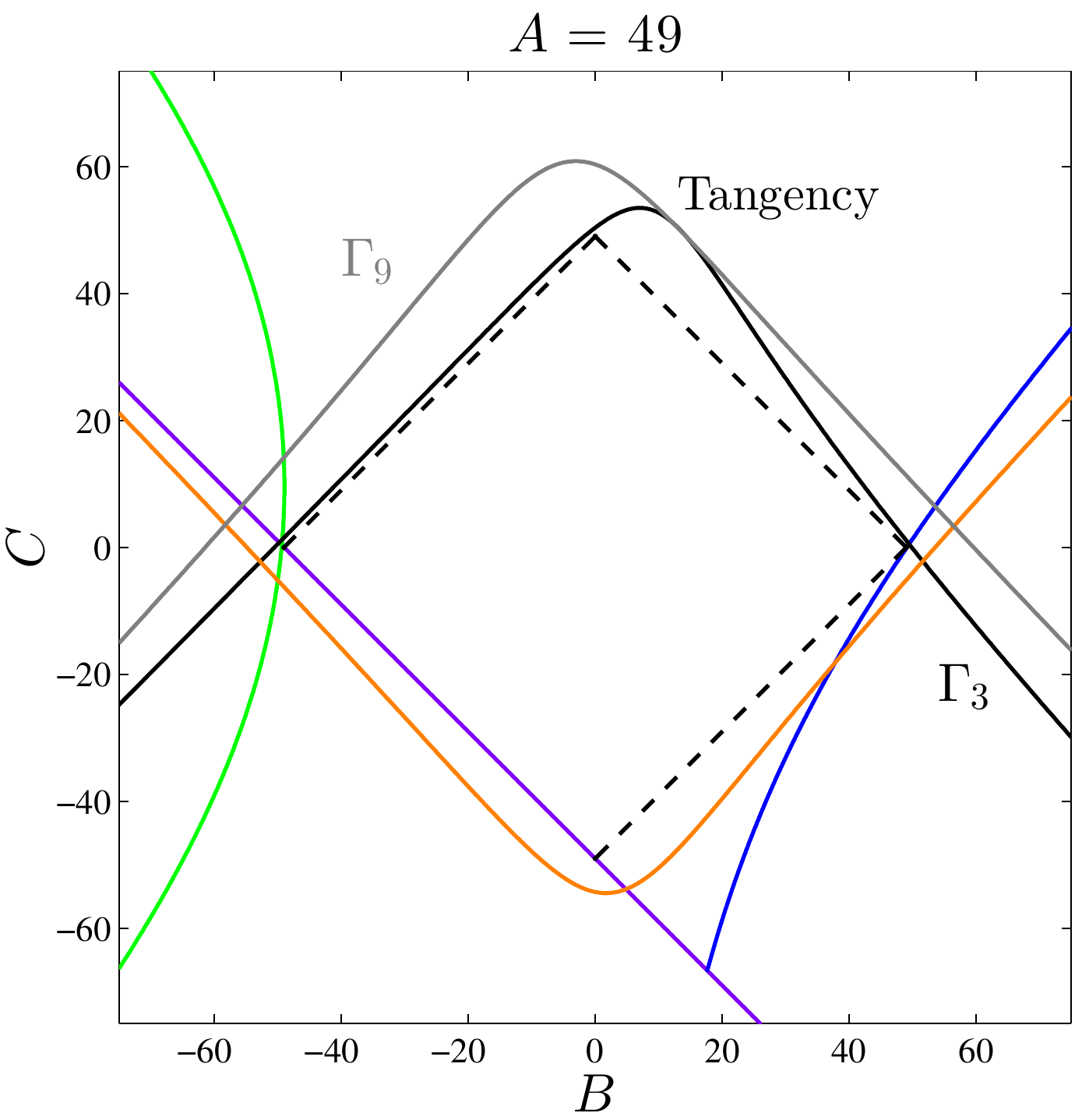}
&
\includegraphics[width=0.3\textwidth]{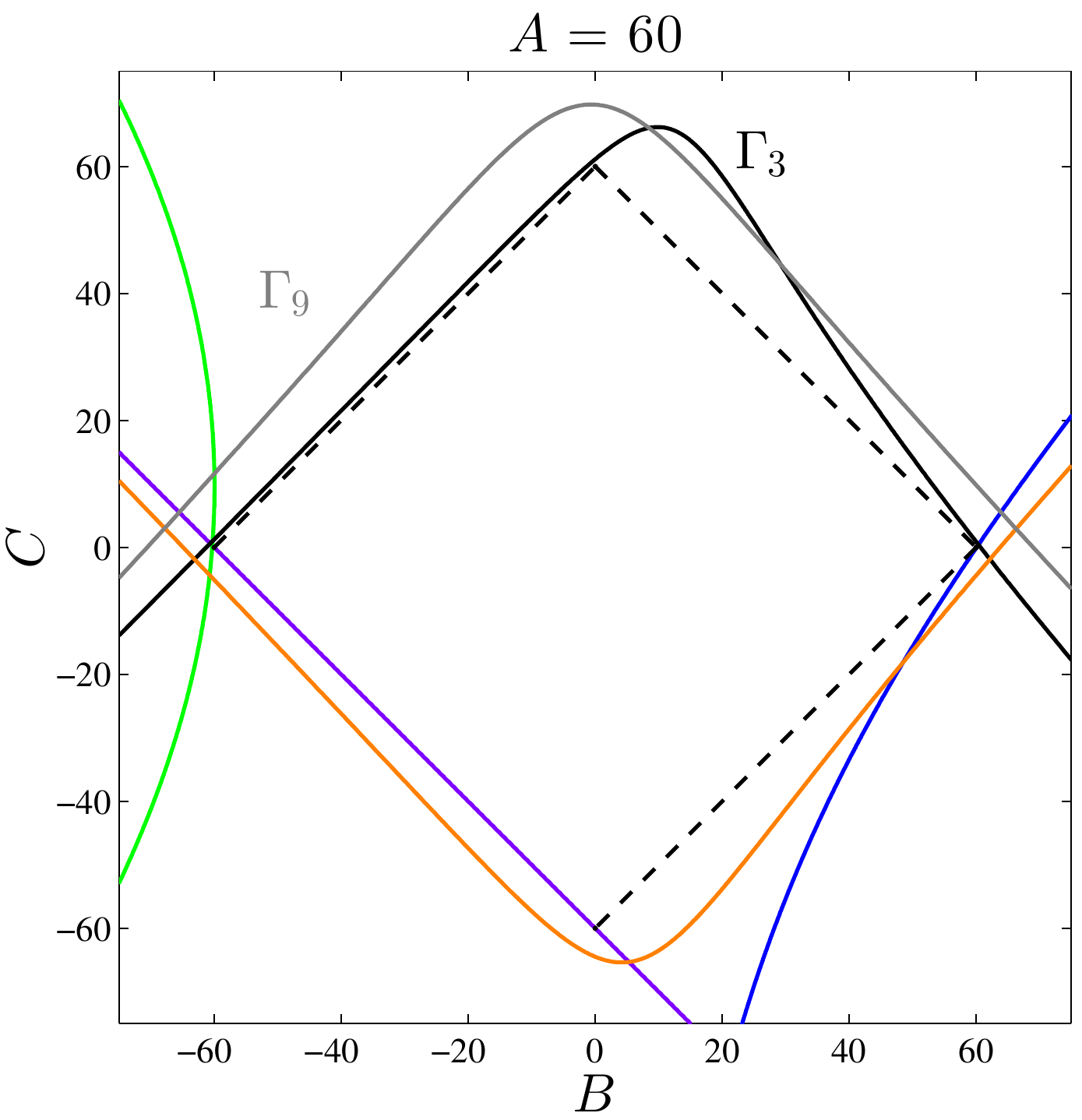}
\end{tabular}
\caption{Tangency, $A^t_{3,9}$. The  stability boundary in the $BC$ plane for $R=0.25$ is given before, during and after $A_{3,9}^{t}\approx 49$. \,The color scheme for the curves is: $\Lambda_0$ (violet), $\Gamma_1$ (blue), $\Gamma_2$ (green), $\Gamma_3$ (black), $\Gamma_6$ (orange), and $\Gamma_9$ (gray). The boundary is the closed region composed of portions of various competing curves that enshroud the MRS (bold dashed line).}
\label{fig:tang}
\end{center}
\end{figure}

Fig.~\ref{fig:transfer} shows an example of the transferral, $A^z_{1,6}$, where bifurcation curve, $\Gamma_6$, enters the stability surface for $A > A^z_{1,6}$ when $R = \frac{1}{4}$. We can readily see this change in the stability region near where $\Gamma_0$ and $\Gamma_1$ intersect. Fig.~\ref{fig:tang} shows an example of a tangency, $A^t_{3,9}$, where bifurcation curve, $\Gamma_9$, enters the stability surface for $A > A^t_{3,9}$ when $R = \frac{1}{4}$. In this case, $\Gamma_9$ becomes tangent to $\Gamma_3$ and for larger $A$ values becomes part of the stability boundary. We note that the majority of the changes to the stability surface come from tangencies (or reverse tangencies).

\setcounter{equation}{0}
\setcounter{theorem}{0}
\setcounter{figure}{0}
\setcounter{table}{0}
\section{Examples from Numerical Studies}

\begin{figure}[htb]
\begin{center}
\begin{tabular}{cc}
\includegraphics[width=0.4\textwidth,height=.4\linewidth]{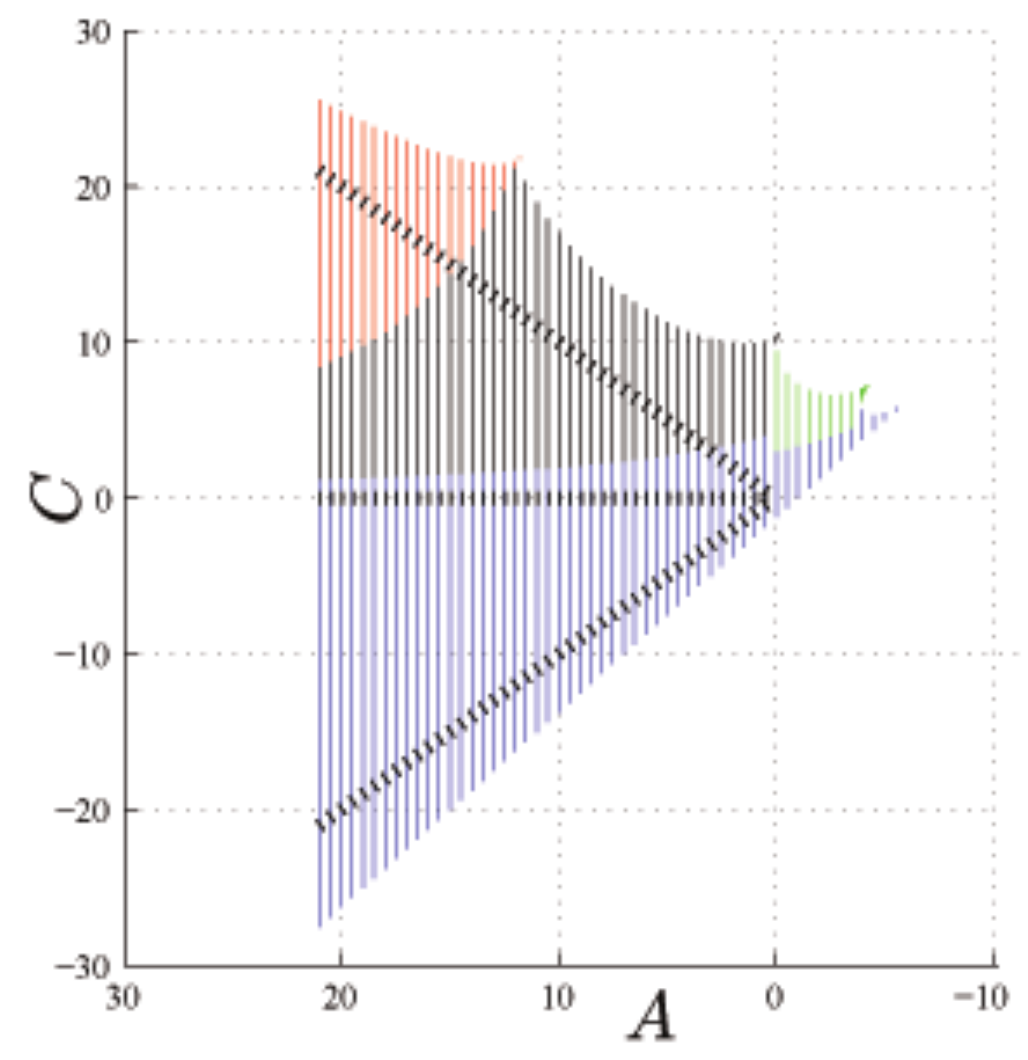}
&
\includegraphics[width=0.5\textwidth]{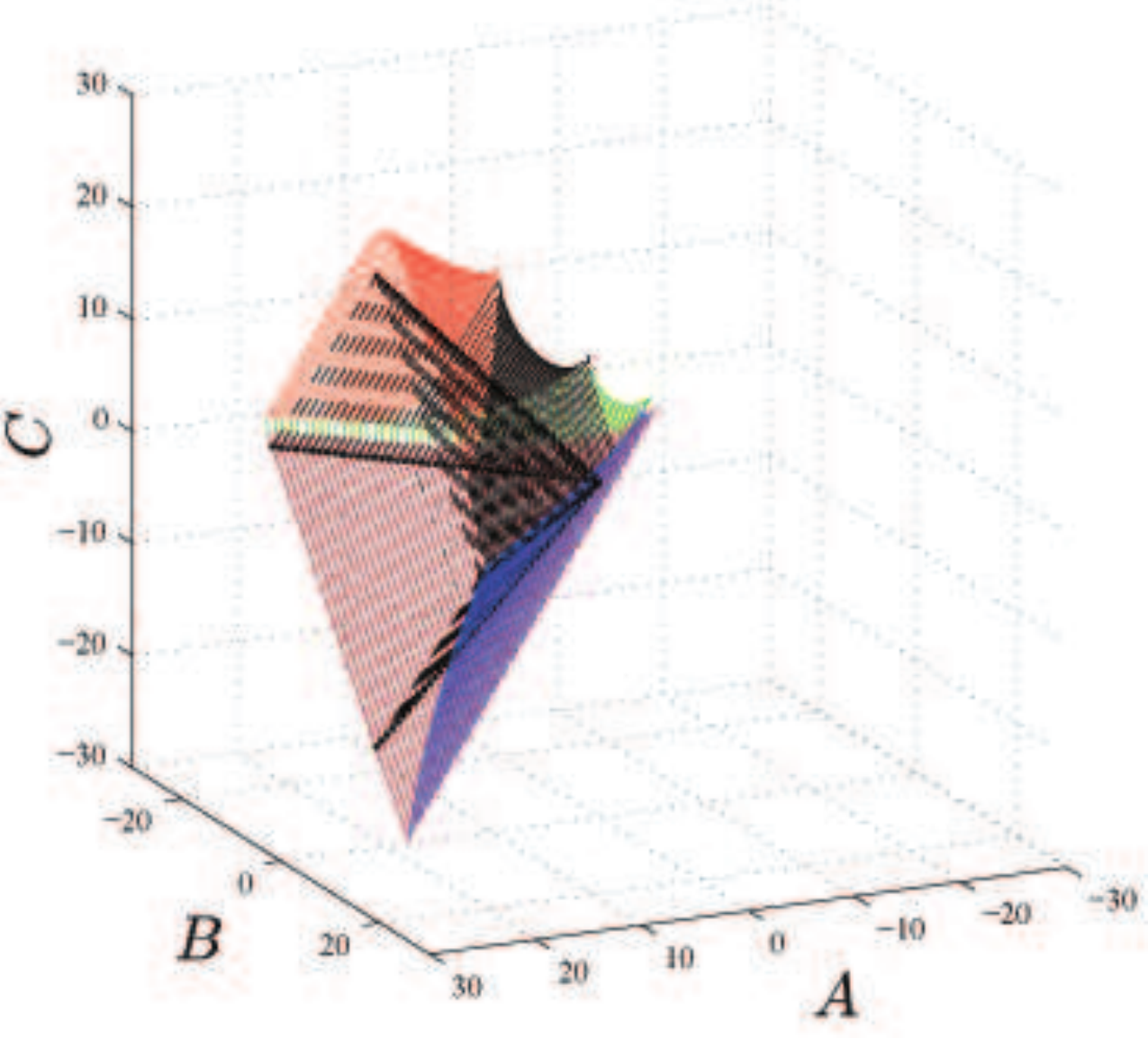}
\end{tabular}
\caption{Pictures of the portion of the stability surface comprised of
$A\in[-6,21]$ for $R = \frac{1}{5}$.}
\label{fig:Rfifth}
\end{center}
\end{figure}

We have developed a number of tools in MatLab to facilitate our stability studies
of (\ref{DDE2}). The ability to rapidly generate bifurcation curves has led to
significant insight into how the stability region evolves in $A$ and $R$ as viewed
in the $BC$-parameter plane. In this section we begin with some 3D stability surfaces
for $R = \frac{1}{5}$ to illustrate how the region of stability varies with $A$. We
present a diagram, which was numerically generated, to illustrate the systematic ordering of transitions, transferrals, and tangencies for a range of $R$ values. Finally, we end this
section with detailed numerical studies for specific values of $R$ as $R \to \frac{1}{4}$.

Fig.~\ref{fig:Rfifth} provides two views of the stability region of (\ref{DDE2}) with
$R = \frac{1}{5}$ and $A \le 21$. For $R = \frac{1}{5}$ there are only three
transitions with a $4^{th}$ transition occurring at $A = +\infty$.
Fig.~\ref{fig:Rfifth} shows the starting point of the stability surface at
$\left(-6, -\frac{1}{4}, \frac{25}{4}\right)$. Initially,  $\Lambda_1$ (blue) and
$\Lambda_0$ (violet) bound the stability region. Next the first stability
spur (green) enters and attaches $\Lambda_2$ to the stability
region at the first transition, $A_1^*$. Subsequently, $\Lambda_3$ (black) and
$\Lambda_4$ (red) adjoin the boundary of stability. At $A \approx 21$ a transferral occurs bringing $\Lambda_7$ onto the boundary, and around $A \approx 70$, there is a tangency of $\Lambda_{11}$ interrupting $\Lambda_4$. For this stability surface with no other transitions affecting the boundary, we only observe additional tangencies, followed by reverse tangencies, where bifurcation
surfaces leave the boundary, which change the boundary of the stability surface for
larger $A$. Fig.~\ref{fig:Rfifth} shows the MRS (black) interior to the stability surface,
and visually it is apparent how much the transitions and stability spurs distort the
stability surface away from the MRS. It is worth noting that the stability spurs
are shrinking in size as $A$ increases.

\begin{figure}[htb]
\begin{center}
\includegraphics[width=0.9\textwidth]{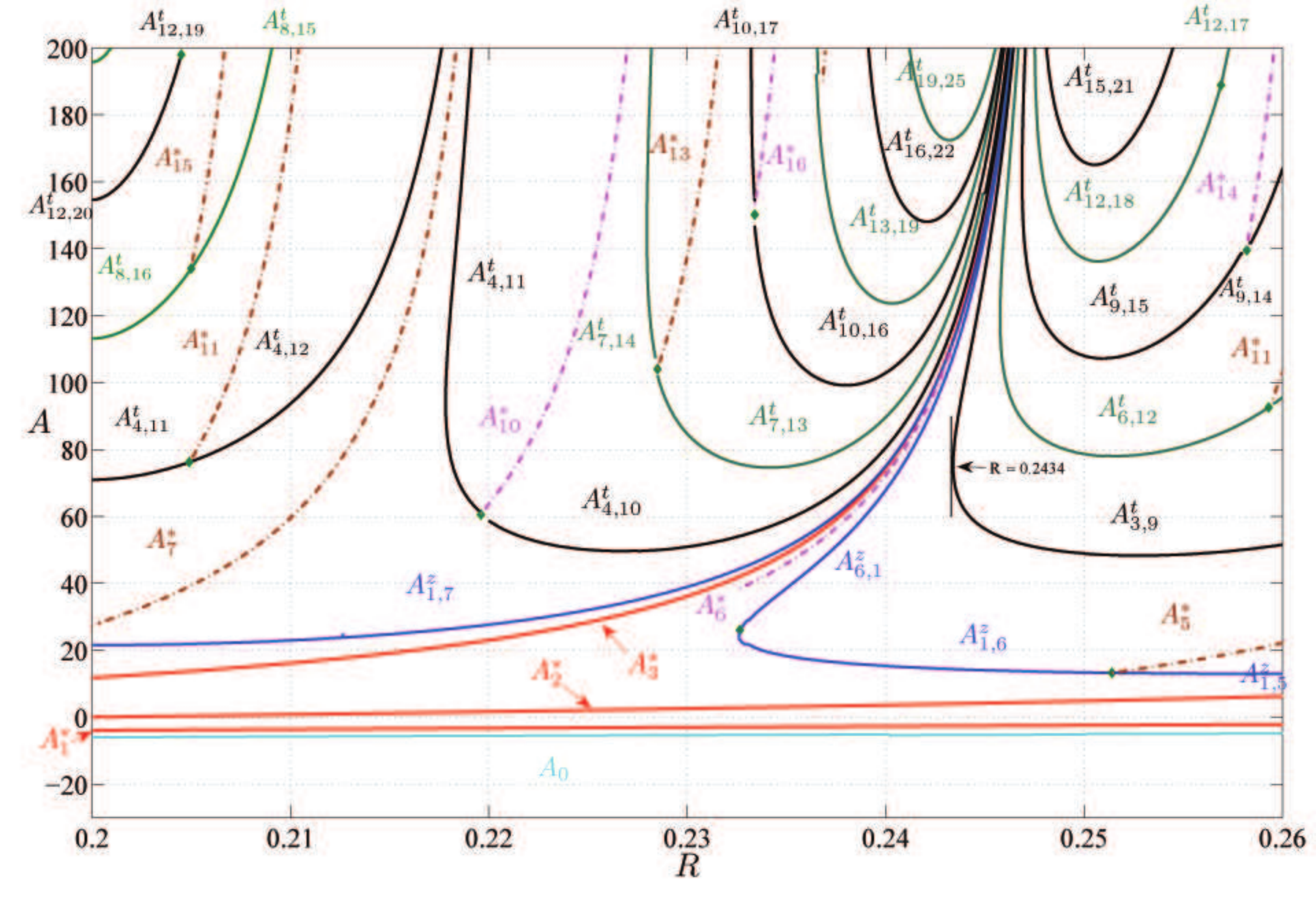}
\caption{The various transitions, transferrals, and tangencies are illustrated along
with $A_0$ for $R\in[0.20, 0.26]$ and $A\leq200$.}
\label{fig:Rall}
\end{center}
\end{figure}

Fig.~\ref{fig:Rall} has a diagram for the range of $R \in [0.2,0.26]$ showing all observed initial points $A_0$, transitions, transferrals, and tangencies. Following a vertical line, {\it i.e.}, fixing $R$, shows exactly which bifurcation curves enter or leave the boundary of the stability surface
as $A$ increases. The transition curves in the $RA$-plane increase monotonically. The $1^{st}$
transition is asymptotic with $A_1^* \to \infty$ at $R = \frac{1}{2}$, while
$A_2^* \to \infty$ at $R = \frac{1}{3}$. When $R = \frac{1}{4}$, $A_3^* \to
\infty$, with the other two transitions, their stability spurs, and one transferral
all occurring before $A = 15$. Afterwards, when $R = \frac{1}{4}$, all remaining changes
to the stability surface occur through tangencies, reverse tangencies, or a
reverse transferral. None of these events dramatically change the shape of the
boundary of the stability region. What is significant is that the $3^{rd}$
transition does not occur until $A_3^* = +\infty$, causing the distortion of this transition
to persist, while nearby delays do not have this distortion for sufficiently large
$A$.

\subsection{Stability region near $R = \frac{1}{4}$}

The characteristic equation (\ref{chareqn}) is an analytic function, so there is continuity of the stability surfaces as the parameters vary. To study what happens to the stability surface for $R = \frac{1}{4}$, we explore in some detail the stability surface for $R = 0.249$. Not surprisingly, there are many similarities between these surfaces until the singularity occurs at the transition, $A_3^* = 749.93$ for $R = 0.249$. We provide details of the evolution of stability surface for $R = 0.249$, which suggests why the region of stability for $R = \frac{1}{4}$ remains larger than the MRS as $A$ increases.

Fig.~\ref{fig:Rall} provides key information on how to determine which bifurcation surfaces compose the boundary of the stability region. Critical changes to the stability region are determined by the intersection of the various curves with a vertical line from a given $R$ as $A$ increases. For $R = 0.249$, this figure shows that the stability region begins at the starting point, $(A_0, B_0, C_0) \approx (-5.016 ,-0.3316, 5.348)$. A stability spur begins at $A_1^p \approx -2.733$ and joins the main stability surface at $A_1^* \approx -2.446$. A second stability spur, which is significantly smaller in length, joins the main stability surface near $A_2^* \approx 4.71$. Continuing vertically in Fig.~\ref{fig:Rall} at $R = 0.249$, we see there is a transferral, $A_{1,6}^z \approx 13.3$. At this stage, the $BC$ cross-section of the stability surface is very similar to the images in Fig.~\ref{fig:transfer}. The boundary at the transferral is comprised of $\Lambda_0$, $\Lambda_1$, $\Lambda_2$, and $\Lambda_3$. Subsequently, $\Lambda_6$ enters the boundary near the intersection of $\Lambda_0$ and $\Lambda_1$.

The next change in the stability surface for $R = 0.249$ is a tangency, which occurs at $A^t_{3,9} \approx 49.4$. This tangency has little effect on the shape of the boundary of the stability region, but allows higher frequency eigenvalues to participate in destabilizing (\ref{DDE2}). This tangency is very similar to the one depicted in Fig.~\ref{fig:tang}, occurring in the $1^{st}$ quadrant. As $A$ increases, a series of tangencies happens, alternating between the $1^{st}$ and $4^{th}$ quadrants. There are a total of 11 tangencies that occur before $A_3^* \approx 749.93$, with the last one being $A^t_{33,39} \approx 462.063$. Following $A^t_{33,39}$, all of these tangencies undergo a reverse tangency (in the reverse sequential order) with higher frequency eigenvalues leaving the boundary of the region of stability. Table~\ref{summarytable249} summarizes all of these events. The onset of reverse tangencies, which occur prior to $A_3^*$, create significant expansion of the stability region, primarily in the $1^{st}$ and $4^{th}$ quadrants of the $BC$-plane. At the same time the stability surface becomes much larger than the MRS. We note that very rapidly after $A_3^*$, there are many tangencies, which occur for increasing $A$ and reduce the size of the boundary of stability for $R = 0.249$. This results in the region of stability shrinking back to being near the MRS for large $A$. Since $R = 0.249$ is rational, we conjecture that the boundary of the stability region never asymptotically approaches the MRS, yet it is substantially closer than for $R = \frac{1}{4}$.

\begin{table}[htb]
 \hspace*{0.01in}
 \parbox{4.9in}{
   \caption{List of 2D boundary changes for $R=0.249$ and $A\in[A_0, 750]$.}
   }
    \begin{center}
   \renewcommand{\arraystretch}{1.25}
 \begin{tabular}{|l|l|l|l|}
   \hline
   surface change & $A$ &  surface change & $A$\\ \hline
   $A_0$ & $-5.016$ &  reverse tangency & $\tilde{A}^t_{39,33}\approx559.216$\\ \hline
   spur 1  & $[A^p_1, A^*_1] \approx [-2.7326, -2.4464]$ & reverse tangency & $\tilde{A}^t_{36,30}\approx622.341$ \\[0.1in] \hline 
    spur 2  & $[A^p_2, A^*_2] \approx [4.7067098, 4.70671]$&reverse tangency & $\tilde{A}^t_{33,27}\approx655.407$  \\[0.1in] \hline  
   transferral & $A^z_{1,6}\approx13.3$& reverse tangency & $\tilde{A}^t_{30,24}\approx678.811$\\[0.1in] \hline
 tangency & $A^t_{3,9}\approx49.4$& reverse tangency & $\tilde{A}^t_{27,21}\approx696.727    $\\[0.1in] \hline
tangency & $A^t_{6,12}\approx80.216$& reverse tangency & $\tilde{A}^t_{24,18}\approx 710.883  $\\[0.1in] \hline 
tangency & $A^t_{9,15}\approx108.4$& reverse tangency & $\tilde{A}^t_{21,15}\approx722.187$\\[0.1in] \hline
tangency & $A^t_{12,18}\approx142.479$& reverse tangency & $\tilde{A}^t_{18,12}\approx731.176$\\[0.1in] \hline 
tangency & $A^t_{15,21}\approx174.915$& reverse tangency & $\tilde{A}^t_{15,9}\approx738.192  $\\[0.1in] \hline
tangency & $A^t_{18,24}\approx 208.787 $& reverse tangency & $\tilde{A}^t_{12,6}\approx743.460 $\\[0.1in] \hline
tangency & $A^t_{21,27}\approx244.699 $& reverse  tangency & $\tilde{A}^t_{9,3}\approx747.134$\\[0.1in] \hline
tangency & $A^t_{24,30}\approx283.613$& reverse transferral & $\tilde{A}^z_{6,1}\approx749.4$\\[0.1in] \hline
tangency & $A^t_{27,33}\approx327.299$& spur~3  & $ A^*_3 \approx 749.93$\\[0.1in] \hline
tangency & $A^t_{30,36}\approx379.973$& transferral & $A^z_{1,7}\approx749.94$\\[0.1in] \hline
tangency & $A^t_{33,39}\approx462.063$&  tangency & $A^t_{4,10}\approx749.953 $\\[0.1in] \hline
&&   tangency &$A^t_{7,13}\approx750.044$\\[0.1in] \hline
 \end{tabular}
  \end{center}
  \label{summarytable249}
\end{table}

\begin{figure}[htb]
  \begin{center}
  \begin{tabular}{cc}
       \includegraphics[width=3.0in]{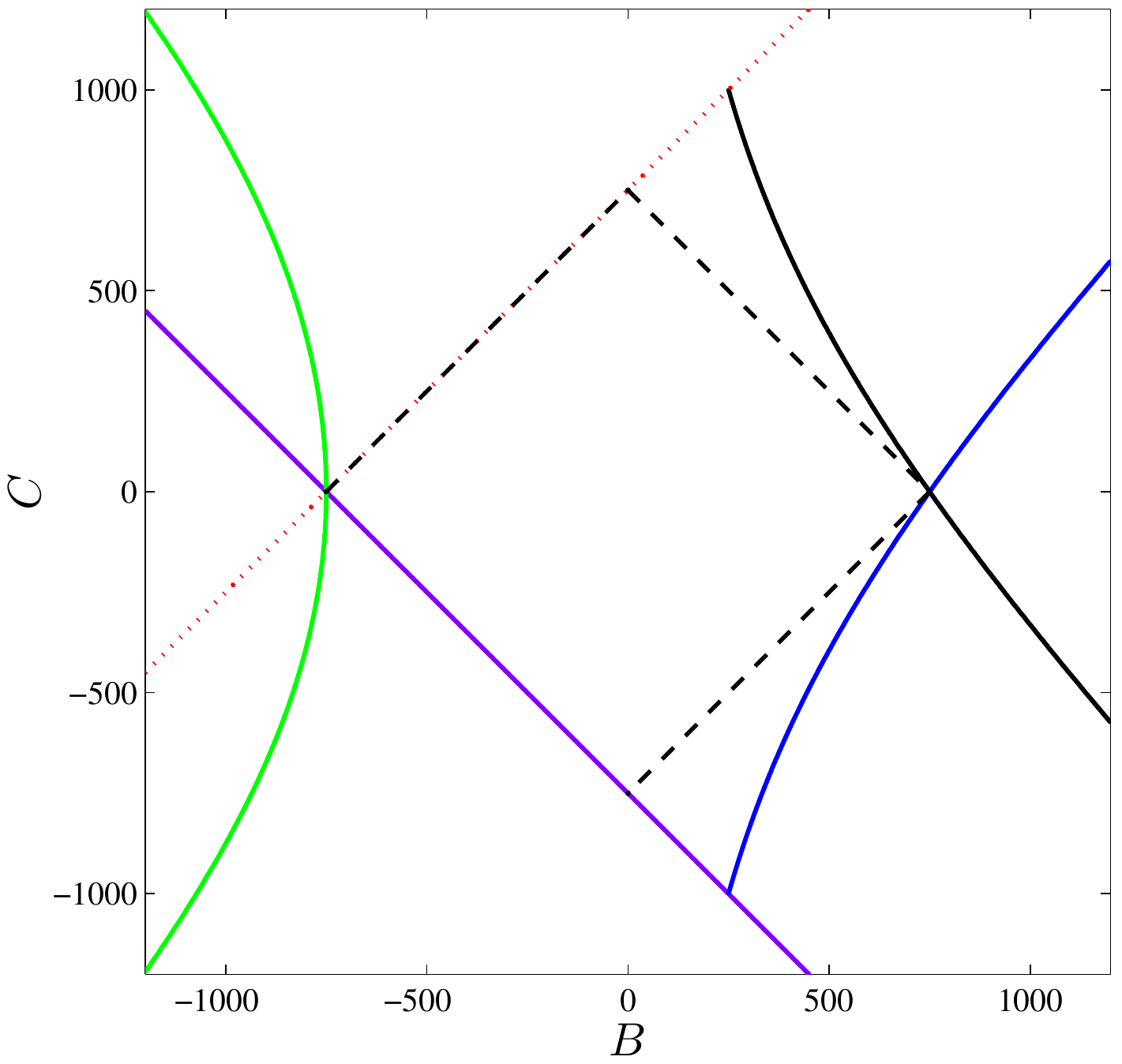} &
       \includegraphics[width=3.0in]{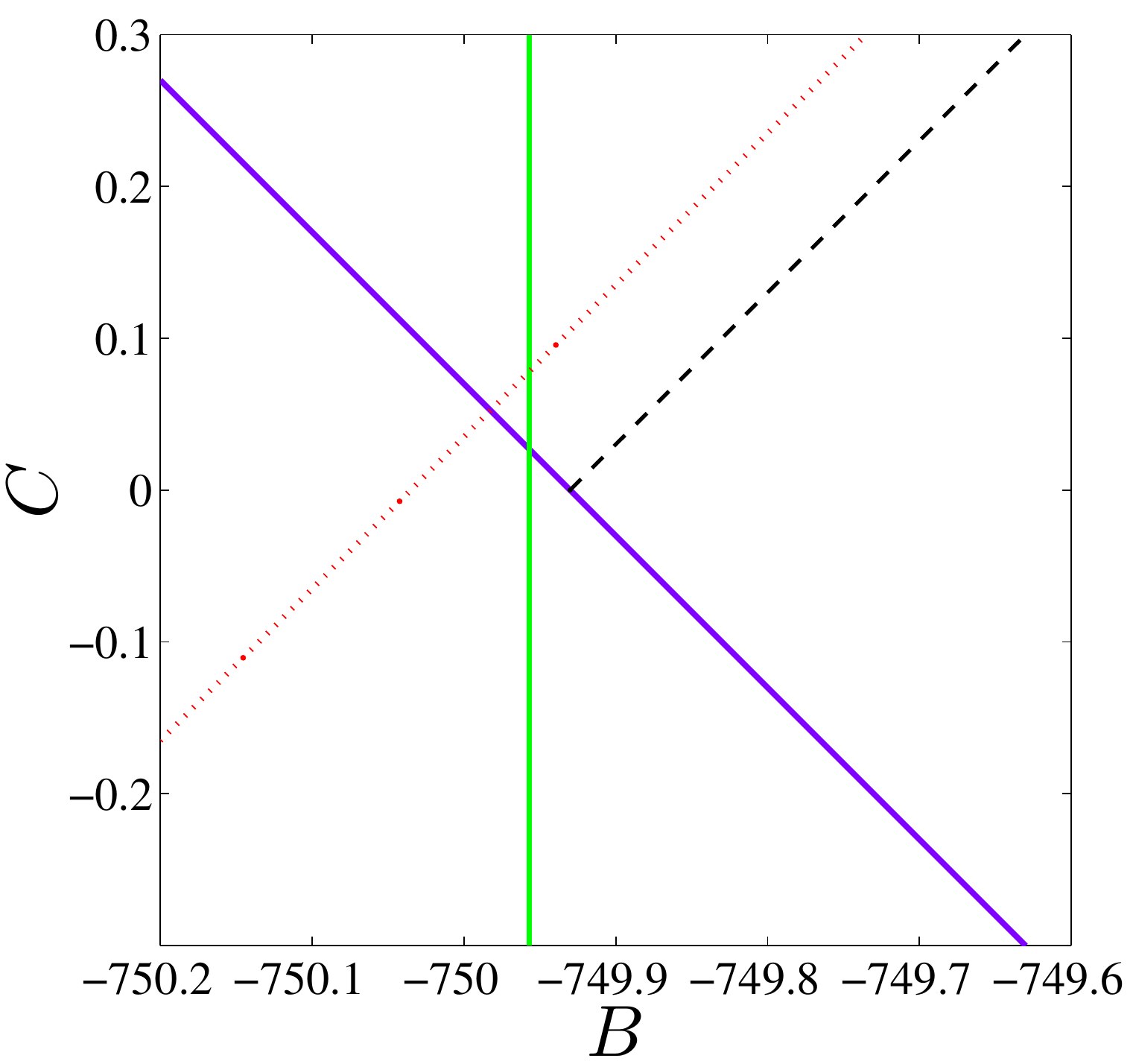}
  \end{tabular}
  \end{center}
       \caption{\small{The five curves on the boundary of the stability region for $R = 0.249$ at
       $A_3^* = 749.93$ are $\Lambda_0$ (violet), $\Gamma_1$ (blue), $\Gamma_3$ (black), $\Delta_3$ (dashed red) and minimally $\Gamma_2$ (green) with the MRS (dashed black).}}
    \label{fig:R249_tran}
\end{figure}

For $R = 0.249$ at $A_3^* = 749.93$, the boundary of the region of stability is reduced to only 5 bifurcation curves. There are lines from $\Lambda_0$ (violet) and the degeneracy line, $\Delta_3$,
with $\lambda = \pm\frac{3i\pi}{0.751}$ at $A_3^*$. The boundaries in the $1^{st}$ and $4^{th}$ quadrants are formed from $\Gamma_3$ (black) and $\Gamma_1$ (blue), respectively. Finally, there is a very small segment of the boundary formed by the $\Gamma_2$ (green). This stability region is shown in Fig.~\ref{fig:R249_tran}. We see distinct bulges away from the MRS in the $1^{st}$ and $4^{th}$ quadrants, increasing the size of the region of stability. This simple
boundary easily allows the computation of the area of the region of stability.
A numerical integration gives that the region of stability outside the MRS is
approximately 26.85\% of the area of the MRS, which is a substantial increase
in the region of stability.

\setcounter{equation}{0}
\setcounter{theorem}{0}
\setcounter{figure}{0}
\setcounter{table}{0}
\section{Analysis}

The objective of this section is to convince the reader that the asymptotic stability region for $R = \frac{1}{n}$ with $A \to +\infty$ reduces to a simple set of curves bounded away from the MRS. As $R \to \frac{1}{n}$, there is a critical transition with $A^*_{n-1} \to +\infty$, which leaves the boundary of the stability region for (\ref{DDE2}) composed of only $\Lambda_0$, $\Gamma_1$, $\Gamma_{n-1}$, and the degeneracy line, $\Delta_{n-1}$. (There may be small segments of other bifurcation curves for any $R < \frac{1}{n}$ at $A^*_{n-1}$.) Mahaffy {\it et al.} \cite{MZJa, MZJ} showed that as $\omega \to 0^+$, $\Gamma_1$ intersects $\Lambda_0$ along the line
\begin{equation}
  \frac{A + 1}{1 - R} = \frac{B - 1}{R} = -C.  \label{plane_bif1}
\end{equation}
For $R = \frac{1}{n}$, it follows that $B = \frac{A + n}{n-1}$ and $C = -\frac{n(A +1)}{n-1}$. Thus, the distance from the MRS to $\Gamma_1$ along $\Lambda_0$ asymptotically (large $A$) extends past the MRS by a length that is a factor of $\frac{1}{n - 1}$ longer than the length of the edge of the MRS. This will provide one measure of the extension of the stability region for (\ref{DDE2}).

The cases for $R = \frac{1}{2n}$ and $\frac{1}{2n-1}$ give different shaped regions, but ultimately as $A \to +\infty$, the regions are bounded by only four curves with two lying on the boundary of the MRS and two bulging away from this region. We prove analytically the position of some key points, then rely on the limited number of families of curves and their distinct ordering to give the simple structure. The observed orderly appearance and disappearance of the tangencies on the boundary of the region of stability provide our argument for the ultimate simple structure of the stability region and the continued bulge away from the MRS for these particular rational delays. From the monotonicity of the transition curves, we note that all limiting arguments require $R$ approaching $\frac{1}{n}$ from below.

\subsection{Stability Region for $R = \frac{1}{2n}$}

In this section we  provide more details on the increased size of the region of stability for delays of the form $R = \frac{1}{2n}$ as $A \to +\infty$. Earlier we presented evidence that as $R \to \frac{1}{4}$, the third transition, $A_3^*$, tended to infinity, and at that transition the boundary of the stability region in the $BC$-plane reduced to just four curves. $\Lambda_0$ provides the lower left boundary along the MRS in the $3^{rd}$ quadrant. With $R < \frac{1}{4}$, $\Delta_3$, which occurs at $A_3^*$ with eigenvalues $\lambda = \frac{3\pi}{1-R}$, creates a boundary parallel to the MRS in the $2^{nd}$ quadrant. This line approaches the MRS as $R \to \frac{1}{4}$ from below. $\Gamma_1$ forms a boundary in the $4^{th}$ quadrant, which significantly bulges away from the MRS. Its intersection at $\Lambda_0$ extends $\frac{1}{3}$ times the length of the MRS along the line of this plane into the $4^{th}$~quadrant. $\Gamma_3$ creates an almost mirror image across the $C$-axis in the $1^{st}$ quadrant, bulging away from the MRS the same distance. Fig.~\ref{fig:R249_tran} illustrates this expanded stability region very well. Below we prove some results about the four curves on the boundary of the stability region and give additional information on why we believe the stability region has its enlarged character.

The case $R = \frac{1}{4}$ extends generically to the case $R = \frac{1}{2n}$. $\Lambda_0$ is always one part of the stability boundary. Symmetric to this boundary across the $B$-axis is $\Delta_{2n-1}$ at $A_{2n-1}^*$, which occurs with $A_{2n-1}^* \to \infty$ as $R \to \frac{1}{2n}$. Below we demonstrate that $\Delta_{2n-1}$ approaches the line $C - B = A_{2n-1}^*$ in the $BC$-plane as $A_{2n-1}^* \to \infty$, which is one side of the MRS. The other two sides are composed of $\Gamma_1$ in the $4^{th}$ quadrant and symmetric to this, $\Gamma_{2n-1}$ in the $1^{st}$ quadrant.

To help obtain the symmetrical shape discussed above (and exclude the bifurcation
curves that pass near the point $(B, C) = (-A_{2n-1}^*,0)$), there are alternate
forms of the Eqns.~(\ref{bctrans}). We multiply and divide the cosecant terms by
$\cos\left(\frac{jR\pi}{1-R}\right)$, then use the definition of $A_j^*$ to obtain
\begin{eqnarray}
 B_j^*(R) & = & (-1)^j\left(\frac{1}{1-R}\cos\left(\textstyle{\frac{jR\pi}{1-R}}\right)
             + \frac{RA_j^*}{(1-R)\cos\left(\textstyle{\frac{jR\pi}{1-R}}\right)}\right)
              \nonumber \\
 C_j^*(R) & = & \frac{-A_j^*}{(1-R)\cos\left(\textstyle{\frac{jR\pi}{1-R}}\right)}
             - \frac{1}{1-R}\cos\left(\textstyle{\frac{jR\pi}{1-R}}\right) \label{bc_alt}
\end{eqnarray}

\begin{lemma}\label{lem_lim_line_even} For $R = \frac{1}{2n}$, one boundary of the region of stability
is the limiting line
\[
 C - B = A_{2n-1}^*,
\]
which lies on the MRS.
\end{lemma}

\proof
For $R < \frac{1}{2n}$ and $R \to \frac{1}{2n}$, we consider the transition $A_{2n-1}^*$. The degeneracy line, $\Delta_{2n-1}$, satisfies
\[
  A = A_{2n-1}^*, \quad C - B = C_{2n-1}^* - B_{2n-1}^*.
\]
Since $j = 2n-1$ and $A_{2n-1}^* = -\frac{(2n-1)\pi}{1-R} \cot\left(\frac{(2n-1)R\pi}{1-R}\right)$, Eqns.~(\ref{bc_alt}) give
\[
 C_{2n-1}^* - B_{2n-1}^* = \textstyle{\frac{(2n-1)\pi}{1-R}\csc\left(\frac{(2n-1)R\pi}{1-R}\right)}.
\]
Now consider
\[
 \lim_{R \to \frac{1}{2n}^-} \frac{C_{2n-1}^* - B_{2n-1}^*}{A_{2n-1}^*}
   = \lim_{R \to \frac{1}{2n}^-}\textstyle{-\sec\left(\frac{(2n-1)R\pi}{1-R}\right)} = 1.
\]
Thus, $C_{2n-1}^* - B_{2n-1}^* \to A_{2n-1}^*$ for $R < \frac{1}{2n}$ as $R \to \frac{1}{2n}$, so $\Delta_{2n-1}$ tends towards one edge of the MRS. We note that the small distance between $\Delta_{2n-1}$ and the MRS may allow a very small segment of $\Gamma_2$ to be part of the stability region for all $R < \frac{1}{2n}$.
\qed

For $R = \frac{1}{2n}$, we showed that $\Gamma_1$ intersects $\Lambda_0$ at $(B,C) = \left(\frac{A+2n}{2n-1}, -\frac{2n(A+1)}{2n-1}\right)$. We now show that as $R \to \frac{1}{2n}$ from below and $A_{2n-1}^* \to +\infty$, the point $(B_{2n-1}^*, C_{2n-1}^*)$ tends to the ordered pair that has $B$-axis symmetry to the intersection of $\Gamma_1$ and $\Lambda_0$.

\begin{lemma}\label{lem_even_pt}
For $R < \frac{1}{2n}$ and $R \to \frac{1}{2n}$, the bifurcation curve $\Gamma_{2n-1}$ comes to the point
\[
  (B_{2n-1}^*, C_{2n-1}^*) = \left(\frac{A_{2n-1}^*+2n}{2n-1}, \frac{2n(A_{2n-1}^*+1)}{2n-1}\right)
\]
with $A_{2n-1}^* \to +\infty$.
\end{lemma}

\proof
For $R \to \frac{1}{2n}$ from below with $j = 2n - 1$, Eqns.~\ref{bc_alt} yield
\begin{eqnarray*}
  \lim_{R \to \frac{1}{2n}} B_{2n-1}^*(R) & = &
   \lim_{R \to \frac{1}{2n}}\left(-\frac{\cos\left(\frac{(2n-1)R\pi}{1-R}\right)}{1-R}
    - \frac{RA_{2n-1}^*}{(1-R)\cos\left(\frac{(2n-1)R\pi}{1-R}\right)}\right) \\
    & = & \frac{2n}{2n-1} + \frac{A_{2n-1}^*}{2n-1},
\end{eqnarray*}
since $1 - R \to \frac{2n-1}{2n}$ and $\cos\left(\frac{(2n-1)R\pi}{1-R}\right) \to \cos(\pi) = -1$. Similarly,
\begin{eqnarray*}
  \lim_{R \to \frac{1}{2n}^-} C_{2n-1}^*(R) & = &
   \lim_{R \to \frac{1}{2n}^-}\left(- \frac{A_{2n-1}^*} {(1-R)\cos\left(\frac{(2n-1)R\pi}{1-R}\right)} - \frac{\cos\left(\frac{(2n-1)R\pi}{1-R}\right)}{1-R}
    \right) \\
    & = & \frac{2n\,A_{2n-1}^*}{2n-1} + \frac{2n}{2n-1}.
\end{eqnarray*}
({\bf Note:} The limit as $R \to \frac{1}{2n}^+$, we have $A^*_{2n-1} \to -\infty$, which is not of interest to our study.)
\qed

The next step in our analysis is to show that the bifurcation curves, $\Gamma_1$ and $\Gamma_{2n-1}$, pass arbitrarily close to the point $(A_{2n-1}^*,0)$ in the $BC$-plane as $R \to \frac{1}{2n}$. Note that this is the point on the MRS, which is opposite the point of intersection of $\Lambda_0$ and $\Delta_{2n-1}$ at $A^*_{2n-1}$.

\begin{lemma} \label{C_cross_even}
For $R < \frac{1}{2n}$ and $R \to \frac{1}{2n}$, the bifurcation curves, $\Gamma_1$ and $\Gamma_{2n-1}$, pass arbitrarily close to the point $(A_{2n-1}^*,0)$ in the $BC$-plane with $A_{2n-1}^* \to +\infty$.
\end{lemma}

\proof
The bifurcation curves cross the $B$-axis whenever $C(w) = 0$. From (\ref{bifC}), $C(w) = 0$ implies
\begin{equation}\label{Czero}
  \omega = -\frac{A \sin(\omega)}{\cos(\omega)} \qquad {\rm or} \qquad
  -\omega\cos(\omega) = A\sin(\omega).
\end{equation}
With this information it follows from (\ref{bifB}) that
\[
 B(\omega) = -\frac{A}{\cos(\omega)}.
\]
For $\Gamma_1$, $0 < \omega < \frac{\pi}{1-R}$, which tends to the interval $0 < \omega < \frac{2n\pi}{2n-1}$ as $R \to \frac{1}{2n}$. When $C(\omega) = 0$, the expression $-\omega\cos(\omega)$ is bounded near $\pi$ as $\omega \to \pi$. For $C(\omega) = 0$ as $A$ becomes arbitrarily large, then $\sin(\omega) \to 0$ or $\omega \to \pi^-$. Thus, as $A \to \infty$, $\omega \to \pi^-$ for $\Gamma_1$ to cross the $B$-axis, and
\[
 \lim_{\omega \to \pi^-} B(\omega) = A.
\]
It follows that as $R \to \frac{1}{2n}$ from below, then $A_{2n-1}^* \to \infty$ and $\Gamma_1$ passes arbitrarily close to the point $(A_{2n-1}^*,0)$.

From the definition of $\Gamma_{2n-1}$, $\frac{(2n-2)\pi}{1-R} < \omega < \frac{(2n-1)\pi}{1-R}$, which tends to the interval $\frac{2n(2n-2)\pi}{2n-1} < \omega < 2n\pi$ as $R \to \frac{1}{2n}$. It is easy to see that $\omega = (2n-1)\pi$ is inside this interval for $\Gamma_{2n-1}$. This $\omega$ being an odd multiple of $\pi$ and $n$ being fixed and finite, the same arguments above for $\Gamma_1$ hold, which implies
\[
 \lim_{\omega \to (2n-1)\pi^-} B(\omega) = A.
\]
It follows that as $R \to \frac{1}{2n}$ from below, then $A_{2n-1}^* \to \infty$ and $\Gamma_{2n-1}$ passes arbitrarily close to the point $(A_{2n-1}^*,0)$.
\qed

The lemmas above prove some of the features of the stability region in Fig.~\ref{fig:R249_tran} for $R = \frac{1}{2n}$ with $A \to +\infty$. In particular, we see the left boundaries aligning with the MRS in the $2^{nd}$ and $3^{rd}$ quadrants. We also proved that $\Gamma_1$ and $\Gamma_{2n-1}$ intersect near $B = A^*_{2n-1}$ when $C = 0$. Finally, we proved that $\Gamma_1$ and $\Gamma_{2n-1}$ intersect $\Lambda_0$ and $\Delta_{2n-1}$, respectively, in a symmetric manner at
\[
  (B, C) = \left(\frac{A_{2n-1}^*+2n}{2n-1}, \pm\frac{2n(A_{2n-1}^*+1)}{2n-1}\right),
\]
as $R \to \frac{1}{2n}$ from below and $A^*_{2n-1} \to +\infty$.

It remains to show that $\Gamma_1$ and $\Gamma_{2n-1}$ are the only bifurcation curves on the boundary as $R \to \frac{1}{2n}$ from below and $A \to +\infty$. In this work we will only present numerical evidence for this result, providing some ideas for how a rigorous proof might proceed. The numerical results will include how much larger the asymptotic region of stability is for rational delays of the form $R = \frac{1}{2n}$.

As noted earlier, we developed a MatLab code for easily generating and observing bifurcation curves at various values of $A$ and $R$ in the $BC$-plane. Mahaffy {\it et al.} \cite{MZJ} showed that the rational delays, $R$, due to the periodic nature of the sinusoidal functions, result in the bifurcation curves ordering themselves into \textit{\textbf{families of curves}}.

\begin{definition}\label{family_defn}
For $A$ \,fixed, take $R=\frac{k}{n}$ and $j=n-k$. From Eqns.~(\ref{bifB}) and (\ref{bifC}), one can see that the singularities occur at $\frac{ni\pi}{j}, i=0,1,...$ . The bifurcation curve $i$, $\Gamma_i$,  with $\frac{n(i-1)\pi}{j}<\omega<\frac{ni\pi}{j}$ satisfies: \vspace*{-0.15in}
\begin{equation*}
    B_{i}(\omega)=\frac{A\sin(\frac{k\omega}{n})+\omega\cos(\frac{k\omega}{n})}{\sin(\frac{j\omega}{n})}, \hspace{0.3in}
      C_{i}(\omega)=-\frac{A\sin(\omega)+\omega\cos(\omega)}{\sin(\frac{j\omega}{n})}
\end{equation*}
\noindent Now consider $\Gamma_{i+2j}$ with $\mu=\omega+2n\pi$, then
\begin{align*}
     B_{i+2j}(\mu)&=\frac{A\sin(\frac{k\mu}{n})+\mu\cos(\frac{k\mu}{n})}{\sin(\frac{j\mu}{n})}=
     \frac{A\sin(\frac{k\omega}{n})+(\omega+2n\pi)\cos(\frac{k\omega}{n})}{\sin(\frac{j\omega}{n})} \\[0.2in]
     C_{i+2j}(\mu)&= -\frac{A\sin(\omega)+(\omega+2n\pi)\cos(\omega)}{\sin(\frac{j\omega}{n})}
\end{align*}
These equations show that $B_{i+2j}(\mu)$ follows the same trajectory as $B_{i}(\omega)$ with a shift of $2n\pi\cos(\frac{k\omega}{n})/\sin(\frac{j\omega}{n})$  \,for $\omega\in\Big( \frac{(j-1)\pi}{1-R}, \frac{j\pi}{1-R}\Big)$, while  $C_{i+2j}(\mu)$ \,follows the same trajectory as $C_{i}(\omega)$ with a shift of \,$2n\pi\cos(\omega)/\sin(\frac{j\omega}{n})$  \,over the same values of $\omega$. \ This related behavior of bifurcation curves separated by $\omega=2n\pi$ creates $2j$ \textbf{families of curves} in the $BC$ plane for fixed $A$. \,Thus, there is a quasi-periodicity among the bifurcation curves when $R$ is rational.
\end{definition}

The organization of these families of curves and systematic transitions allow one to observe how the different bifurcation curves enter and leave the boundary of the stability region. (See Fig.~\ref{fig:Rall}.) We provide more details following some analytic results for delays of the form $R = \frac{1}{2n+1}$.

\subsection{Stability region for $R = \frac{1}{2n+1}$}

\begin{figure}[htb]
  \begin{minipage}[b]{0.55\linewidth}
    \begin{center}
      \includegraphics[width=\textwidth]{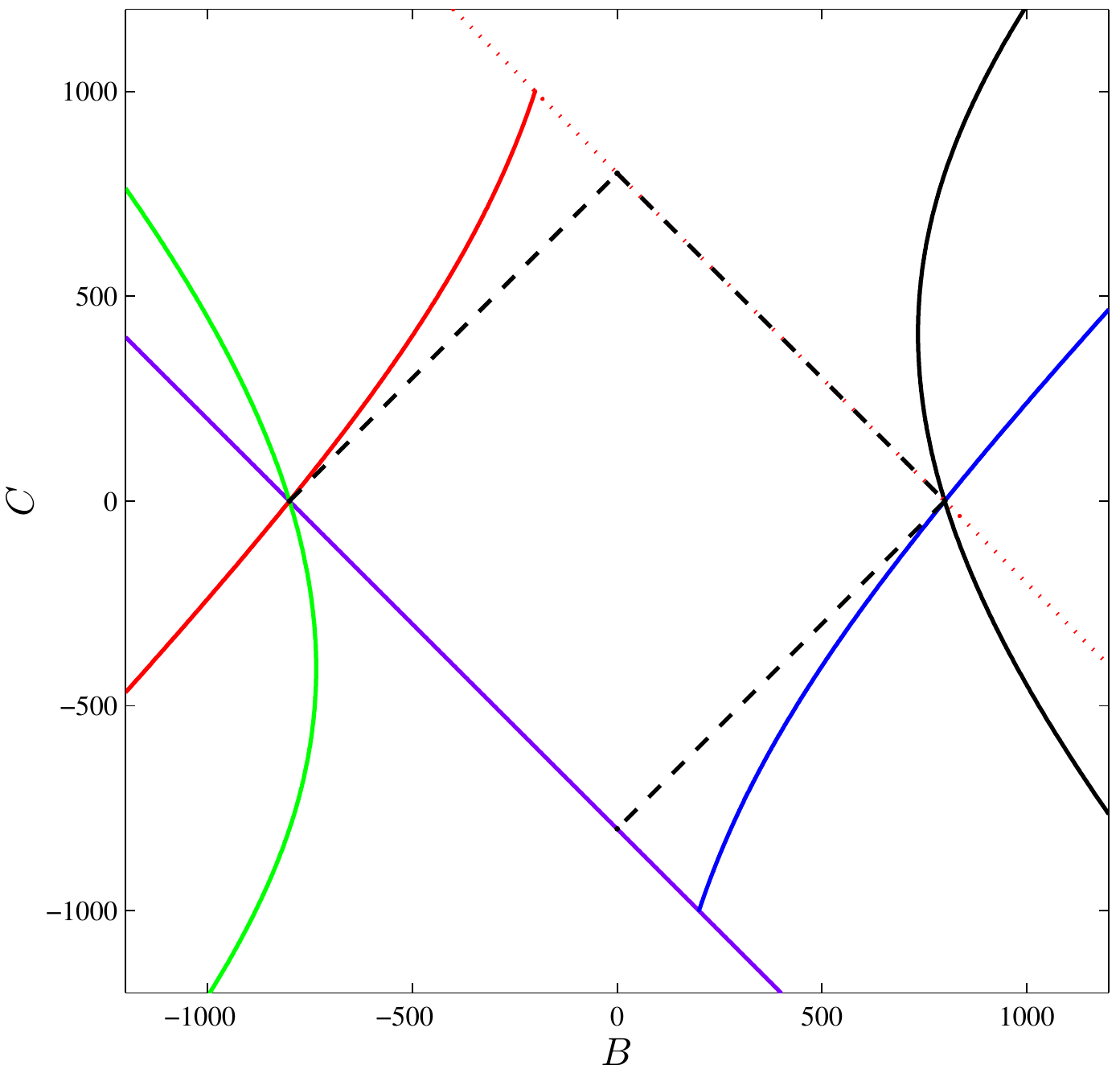}
    \end{center}
  \end{minipage}
  \hspace{0.5cm}
  \begin{minipage}[b]{0.4\linewidth}
    \begin{center}
      \includegraphics[width=0.7\textwidth]{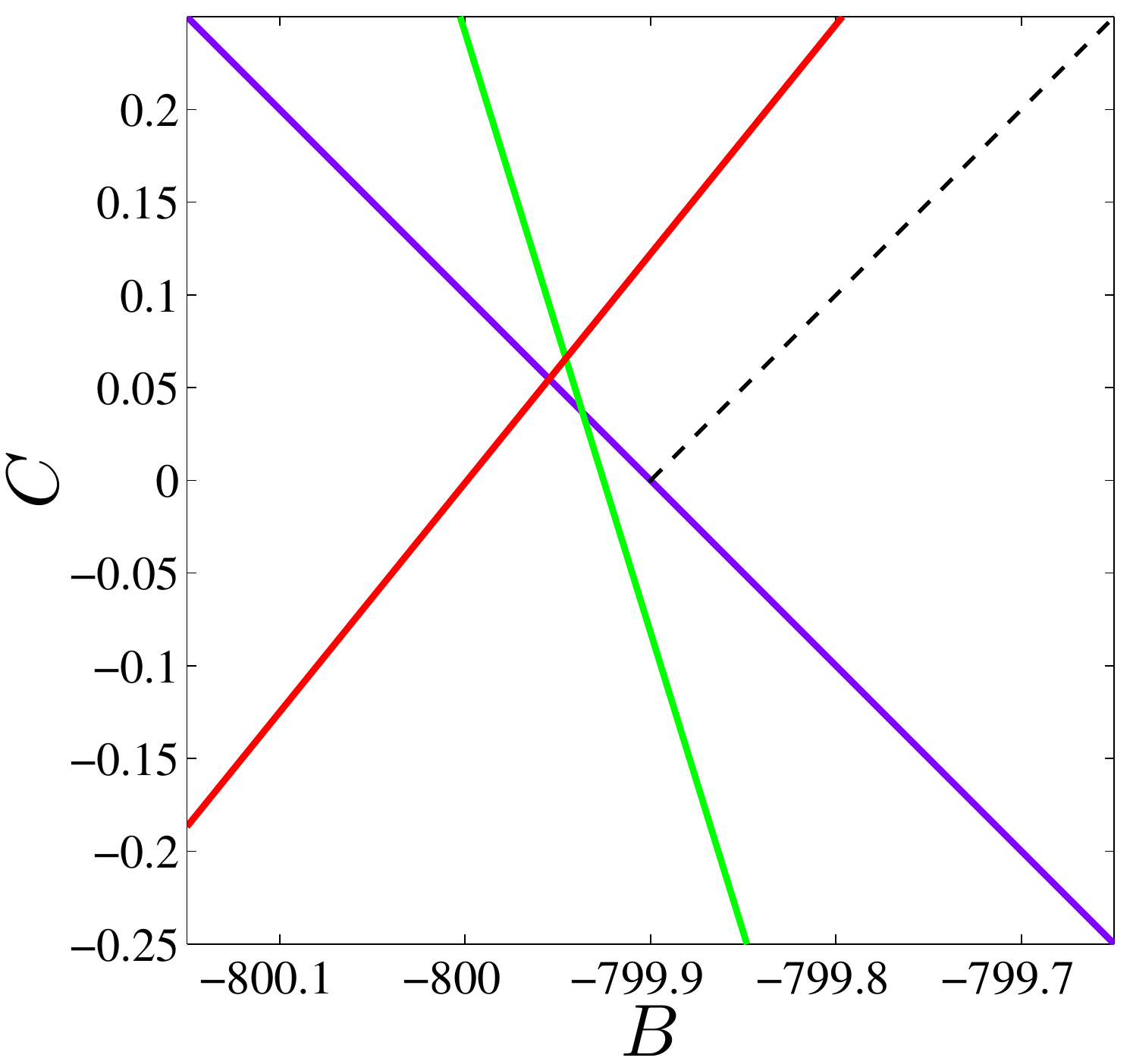} \\
      \includegraphics[width=0.7\textwidth]{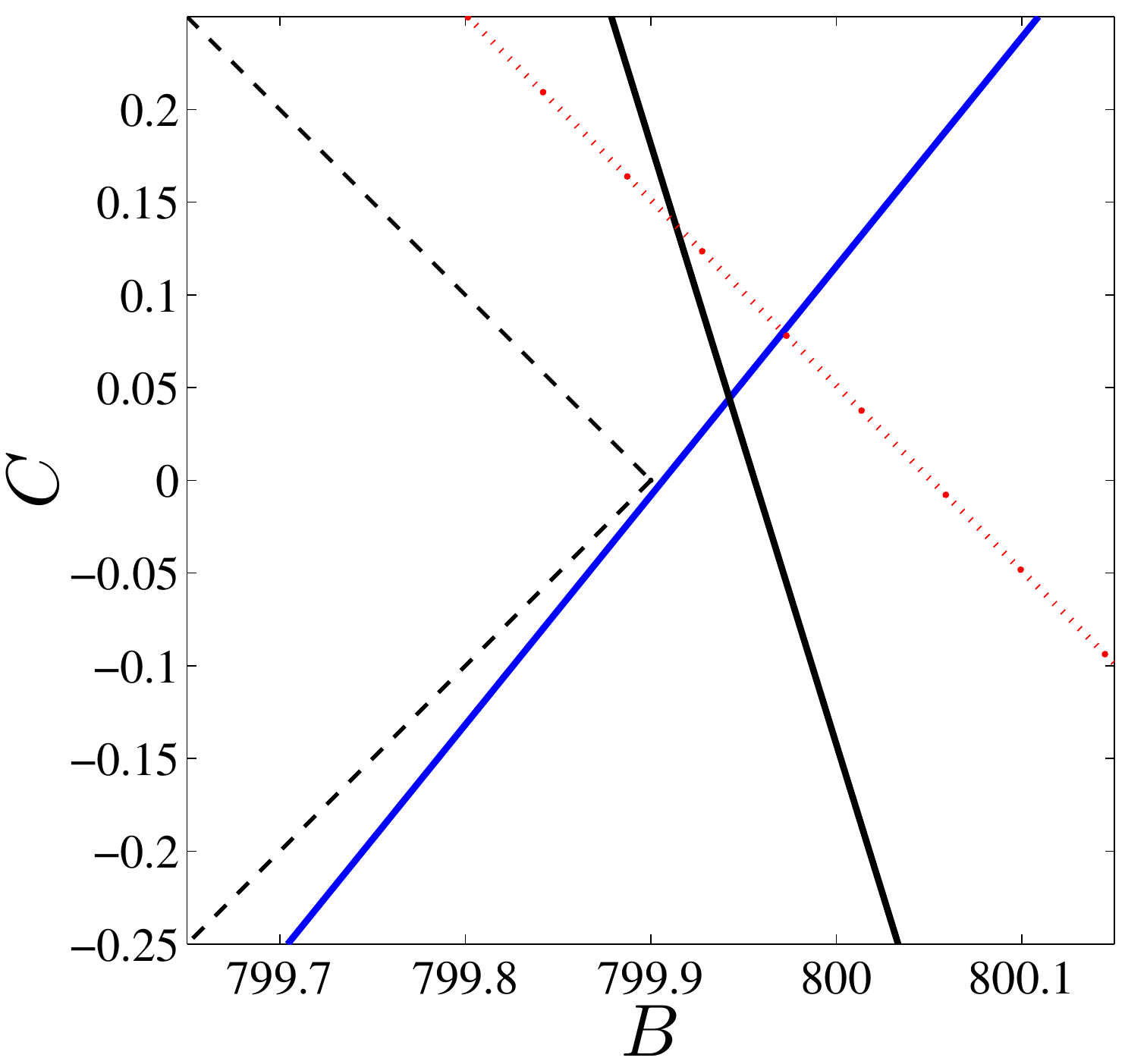}
    \end{center}
  \end{minipage}
  \caption{\small{Typical region for $R = \frac{1}{2n+1}$. Example shows $R = 0.199 \approx \frac{1}{5}$ with $A^*_4 = 799.9$. The boundary of the stability region consists of $\Lambda_0$ (violet), $\Gamma_1$ (blue), $\Delta_4$ (dashed red), $\Gamma_4$ (red), and very small segments of $\Gamma_2$ (green) and $\Gamma_3$ (black).}}
  \label{fig:R333_tran}
\end{figure}

When $R = \frac{1}{2n+1}$ and $A \to +\infty$, the stability region again bulges away from the MRS. The stability region again asymptotically reduces to just four bifurcation curves. However, the regions of stability, which lie outside the MRS, now appear in the $2^{nd}$ and $4^{th}$ quadrants. Fig.~\ref{fig:R333_tran} shows $R = 0.199$ and $A^*_4 = 799.9$ with primarily four bifurcation surfaces comprising the boundary of the region of stability. This figure is quite representative of any $R \to \frac{1}{2n+1}$ from below as $A^*_{2n} \to +\infty$. Below we present several lemmas, which analytically show results similar to Lemmas~4.1-3 for $R = \frac{1}{2n}$. The next section will complete the study with numerical results.

For any delay, $R$, $\Lambda_0$ is always one part of the stability boundary, appearing in the $3^{rd}$ quadrant of the $BC$-plane along the MRS with $B + C = -A$. When $R \to \frac{1}{2n+1}$ from below, there is a transition $A^*_{2n} \to \infty$, which results in a degeneracy line that approaches
\[
 B + C = A^*_{2n},
\]
is parallel to $\Lambda_0$, and lies on the opposite side of the MRS. As in the previous case ($R = \frac{1}{2n}$), $\Gamma_1$ in the $4^{th}$ quadrant creates another edge of the stability region. Finally, the stability region for $R = \frac{1}{2n+1}$, asymptotically finds $\Gamma_{2n}$ mirroring $\Gamma_1$ in the $2^{nd}$ quadrant to complete the simple enlarged stability region.

Below we present a few lemmas to prove some of our claims.

\begin{lemma} For $R = \frac{1}{2n+1}$, one boundary of the region of stability
is the limiting line
\[
 B + C = A_{2n}^*,
\]
which lies on the MRS.
\end{lemma}

\proof
The proof of this lemma closely parallels the proof of Lemma~\ref{lem_lim_line_even}. With $R$ near $\frac{1}{2n+1}$, we consider the transition $A^*_{2n}$, and $\Delta_{2n}$ satisfies
\[
  A = A_{2n}^*, \quad B + C = B_{2n}^* + C_{2n}^*.
\]
Since $A_{2n}^* = -\frac{(2n)\pi}{1-R} \cot\left(\frac{(2n)R\pi}{1-R}\right)$, Eqns.~(\ref{bc_alt}) give
\[
  B_{2n}^* + C_{2n}^* = \textstyle{\frac{(2n)\pi}{1-R}\csc\left(\frac{(2n)R\pi}{1-R}\right)}.
\]
It follows that
\[
 \lim_{R \to \frac{1}{2n+1}^-} \frac{B_{2n}^*+C_{2n}^*}{A_{2n}^*}
   = \lim_{R \to \frac{1}{2n+1}^-}\textstyle{-\sec\left(\frac{(2n)R\pi}{1-R}\right)} = 1.
\]
Thus, $B_{2n}^* + C_{2n}^* \to A_{2n}^*$ for $R < \frac{1}{2n+1}$ as $R \to \frac{1}{2n+1}$.
\qed

For $R = \frac{1}{2n+1}$, $\Gamma_1$ intersects the $\Lambda_0$ at
\[
 (B,C) = \left(\frac{A+2n+1}{2n}, -\frac{(2n+1)(A+1)}{2n}\right).
\]
The next lemma shows that as $R \to \frac{1}{2n+1}$ from below and $A_{2n}^* \to +\infty$, the point $(B_{2n}^*, C_{2n}^*)$ tends to a value symmetric with the origin to the intersection of $\Gamma_1$ and $\Lambda_0$.

\begin{lemma}\label{lem_odd_pt}
For $R < \frac{1}{2n+1}$ and $R \to \frac{1}{2n+1}$, the bifurcation curve $\Gamma_{2n}$ comes to the point
\[
  (B_{2n}^*, C_{2n}^*) = \left(-\frac{A_{2n}^*+2n+1}{2n}, \frac{(2n+1)(A_{2n}^*+1)}{2n}\right)
\]
with $A_{2n}^* \to +\infty$.
\end{lemma}

\proof
This proof is very similar to the proof of Lemma~\ref{lem_even_pt}, so we omit it here.
\qed

Figure~\ref{fig:R333_tran} shows that $\Gamma_1$ and $\Gamma_{2n}$ cross the $C$-axis near the corners of the MRS. The next lemma proves this asymptotic limit, giving more information about the symmetric shape of the region.

\begin{lemma} \label{C_cross_odd}
For $R < \frac{1}{2n+1}$ and $R \to \frac{1}{2n+1}$, the bifurcation curves, $\Gamma_1$ and $\Gamma_{2n}$, pass arbitrarily close to the point $\left(A_{2n}^*,0\right)$ and $\left(-A_{2n}^*,0\right)$, respectively, in the $BC$-plane with $A_{2n}^* \to +\infty$.
\end{lemma}

\proof
The argument for $\Gamma_1$ passing arbitrarily close to $(A_{2n}^*,0)$ is almost identical to the argument given in Lemma~\ref{C_cross_even}. Similarly, $\Gamma_{2n}$ has $\frac{(2n-1)\pi}{1-R} < \omega < \frac{(2n)\pi}{1-R}$, which has $\omega = 2n\pi$ inside the interval. Since this is an even multiple of $\pi$, the $\cos(\omega) \to 1$. Otherwise, the arguments parallel Lemma~\ref{C_cross_even} and so $\Gamma_{2n}$ passes arbitrarily close to $\left(-A_{2n}^*,0\right)$.
\qed

This completes our analytic results to date. The next section provides more numerical details to support our claims of increased stability regions for delays of the form $R = \frac{1}{n}$ with $A$ large.

\setcounter{equation}{0}
\setcounter{theorem}{0}
\setcounter{figure}{0}
\setcounter{table}{0}
\section{Asymptotic Stability Region for $R = \frac{1}{n}$}

The previous section presents the simple stability regions for $R$ near $\frac{1}{n}$, and we gave some details for $R = 0.249$ on how the stability region evolves as $A$ increases. In this section we present more details from numerical studies to convince the reader of the enlarged stability region for these specific rational delays and provide some measure of the size increase compared to the MRS.

Geometrically, for a fixed $A$ the MRS is a square region with one side bounded by $\Gamma_0$. Definition~\ref{family_defn} shows that rational delays, $R$, create families of smooth curves with similar properties and following similar trajectories. The limited number of family members for a fixed $R = \frac{1}{n}$ creates a type of harmonic, which prevents the complete collection of bifurcation curves from asymptotically approaching all sides of the MRS. Hence, the stability region for $R = \frac{1}{n}$ is enlarged with larger regions for smaller values of $n$. This section provides some details and numerical results for the structure and size of the asymptotic stability region.

\subsection{$R$ near $\frac{1}{4}$}

To illustrate our analysis we concentrate on the case $R = \frac{1}{4}$ and note that the arguments generalize to other delays of the form $R = \frac{1}{n}$. For $R = \frac{1}{n}$, Definition~\ref{family_defn} gives $2(n-1)$ distinct family members, which implies $R = \frac{1}{4}$ has six distinct family members. As noted above the key bifurcation curves on the boundary of the stability region asymptotically as $A^*_3 \to +\infty$ with $R \to \frac{1}{4}$ from below are $\Gamma_1$ and $\Gamma_3$ with $\Lambda_0$ and $\Delta_3$, creating the other two boundaries along the MRS. $\Gamma_1$ and $\Gamma_3$ are obviously the first bifurcation curves of the first and third families, while $\Delta_3$ arises from the singular point $\omega = \frac{3\pi}{1-R}$ between $\Gamma_3$ and $\Gamma_4$.

\begin{figure}[htb]
\begin{center}

\includegraphics[width=4.0in]{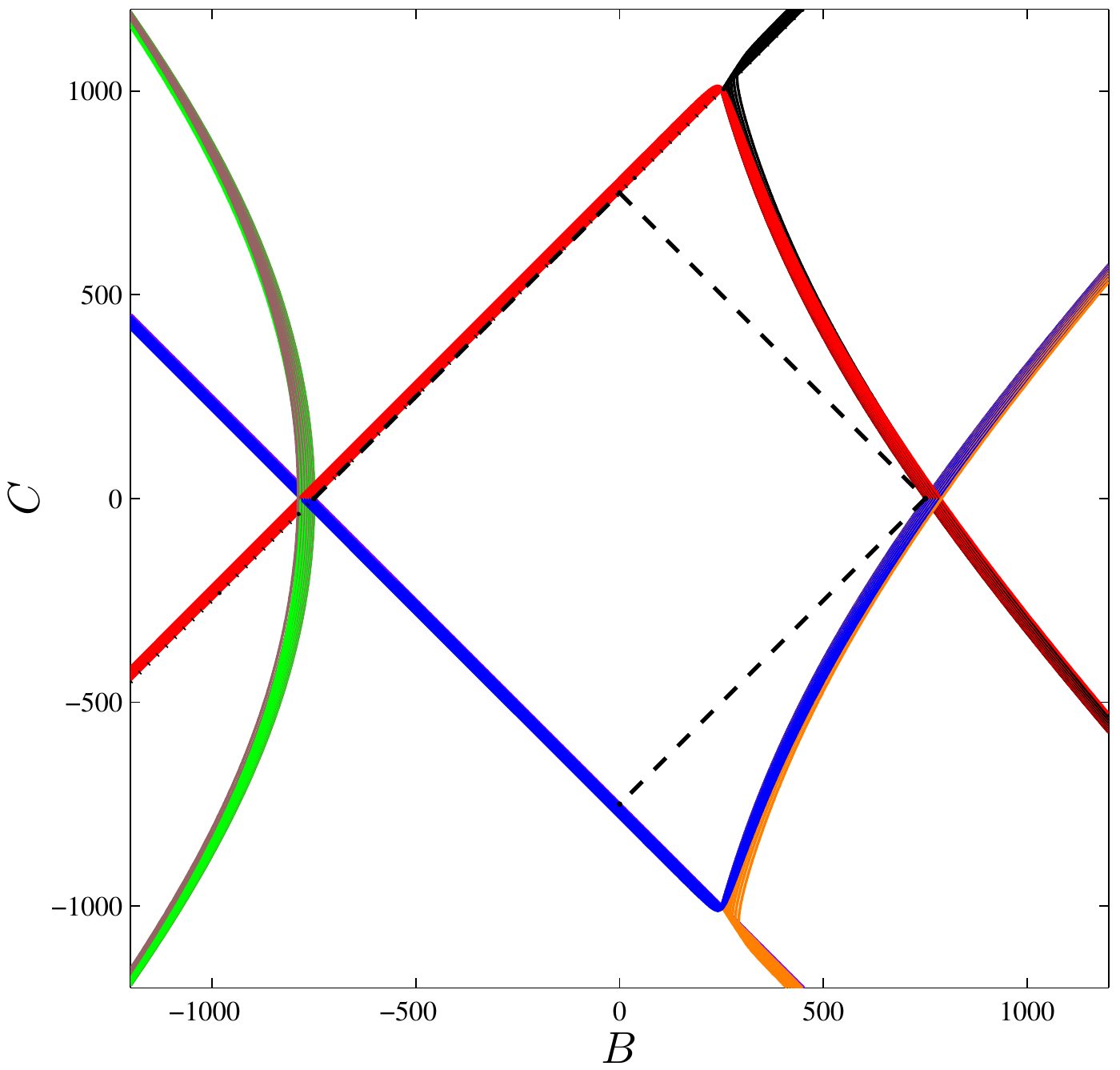}

\begin{tabular}{cccc}
\includegraphics[width=0.23\textwidth]{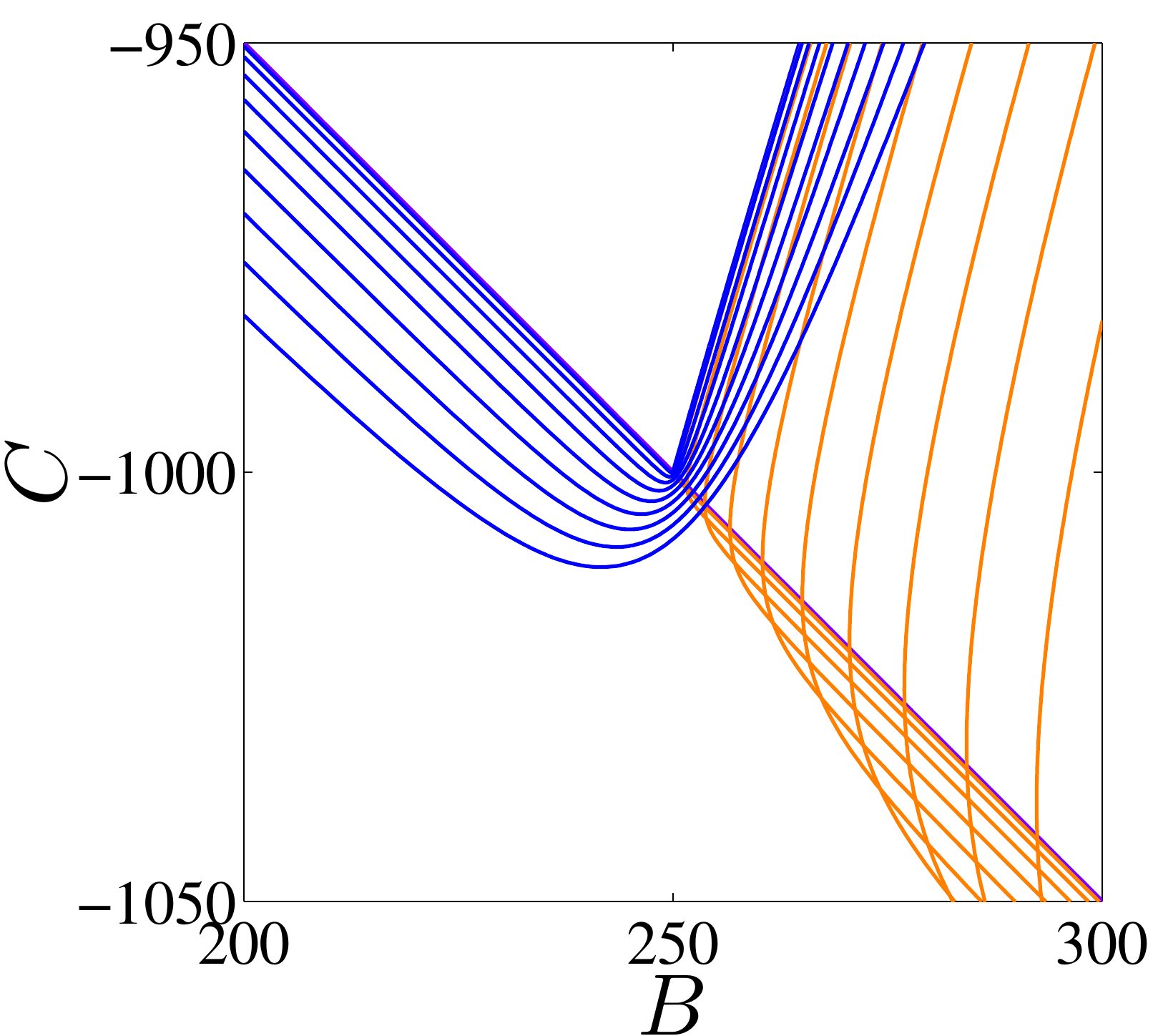}
&
\includegraphics[width=0.22\textwidth]{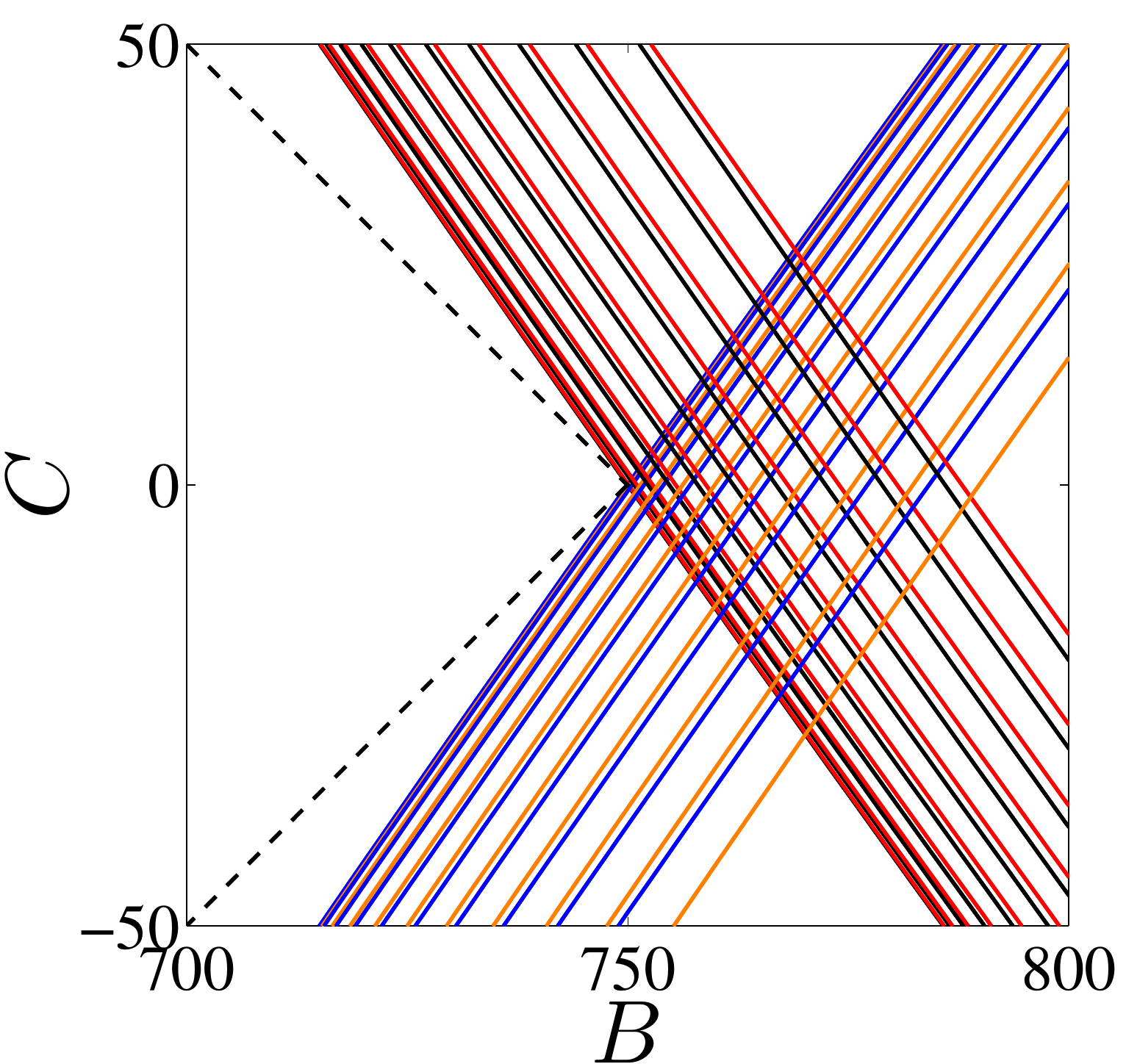}
&
\includegraphics[width=0.22\textwidth]{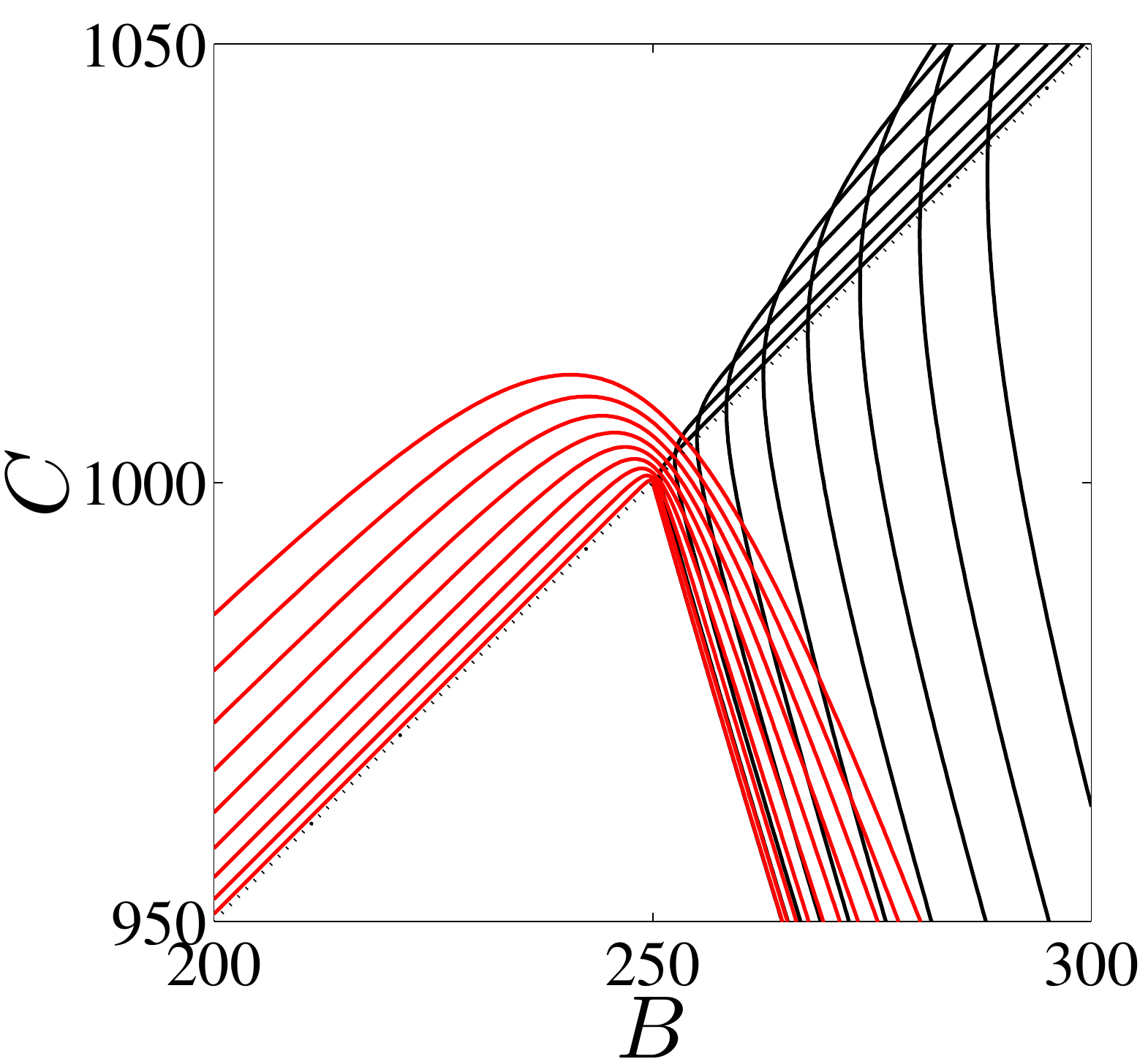}
&
\includegraphics[width=0.22\textwidth]{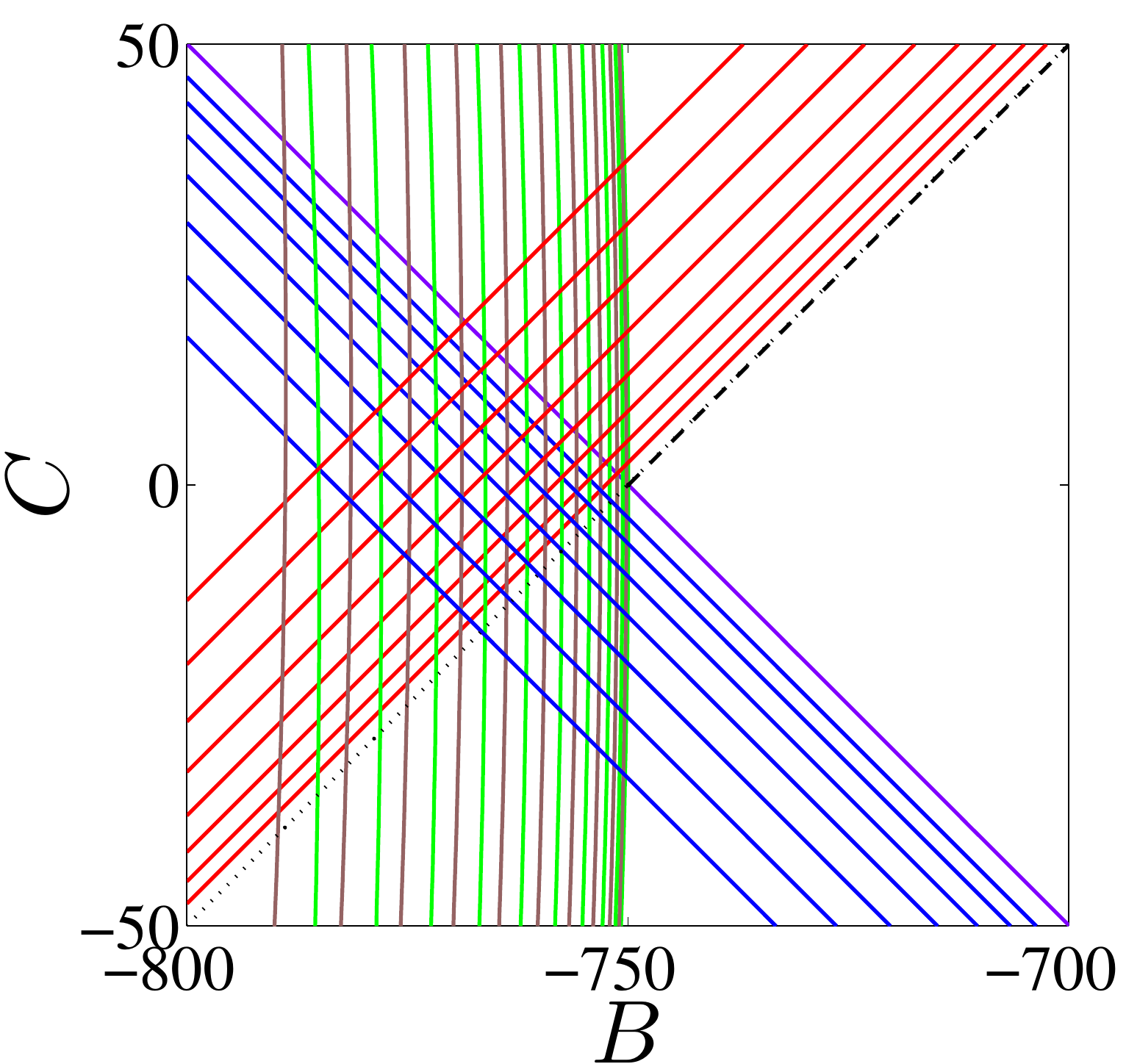}
\end{tabular}
\caption{\small{Ten bifurcation curves for each of the six families for $R = 0.249$ at
       $A_3^* = 749.93$ with close-ups at the corners of the MRS.}}
\label{fig:R249_fam}
\end{center}
\end{figure}

Fig.~\ref{fig:R249_fam} shows 60 bifurcation curves for $R = 0.249$ at $A^*_3 = 749.93$. By continuity of the characteristic equation (\ref{chareqn}), we expect similarities between $R = \frac{1}{4}$ and $R = 0.249$, particularly for $A < A^*_3$. Fig.~\ref{fig:R249_fam} shows clearly the six family structure we predict for $R = \frac{1}{4}$. (Note that $R = 0.249$ is predicted to have 1502 families by Definition~\ref{family_defn}, which will ultimately result in a much closer approach of the stability region to the MRS as $A \to +\infty$.) The figure shows the distinct ordering of the family members within the six families and the characteristic pattern of each of the six families. The coloring pattern in Figure~\ref{fig:R249_fam} follows $\Lambda_0$ in violet and then the six successive families in blue, green, black, red, gray, and orange.

Fig.~\ref{fig:R249_tran} shows the five curves on the boundary of the stability region. As noted before, $\Lambda_0$ is always on the boundary in the $3^{rd}$ quadrant. It connects to $\Gamma_1$, the $1^{st}$ member of the first family. The first close-up in Fig.~\ref{fig:R249_fam} shows the organization of the first family with all members lying outside the boundary of the region of stability with each successive member further away. This first close-up also shows the $6^{th}$ family, which parallels the first family along the boundary in the $4^{th}$ quadrant, then diverges opposite the first family below $\Lambda_0$. Thus, in the $4^{th}$ quadrant the boundary of the stability region consists only of the curves from $\Lambda_0$ and $\Gamma_1$.

In the $2^{nd}$ quadrant, the primary boundary is $\Delta_3$. By Lemma~\ref{lem_lim_line_even}, as $R \to \frac{1}{4}$, $\Delta_3$ approaches the line $C - B = A^*_3$, which is visible in Fig.~\ref{fig:R249_fam}. Since $R = 0.249 < \frac{1}{4}$, there is a small gap between $\Delta_3$ and the boundary of the MRS, so we see a small segment of $\Gamma_2$ on the boundary of the stability region, which is visible in the final close-up of Fig.~\ref{fig:R249_fam}. The line $\Delta_3$ extends into the $1^{st}$ quadrant. We see that outside the tiny segment of $\Gamma_2$ (left of the MRS), the second and fifth families are outside the stability region running symmetrically about the $C$-axis through the $2^{nd}$ and $3^{rd}$ quadrants.

In the $1^{st}$ quadrant, we see $\Gamma_3$ on the boundary of the stability region. The third close-up of Fig.~\ref{fig:R249_fam} shows the ordering of the third family with all members successively outside $\Gamma_3$. This figure also shows all members of the fourth family outside (and intertwined with) the third family, paralleling each other in the first quadrant. The fourth family diverges opposite the third family parallel to $\Delta_3$.

Numerically, our MatLab programs allow us to observe the stability region for any delay $R$ slightly less than $\frac{1}{4}$ at $A^*_3$, and the boundary of the stability region is almost identical to Fig.~\ref{fig:R249_fam}, except that $A^*_3$ increases with $R \to \frac{1}{4}$, expanding the scales of the $B$ and $C$ axes. To further understand the evolution of this stability surface as $R \to \frac{1}{4}$ from below, we detail how Fig.~\ref{fig:Rall} can be used to explain the process. As noted earlier, the key elements causing the bulges in the stability region are the limited number of families of bifurcation curves and the transitions, which distort the boundary. At $R = \frac{1}{4}$, all the transitions $A^*_{3n} \to +\infty$, $n = 1,2,...$, which is easily verified by Definition~\ref{atrans}. Furthermore, it can be shown using perturbation analysis that the transitions have a distinct ordering,
\[
 A^*_3 > A^*_6 > ... > A^*_{3n} > A^*_{3(n+1)} > ...
\]
for delays $R < \frac{1}{4}$. (Details of this proof are not included, but depend on $n$ and $R$, as would be expected.) This sequence allows one to readily determine changes to the boundary of the stability region.

Table~\ref{summarytable249} provides the complete evolution of the stability surface for $R = 0.249$ from its beginning at $A_0 \approx -5.016$ until $A^*_3 \approx 749.93$. We call attention to the sequence of reverse tangencies and the reverse transferral for $A < A^*_3$. These events all occur just prior to one of the transitions, $A^*_{3n}$. The transition $A^*_3$ creates $\Delta_3$ with $C - B = C^*_3 - B^*_3$ at $\omega = \frac{3\pi}{1-R}$, which becomes part of the boundary of the stability region. The transition, $A^*_6 \approx 749.72$, results in $\Delta_6$ passing through the point $(B^*_6,C^*_6) \approx (250.05, -1000.19)$ with the line parallel and below $\Lambda_0$. At $A^*_6$, $\Gamma_1$ intersects $\Lambda_0$ at $(B, C) \approx (247.91, -999.63)$, which lies closer to the stability region than $(B^*_6,C^*_6)$. Thus, this transition pulls $\Gamma_6$ and $\Gamma_7$ outside the stability region, swapping the directions in which the curves go to infinity, giving $\Gamma_7$ its flow paralleling $\Gamma_1$ and $\Lambda_0$ and maintaining its position outside the stability region. This distortion from the transition $A^*_6$ results in the reverse transferral, $\tilde{A}^z_{6,1} \approx 749.4$, just prior to the transition.

In a similar fashion, $A^*_9 \approx 749.37$ creates $\Delta_{9}$, which passes through $(B^*_9, C^*_9) \approx (250.10, 1000.42)$ and is parallel to $\Delta_3$. We note that $\Delta_3$ through $(B^*_3, C^*_3) = (250.01, 1000.05)$. Again $(B^*_9, C^*_9)$ is outside the stability region pulling $\Gamma_9$ and $\Gamma_{10}$ to their positions paralleling, but outside $\Gamma_3$ in the $1^{st}$ quadrant. Subsequently, $\Gamma_{10}$ falls in order with other members of the fourth family, which is seen with the red curves of Fig.~\ref{fig:R249_fam}. The distortion from $A^*_9$ and reordering of the curves is what produces the reverse tangency, $\tilde{A}^t_{9,3} \approx 747.13$, simplifying the composition of the boundary of the stability region.

\begin{figure}[htb]
\begin{center}

\includegraphics[width=3.5in]{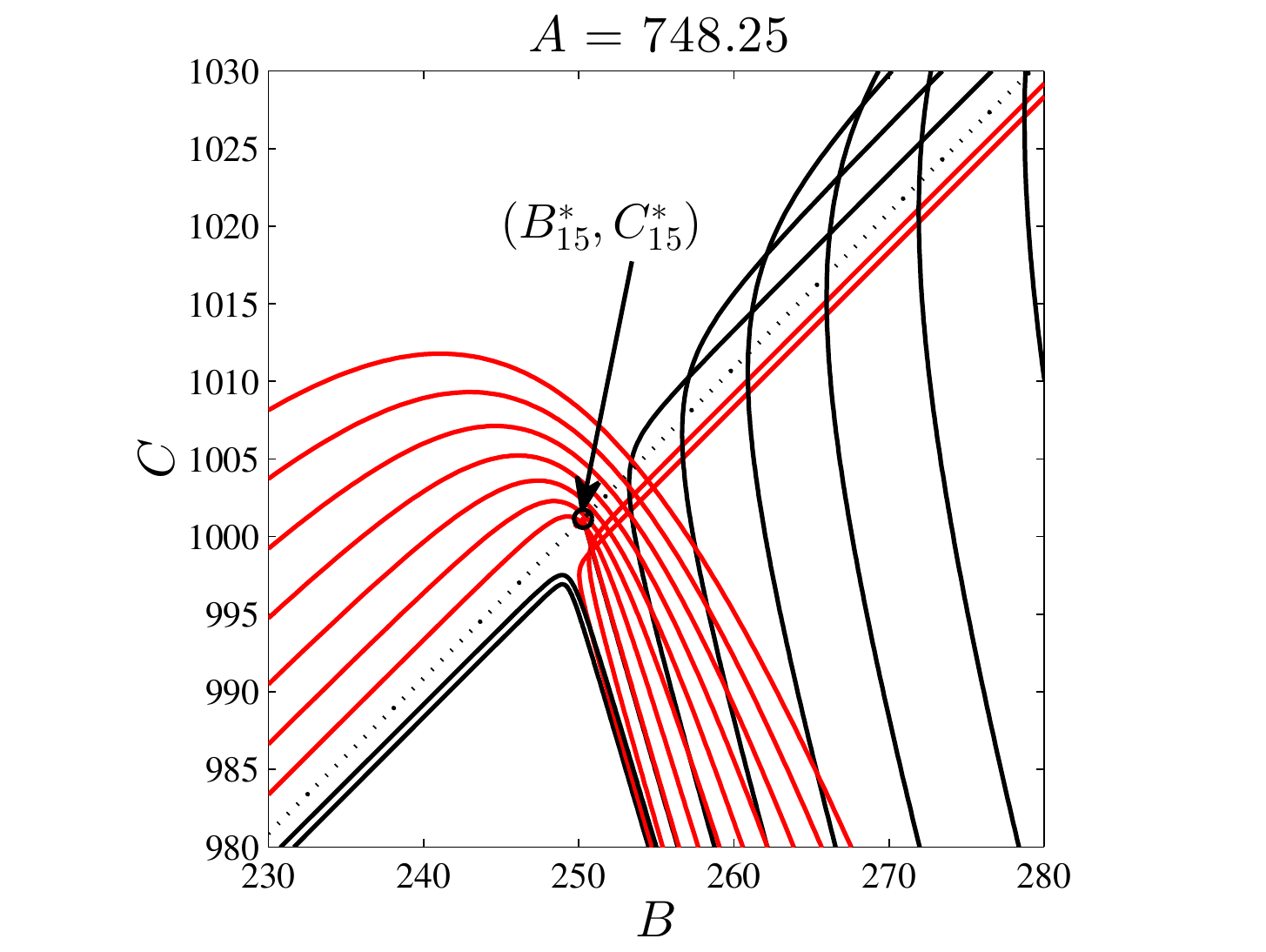}

\caption{\small{Nine and ten bifurcation curves for the third and fourth families, respectively, for $R = 0.249$ at $A_{15}^* = 748.25$ are shown, including the $\Delta_{15}$. Note that $\Gamma_3$ and $\Gamma_9$ remain close to the boundary of the stability region with $\Gamma_3$ constructing this portion of the boundary. $\tilde{A}^t_{9,3} \approx 747.134$ has recently occurred, removing $\Gamma_9$ from the boundary of the stability region.}}
\label{fig:R249_intermed}
\end{center}
\end{figure}

Following the alternating pattern, we find $A^*_{12} \approx 748.88$, producing $\Delta_{12}$ passing through $(B^*_{12}, C^*_{12})$ $\approx (250.18, -1000.74)$, which is parallel to $\Lambda_0$. Now $(B^*_{12}, C^*_{12})$ is outside the stability region pulling $\Gamma_{13}$ to a position paralleling, but outside $\Gamma_7$. The distortion from $A^*_{12}$ and the reorganization of the curves create the reverse tangency, $\tilde{A}^t_{12,6} \approx 743.46$, losing the curve $\Gamma_{12}$ from the boundary of the stability region. Fig.~\ref{fig:R249_intermed} shows the boundary of the stability region at $A^*_{15}$, where $\Delta_{15}$ is formed. This pulls $\Gamma_{15}$ and $\Gamma_{16}$ outside the stability region, which earlier resulted in the reverse tangency, $\tilde{A}^t_{15,9} \approx 738.19$, and $\Gamma_{15}$ leaving the boundary of the stability region. Fig.~\ref{fig:R249_intermed} shows the $4^{th}$ family (red) lined sequentially outside the region of stability for $\Gamma_{22}$, $\Gamma_{28}$, ... with the $3^{rd}$ family (black) paralleling $\Gamma_3$ along the upper right boundary of the stability region, then moving away, except for $\Gamma_{3}$ and $\Gamma_{9}$. Since $\tilde{A}^t_{9,3} \approx 747.13$, $\Gamma_9$ has just left the boundary of the stability region and soon transitions with $\Gamma_{10}$ at $A^*_9 \approx 749.37$, causing $\Gamma_{10}$ to follow the pattern of the other members of the $4^{th}$ family outside the region of stability.

As seen in Table~\ref{summarytable249}, there is an alternating pattern of reverse tangencies as we progress to lower values of $A$, and each reverse tangency simplifies the boundary of the stability region and organizes the families into the pattern seen in Fig.~\ref{fig:R249_fam} because of one of the $A^*_{3n}$ transitions. This same sequence of events occurs for each $R < \frac{1}{4}$ (sufficiently close) with more tangencies and reverse tangencies before $A^*_3$ as $R \to \frac{1}{4}$ from below and $A^*_3$ getting larger. Thus, the geometric orientation of the curves and the sequence of reverse tangencies and transferrals are virtually identical to the figures shown with $R = 0.249$, except the $B$ and $C$ scales increase as $R \to \frac{1}{4}$.

The continuity of the characteristic equation shows that all delays $R < \frac{1}{4}$, yet sufficiently close to $R = \frac{1}{4}$, will generate a simplified stability region very similar to Fig.~\ref{fig:R249_tran} at $A^*_3$. As $R \to \frac{1}{4}$ from below, $\Delta_3$ gets closer to the MRS and the contribution of $\Gamma_2$ on the boundary shrinks. We have been unable to definitively prove whether $\Gamma_2$ persists on the boundary of the stability region for $R < \frac{1}{4}$ at $A_3^*$ or if $R$ sufficiently close to $\frac{1}{4}$ sees $\Gamma_2$ exit the boundary of the stability region. This simple shape of the stability region allows easy numerical computation of how enlarged the stability region is relative to the MRS. Table~\ref{table:R1over4} gives the relative increase of this region of stability for several values of $R \to \frac{1}{4}$, indicating the predicted asymptotic increase in size of the stability region for $R = \frac{1}{4}$ at $A^*_3 = +\infty$. Since the transitions, $A^*_{3n}$, all occur at $+\infty$ when $R = \frac{1}{4}$, continuity of the characteristic equation suggests this enlarged stability region persists for $R = \frac{1}{4}$ and should be 26.86\% larger than the MRS. For any $R < \frac{1}{4}$ and $A > A^*_3$, the six family structure breaks down, leading to new tangencies and a new ordering of the larger families, which results in significantly smaller regions of stability.

\begin{table}[htb]
\begin{center}
\begin{tabular}{|c|c||c|c|} \hline
$R$ & Area Ratio & $R$ & Area Ratio \\ \hline
 0.249 &    1.2687437 & 0.24999 & 1.2686377 \\ \hline
 0.2499 &    1.2686388 & 0.249999 & 1.2686377 \\ \hline
\end{tabular}
\caption{\small{Numerical computation of the region of stability at $A^*_3$ near $R = \frac{1}{4}$. Area given as ratio of stability region to MRS.}}
\label{table:R1over4}
\end{center}
\end{table}

\subsection{Increased Area for $R = \frac{1}{n}$}

Figs.~\ref{fig:R249_tran} and \ref{fig:R333_tran} show the simple bifurcation curve structure on the typical boundary of stability region for $R$ near $\frac{1}{n}$ at $A^*_{n-1}$. The symmetrical shape varies depending on whether $n$ is even or odd, but all these stability regions at $A^*_{n-1}$ reduce to having $\Lambda_0$ on one edge, $\Delta_{n-1}$ on another, and the bifurcation curves $\Gamma_1$ and $\Gamma_{n-1}$ comprising the majority of the remaining two edges of the boundary of the region of stability. (The gap between $\Delta_{n-1}$ and the MRS allows small segments of $\Gamma_2$ and possibly $\Gamma_{n-2}$ to remain on the boundary of the stability region at $A_{n-1}^*$, shrinking as $R \to \frac{1}{n}$.) $\Gamma_1$ and $\Gamma_{n-1}$ bulge out from the MRS, leaving an increased region of stability, which is readily computed.

For $R = \frac{1}{n}$, Def.~\ref{family_defn} gives $2(n-1)$ families of bifurcation curves. These families organize much in the same way as shown in the previous section ($R$ near $\frac{1}{4}$) to help maintain the increased regions of stability for $R = \frac{1}{n}$ over the MRS. The orderly family structure allows one to study each $R = \frac{1}{n}$ much as we did in the previous section and observe similar sequences of transferrals, tangencies, reverse tangencies, and reverse transferrals, which ultimately lead to the simple structure of the stability region seen in Figs.~\ref{fig:R249_tran} and \ref{fig:R333_tran} at $A^*_{n-1}$ for $R < \frac{1}{n}$, sufficiently close. Using the continuity (pointwise) of the characteristic equation, we claim that the bulge in the region of stability persists for $R = \frac{1}{n}$.

\begin{table}[htb]
\begin{center}
\begin{tabular}{|c|c|c||c|c|c|} \hline
$R$ & Area Ratio & Linear Extension & $R$ & Area Ratio & Linear Extension\\ \hline
 $\frac{1}{2}\rule[-0.08in]{0in}{.25in}$  &   2.0000  & 1.0000 &
    $\frac{1}{7}\rule[-0.08in]{0in}{.25in}$  &   1.1084  & 0.1667\\ \hline
 $\frac{1}{3}\rule[-0.08in]{0in}{.25in}$  &   1.4431  & 0.5000 &
    $\frac{1}{8}\rule[-0.08in]{0in}{.25in}$  &   1.0878  & 0.1429 \\ \hline
 $\frac{1}{4}\rule[-0.08in]{0in}{.25in}$  &   1.2686  & 0.3333 &
   $\frac{1}{9}\rule[-0.08in]{0in}{.25in}$  &   1.0729  & 0.1250 \\ \hline
 $\frac{1}{5}\rule[-0.08in]{0in}{.25in}$  &   1.1859  & 0.2500 &
   $\frac{1}{10}\rule[-0.08in]{0in}{.25in}$ &   1.0617  & 0.1111 \\ \hline
 $\frac{1}{6}\rule[-0.08in]{0in}{.25in}$  &   1.1386  & 0.2000 & & &\\ \hline
\end{tabular}
\caption{\small{Increases in the Region of Stability for $R = \frac{1}{n}$. Area ratio is the ratio of the stability area to the MRS. The linear extension is the ratio of how far the stability region extends past the MRS.}}
\label{table:area_bulge}
\end{center}
\end{table}

We showed the region of stability extends linearly along the MRS by a factor of $\frac{1}{n-1}$ for $R = \frac{1}{n}$. The simple bifurcation curve structure allows easy numerical computation of this area bulging from the MRS. Table~\ref{table:area_bulge} gives the size of the increased region of stability for various $R = \frac{1}{n}$. Asymptotically, the region of stability for $R = \frac{1}{2}$ is triangular with only $\Lambda_0$, $\Gamma_1$, and $\Delta_{1}$. This results in a region that is twice the size of the MRS, asymptotically. As the denominator increases, the asymptotic region of stability decreases relative to the MRS, yet it is still over 10\% larger than the MRS when $R = \frac{1}{7}$.

\subsection{Stability Spurs}

The discussion above shows how transitions increase the area of the region of stability for $R = \frac{1}{n}$. Just before a transition, $A_j^*$, bifurcation curve $\Gamma_j$ extends out toward bifurcation curve $\Gamma_{j+1}$, expanding the main region of stability. Numerically, we observe that just prior to $A_j^*$, $\Gamma_{j+1}$ self-intersects creating an island of stability that is disconnected in the $BC$-plane for a fixed $A$. Definition~3.6 describes the stability spurs, which connect to the main stability surface at transitions, $A^*_j$. In this section we provide some details from our numerical simulations about the stability spurs.

\begin{figure}[htb]
    \centering\includegraphics[width=0.55\textwidth]{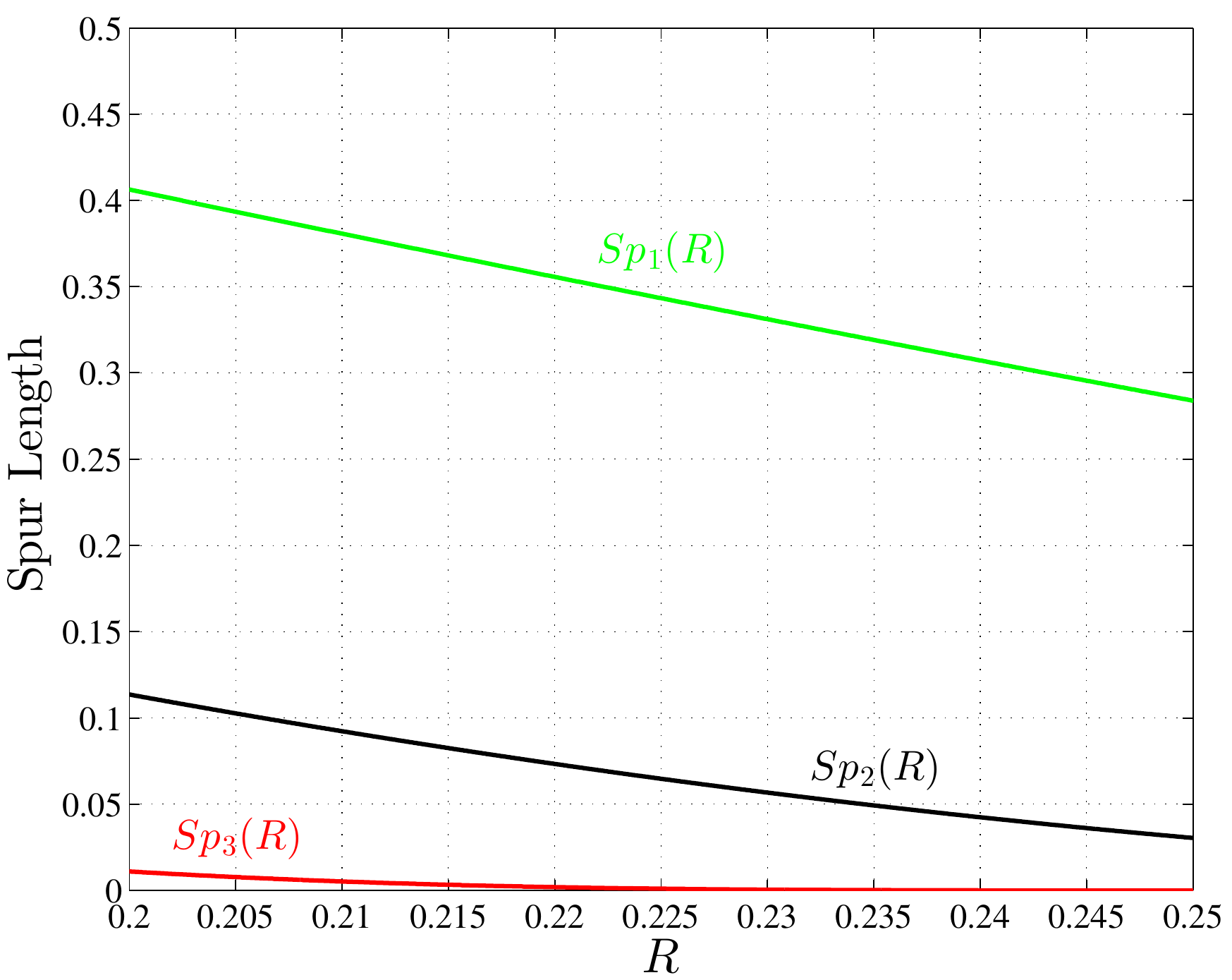}
  \caption{The spur Lengths  associated with $Sp_1(R)$ (green), $Sp_2(R)$ (black) and $Sp_3(R)$ (red) as a function  of $R$.}
  \label{fig:spurlength}
\end{figure}

Numerically, we observe exactly $n-1$ stability spurs for $R \in \left(\frac{1}{n+1}, \frac{1}{n}\right)$. The largest stability spurs, $Sp_1(R)$, correspond to $A_1^*$ for at least a limited range of $R$. We performed an extensive study of the stability spurs for $R \in \left(\frac{1}{5},\frac{1}{4}\right)$. Over this range there appear to be exactly three stability spurs. Fig.~(\ref{fig:spurlength}) shows the variation in length of the three stability spurs, which are all monotonically decreasing in $R$. The length of the first spur for $R \in \left[\frac{1}{5},\frac{1}{4}\right]$ ranges from $|Sp_1(0.2)| = 0.4064$ adjoining the main stability surface at $A_1^*(0.2) = -3.927$ to $|Sp_1(0.25)| = 0.2840$ adjoining the main stability surface at $A_1^*(0.25) = -2.418$. Where this first stability spur adjoins the main stability surface, the cross-sectional area in the $BC$-plane of $Sp_1(0.2)$ is 17.65\% of the total stable region, while the area of $Sp_1(0.25)$ at $A_1^*(0.25)$ is 7.08\% of the total stable cross-section. (See Fig.~\ref{fig:spur}.)

The sizes of the second and third stability spurs are significantly smaller. The length of $Sp_2(0.2)$ is 0.1136, while the length of $Sp_2(0.25)$ is 0.0305. The cross-sectional area of $Sp_2(0.2)$ is 0.3003\% of the total stable region at $A_2^*(0.2) = 0$, while the area of $Sp_2(0.25)$ is 0.01483\% of the total stable region at $A_2^*(0.25) = 4.837$. the length of $Sp_3(0.2)$ is 0.0110, and this spur adjoins the main stability region at $A_3^*(0.2) = 11.78$. Clearly, no spur occurs at $R = \frac{1}{4}$, as $A_3^*(0.25) = +\infty$.

The stability spurs are easy to observe in Fig.~\ref{fig:Rfifth} for $R = \frac{1}{5}$. From a numerical perspective the length of a stability spur is difficult to accurately compute because of the cusp point, $A_j^p$, which is a singularity. Both the length and the self-intersection of $Sp_3(R)$ as $R \to \frac{1}{4}$ becomes impossible to compute, and effectively, $Sp_3(R)$ vanishes as $R \to \frac{1}{4}$. However, as we have already seen, the main stability region bulges out in the $1^{st}$ quadrant because of $A_3^*(R)$, but we have been unable to confirm that $\Gamma_4$ always self-intersects as $A_3^*(R) \to + \infty$ with $R \to \frac{1}{4}$.

\section{Example Continued}

\begin{figure}[htb]
 \begin{center}
 \begin{tabular}{cc}
 \includegraphics[width=0.35\textwidth]{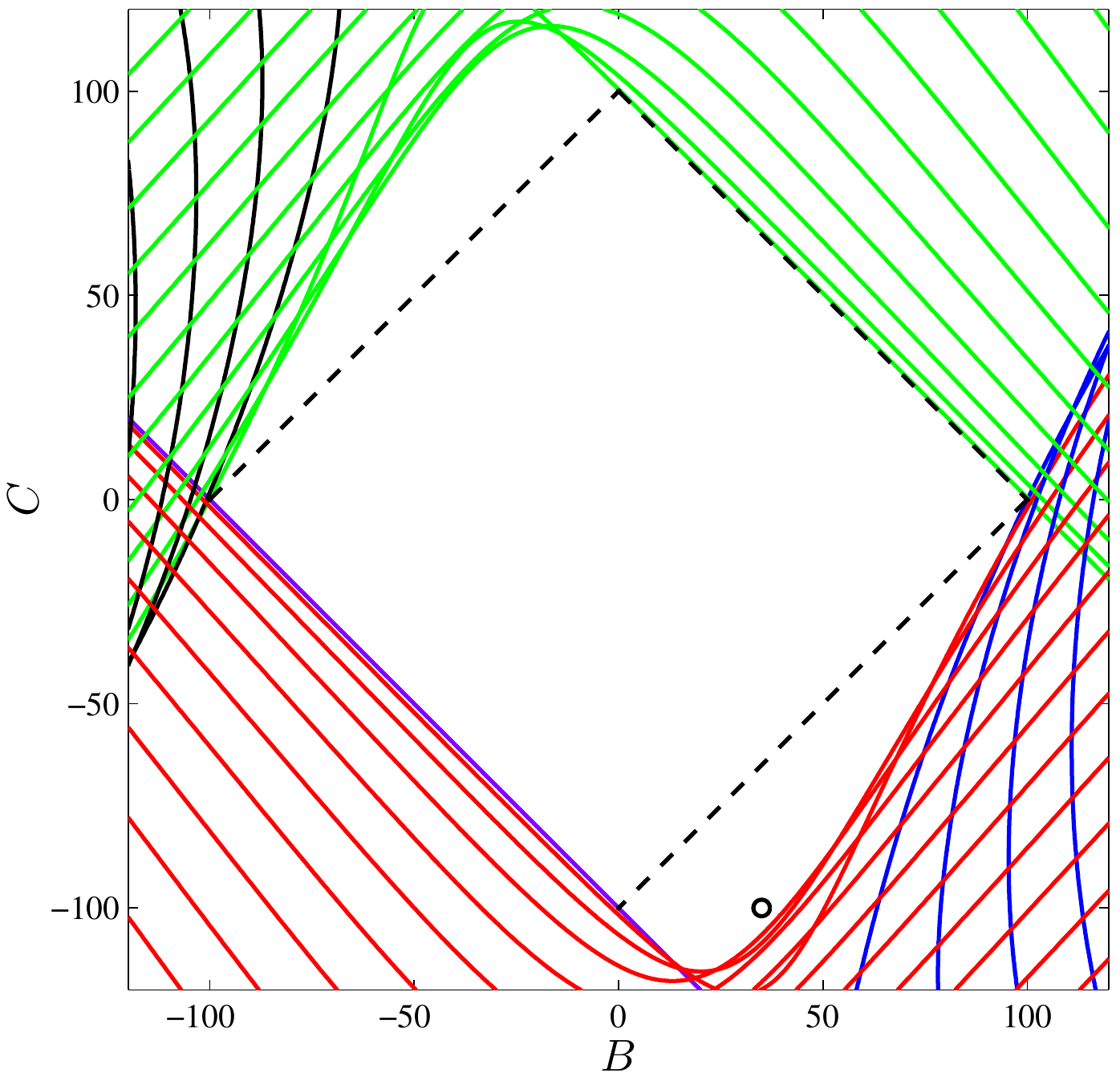}
  &
 \includegraphics[width=0.35\textwidth]{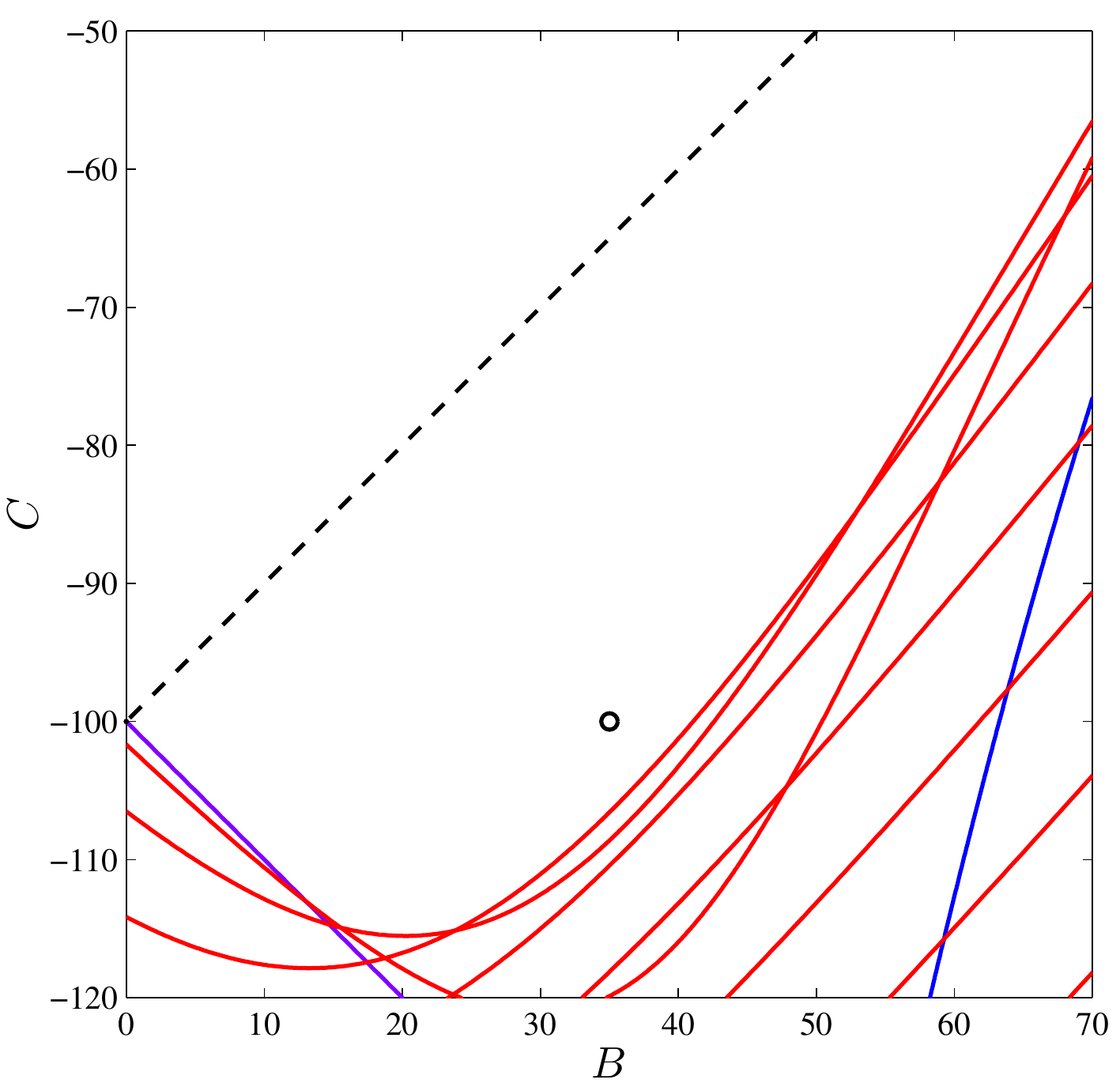}
 \end{tabular}
 \end{center}
  \caption{Bifurcation curves for the modified platelet model with $R = \frac{1}{3}$ and $A = 100$, showing the point where the model exists in the parameter space. This point is clearly in the stable region.}
  \label{fig:model_13}
\end{figure}

In Section~2, we examined an example motivated by work of B{\'e}lair and Mackey \cite{BEM}, and Fig.~\ref{fig:platelet} showed how rational delays stabilized the model. Here we use some of the information above to provide more details about the complex behavior observed in Fig.~\ref{fig:platelet}. With the parameters given in Section 2, the model is linearized about the nontrivial equilibrium, $P_e \approx 5.565$. If $y(t) = P(t) - P_e$, then the approximate linearized model becomes:
\begin{equation}
  \frac{dy}{dt} = -100\,y(t) + 100\,y(t-R) - 35\,y(t-1). \label{lin_model}
\end{equation}
This is Eqn.~(\ref{DDE2}) with $(A, B, C) = (100, 35, -100)$. We use our MatLab program to generate plots of 40 bifurcation curves in the $BC$-plane with $A = 100$ for various values of $R$. Fig.~\ref{fig:model_13} shows the distinct four family feature for the delay $R = \frac{1}{3}$ and the enlarged region of stability. The linearized model is in the region of stability, which agrees with the numerical simulation in Fig.~\ref{fig:platelet}, where the delay $R = \frac{1}{3}$ gives a stable solution.

\begin{figure}[htb]
 \begin{center}
 \begin{tabular}{cc}
 \includegraphics[width=0.35\textwidth]{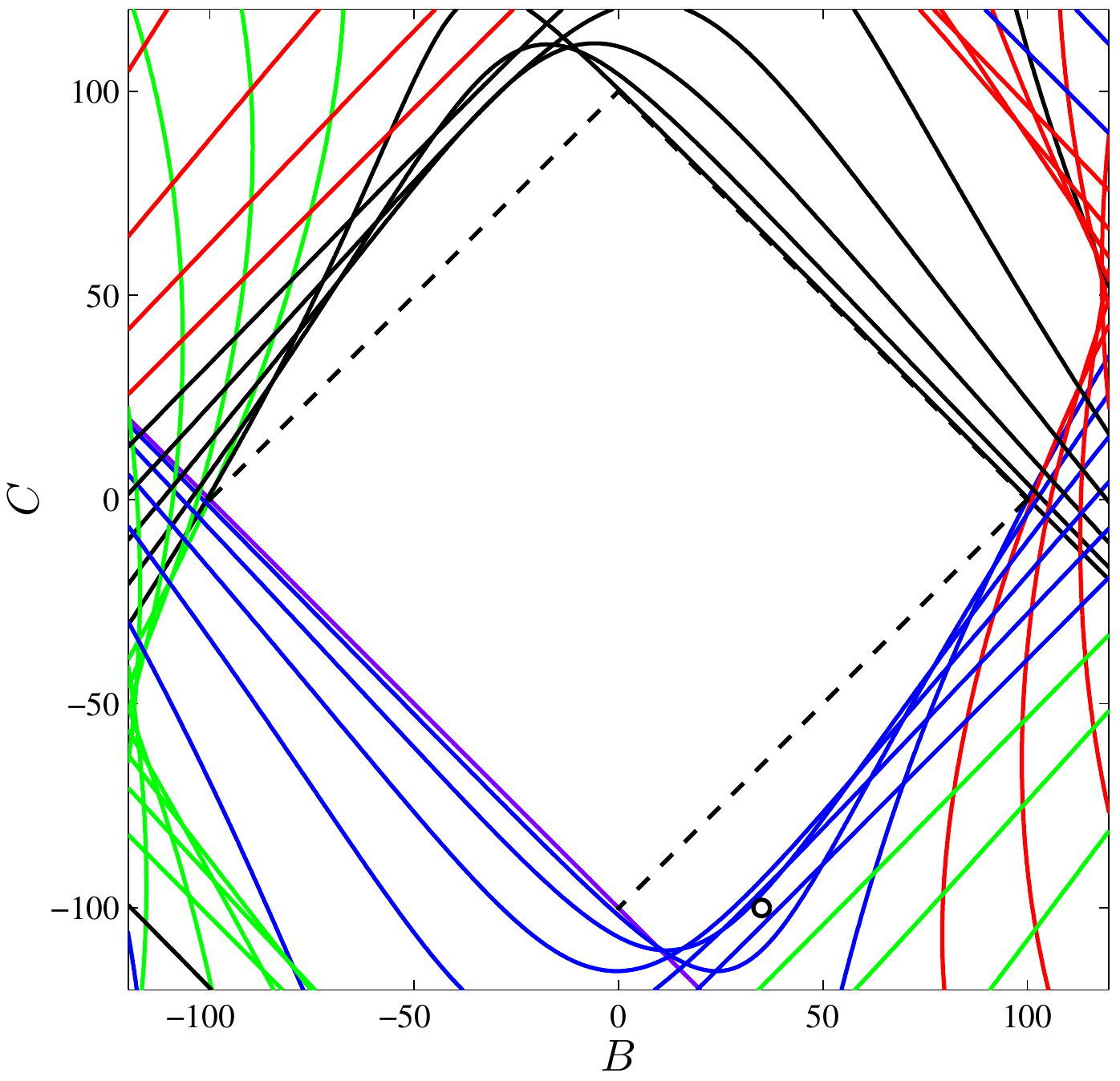}
  &
 \includegraphics[width=0.35\textwidth]{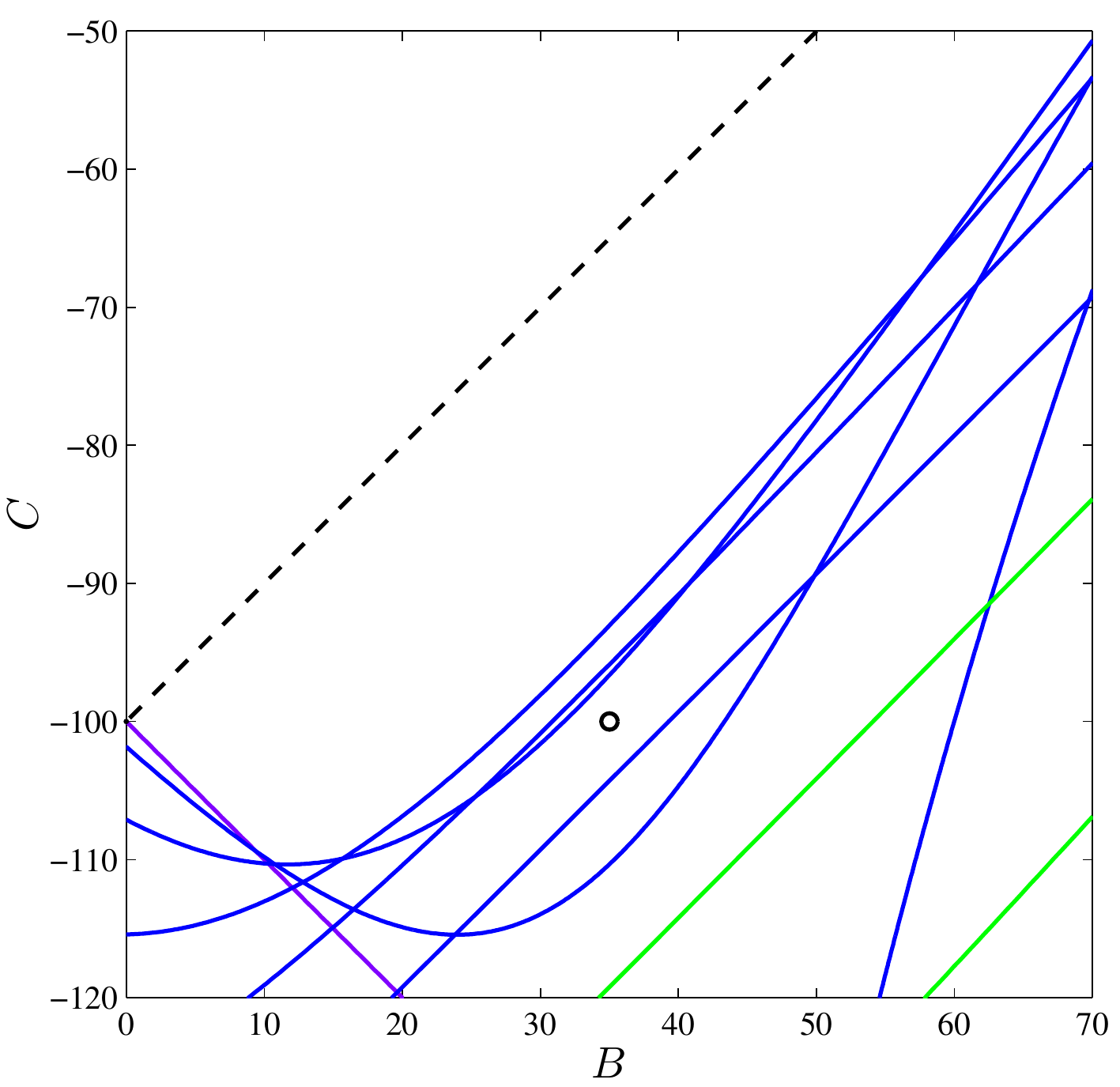}
 \end{tabular}
 \end{center}
  \caption{Bifurcation curves for the modified platelet model with $R = 0.318$ and $A = 100$, showing the point where the model exists in the parameter space. This point is outside the region of stability.}
  \label{fig:model_318}
\end{figure}

When the delay is decreased to $R = 0.318$, there is a transition at $A_2^* = 42.8$. The four family structure rapidly unravels, and the bifurcation curves begin approaching the MRS more closely. Fig.~\ref{fig:model_318} shows that Eqn.~(\ref{lin_model}) has its equilibrium outside the curves $\Gamma_9$, $\Gamma_{13}$, and $\Gamma_{17}$. With the help of Maple (using information from Fig.~\ref{fig:model_318}), the eigenvalues of Eqn.~(\ref{chareqn}) with positive real part are computed. These eigenvalues are:
\[
  \lambda_1 = 0.1056\pm 58.36\,i \qquad \lambda_2 = 0.06238\pm 77.43\,i \qquad \lambda_3 = 0.04914\pm 39.32\,i.
\]
The leading eigenvalue comes from the equilibrium point being furthest from $\Gamma_{13}$, and its frequency suggests a period of $2\pi/58.36 \approx 0.108$, which is close to the period of oscillation seen in Fig.~\ref{fig:platelet} for $R = 0.318$.

\begin{figure}[htb]
 \begin{center}
 \begin{tabular}{cc}
 \includegraphics[width=0.35\textwidth]{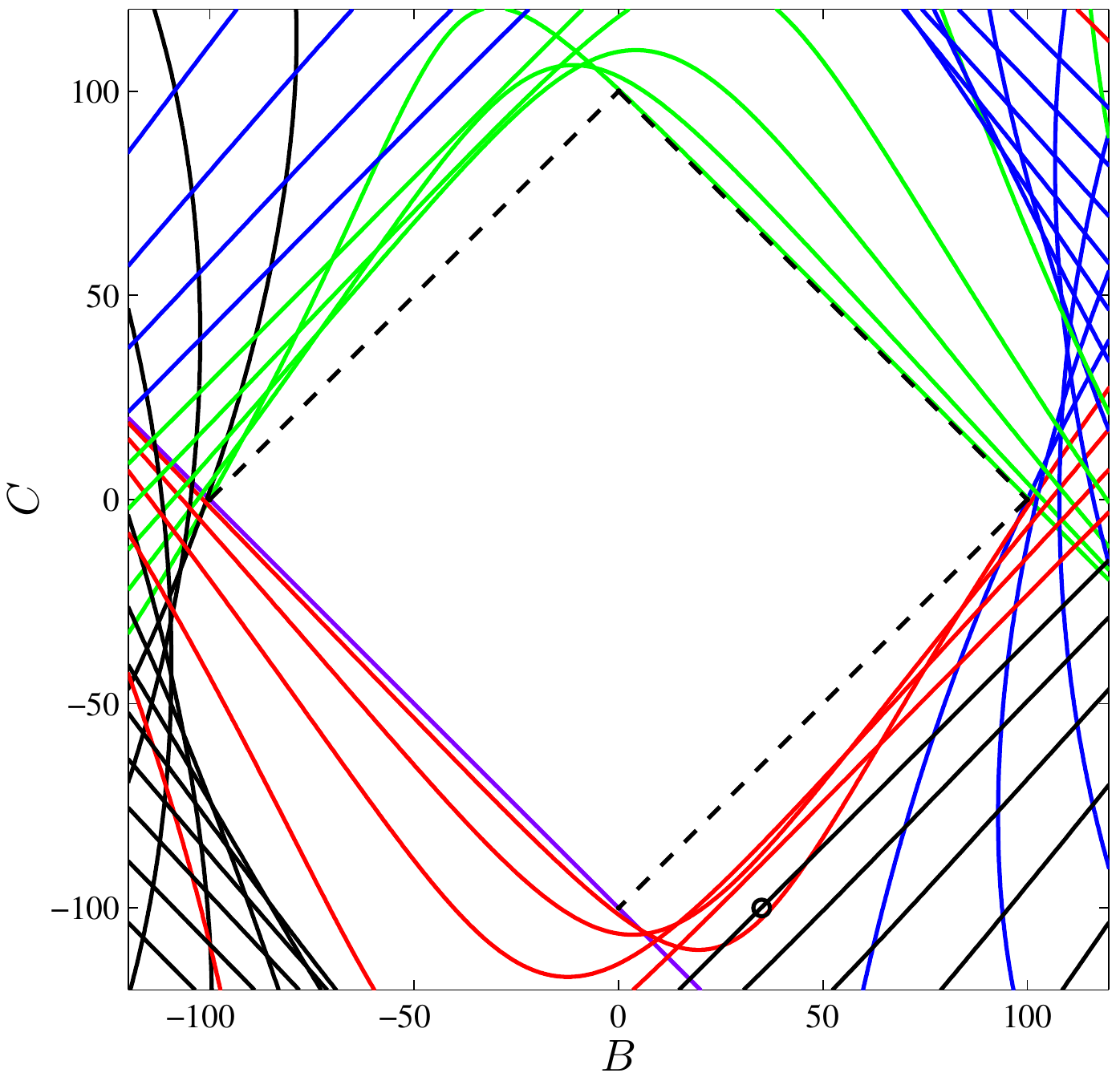}
  &
 \includegraphics[width=0.35\textwidth]{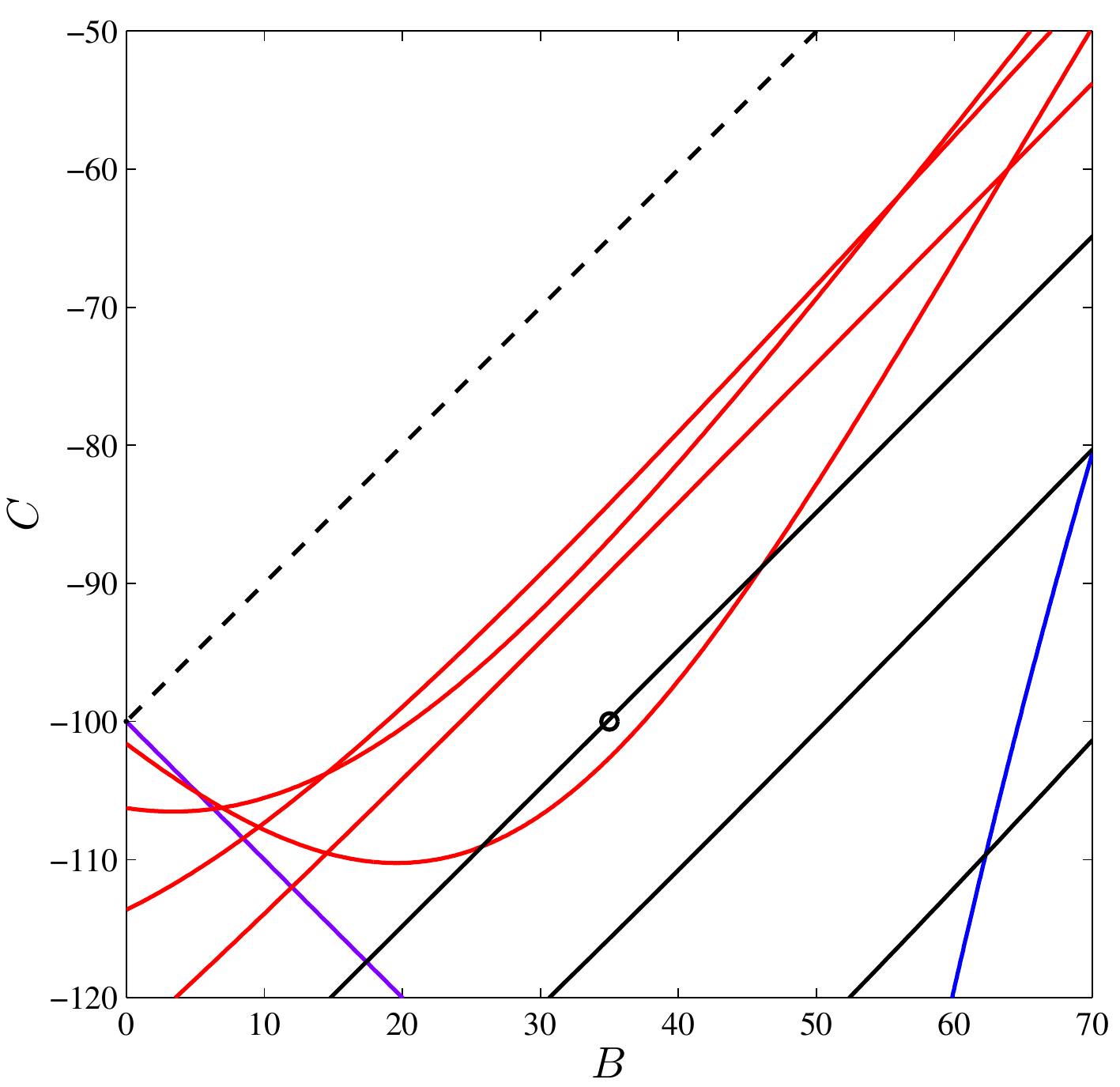}
 \end{tabular}
 \end{center}
  \caption{Bifurcation curves for the modified platelet model with $R = 0.34$ and $A = 100$, showing the point where the model exists in the parameter space. This point is outside the region of stability.}
  \label{fig:model_34}
\end{figure}

A similar analysis can be performed for $R = 0.34$. Since this delay is larger than $R = \frac{1}{3}$, there is not the four family structure. Fig.~\ref{fig:model_34} shows four bifurcation curves, $\Gamma_{8}$, $\Gamma_{12}$, $\Gamma_{16}$, and $\Gamma_{19}$, between the region of stability and where the linearized model is located. As before, it is easy to use this information to determine the eigenvalues:
\[
  \lambda_1 = 0.2424\pm 53.51\,i \qquad \lambda_2 = 0.1988\pm 35.21\,i \qquad
  \lambda_3 = 0.1625\pm 71.87\,i \qquad \lambda_4 = 0.002273\pm 90.31\,i.
\]
The equilibrium point is furthest from $\Gamma_{12}$, and the period from $\lambda_1$ is $2\pi/53.51 \approx 0.117$, which is similar to the period of oscillation seen in Fig.~\ref{fig:platelet} for $R = 0.34$. The next two eigenvalues are moderately large, resulting in the additional irregular structure observed in the simulation. Note that the last eigenvalue with positive real part has just barely crossed the imaginary axis, which again is apparent from Fig.~\ref{fig:model_34}.

\begin{figure}[htb]
 \begin{center}
 \begin{tabular}{ccc}
 \includegraphics[width=0.32\textwidth]{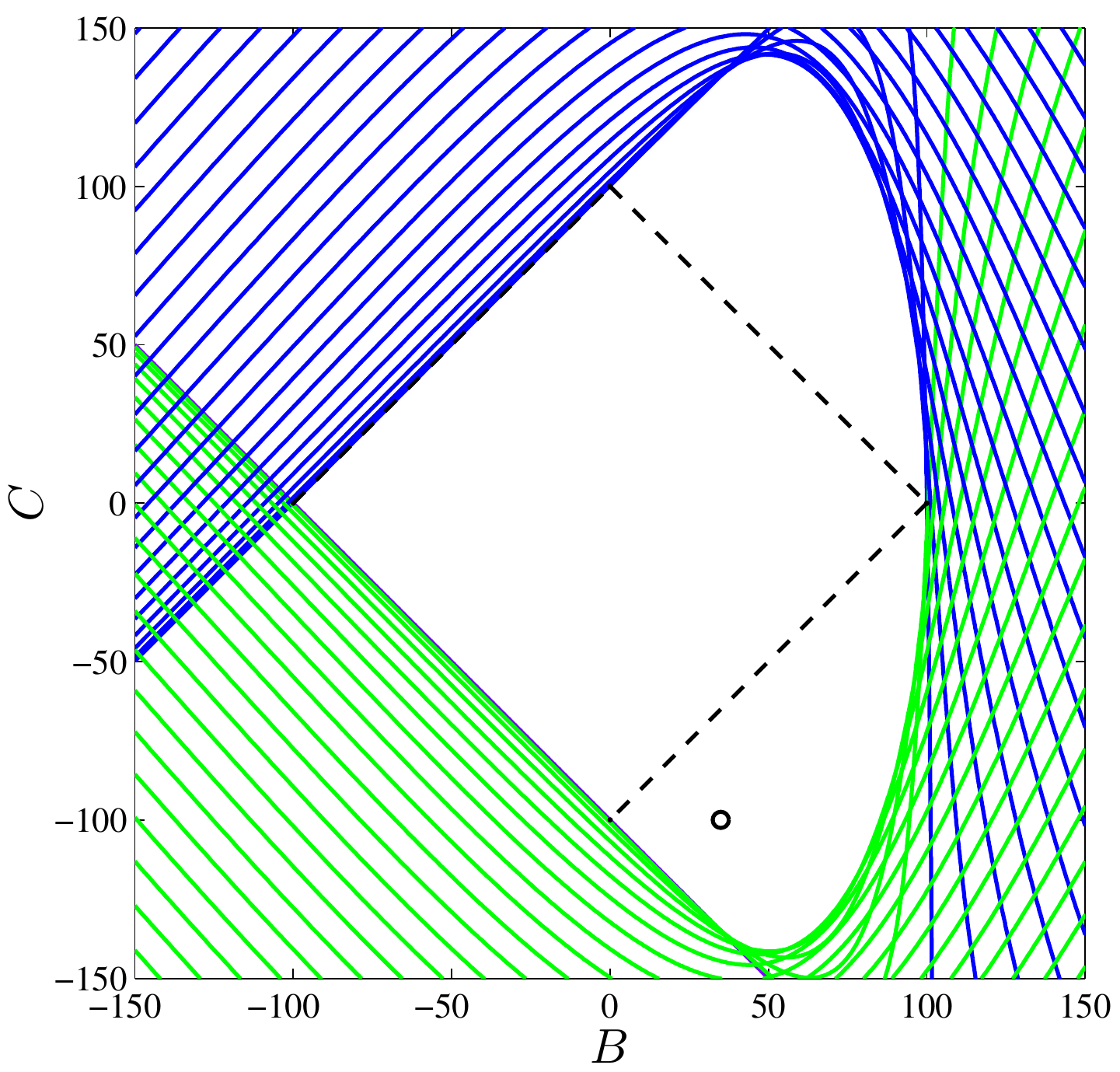}
  &
 \includegraphics[width=0.32\textwidth]{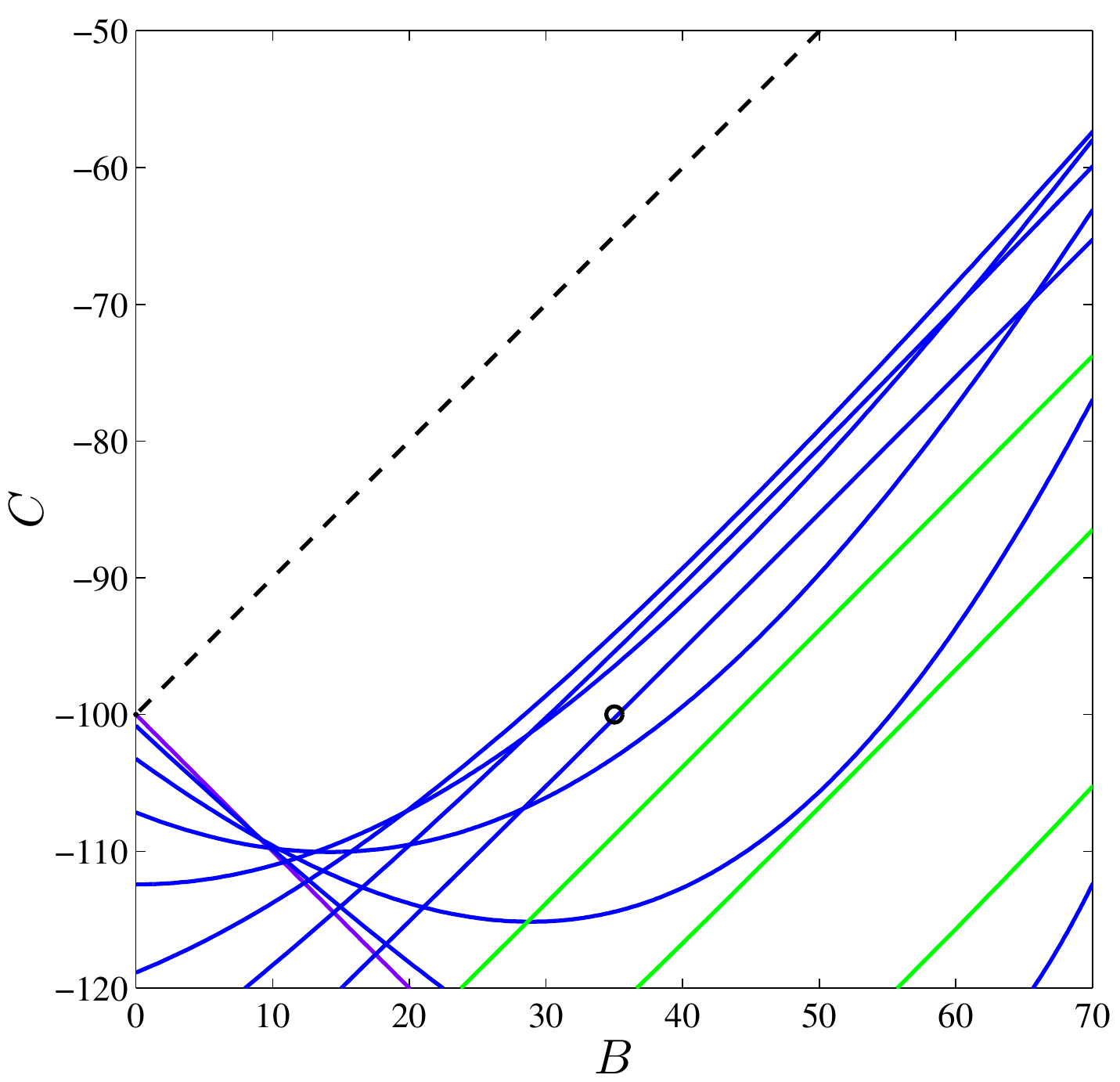}
   &
 \includegraphics[width=0.32\textwidth]{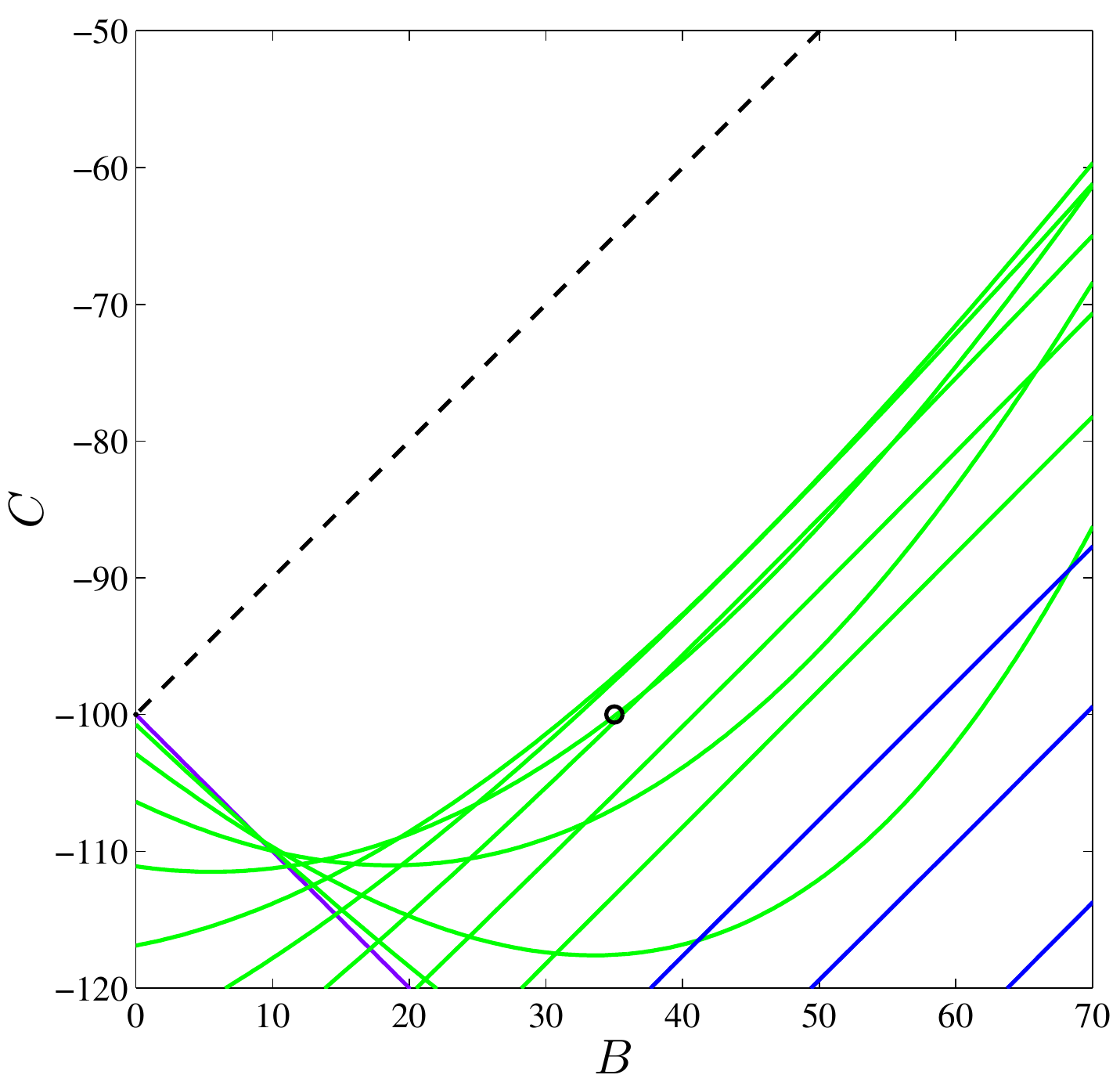}
 \end{tabular}
 \end{center}
  \caption{Bifurcation curves for the modified platelet model near $R = 0.5$ and $A = 100$, showing the point where the model exists in the parameter space. The left figure is at $R = \frac{1}{2}$ with only two families of bifurcation curves, showing the model parameters inside the stable region. The middle figure is for $R = 0.48$, which has $\Gamma_{9}$, $\Gamma_{11}$, and $\Gamma_{13}$ between the model point and the stable region. The right figure with $R = 0.51$ has only $\Gamma_{10}$ and $\Gamma_{12}$ between the model point and the stable region.}
  \label{fig:model_12}
\end{figure}

The modified platelet model also considered delays near $R = \frac{1}{2}$. The simulation in Fig.~\ref{fig:platelet} showed the stability of the solution at $R = \frac{1}{2}$. For $R = \frac{1}{2}$, there are only two families of bifurcation curves, leaving a fairly large region of stability. This is apparent in the leftmost graph of Fig.~\ref{fig:model_12}, where the point for the model parameters is clearly inside the region of stability. When the delay is decreased to $R = 0.48$, there is a transition, $A_1^* = 24.5$. Thus, when $A = 100$, the two family structure has broken down, and the bifurcation curves form a very different pattern. The result is shown in the middle graph of Fig.~\ref{fig:model_12}, where the bifurcation curves $\Gamma_{9}$, $\Gamma_{11}$, and $\Gamma_{13}$ are visible between the model parameter point and the region of stability. Once again, it is easy to compute the eigenvalues with positive real part for this case giving:
\[
  \lambda_1 = 0.07503\pm 64.71\,i \qquad \lambda_2 = 0.05662\pm 77.46\,i \qquad \lambda_3 = 0.04591\pm 51.96\,i.
\]
The leading eigenvalue, $\lambda_1$, has a frequency of 64.71, which suggests a period of 0.09710. This is consistent with the observed period of oscillation in Fig.~\ref{fig:platelet}. We observe that this oscillatory behavior is irregular, which reflects a strong contribution from $\lambda_2$ and $\lambda_3$, which have higher and lower frequencies, respectively. For $R = 0.51$, the rightmost graph of Fig.~\ref{fig:model_12} shows the bifurcation curves $\Gamma_{10}$ and $\Gamma_{12}$ between the model point and the region of stability. When the eigenvalues are computed, we obtain:
\[
  \lambda_1 = 0.03415\pm 60.02\,i \qquad \lambda_2 = 0.02930\pm 72.28\,i.
\]
The frequency of the leading eigenvalue is 60.02, which yields a period of 0.1047. Again, this is consistent with the observed oscillations in Fig.~\ref{fig:platelet}.

\setcounter{equation}{0}
\setcounter{theorem}{0}
\setcounter{figure}{0}
\setcounter{table}{0}

\section{Discussion and Conclusion}

Delay differential equations (DDEs) with multiple time delays are used in a wide
array of applications. When studying the stability of two delay models, our results
show the high sensitivity of the DDE for certain delays (rationally dependent) for
some ranges of the parameters. In particular, the DDE (\ref{DDE2}) shows unusually
large regions of stability for $R = \frac{1}{n}$, when $n$ is a small integer. For
example, when $R = \frac{1}{2}$, the region of stability doubles over the Minimum
Region of Stability (MRS), which is independent of the delay. Our geometric approach
allows a systematic method for visualizing the region of stability and provides a
simplified understanding for how the stability region evolves. In our motivating
example, we demonstrated how easily the leading eigenvalues could be found, which helped
explain some of the observed behavior in the nonlinear problem.

The characteristic equation for (\ref{DDE2}) is an exponential polynomial, which is
deceivingly complex to analyze. For rational delays, $R$, this characteristic equation
organizes into families of curves, which undergo only a few types of reorderings.
The most significant changes occur at values of the parameter, $A$, which we defined
as transitions, $A^*_j$. One interesting phenomenon that can occur at a transition
is a ``stability spur,'' where a region of stability outside the main region of
stability joins, distorting the stability region to become larger. These ``stability
spurs'' also lead to interesting disconnected regions of stability in the $BC$-cross
sectional parameter space. More significantly, as $R \to \frac{1}{n}^-$, the transition
$A_{n-1}^* \to +\infty$, moving the accompanying distortion further away and maintaining
an increased region of stability. We showed that for any $R < \frac{1}{n}$, but close,
the boundary of the stability region in the $BC$-cross section at $A_{n-1}^*$ reduces
primarily to just four simple curves.

The regions of stability, as depicted in Figs.~\ref{fig:R249_tran} and \ref{fig:R333_tran} with
primarily only four curves, allowed us to easily compute the increased area of stability
for delays $R = \frac{1}{n}$ as $A \to +\infty$. The evolution of the stability region
had limited, yet very orderly ways of changing for $R$ near $\frac{1}{n}$. This is
the quasi-periodicity of the families of bifurcation curves, which could self-intersect
mostly through tangencies. The organization of the curves, as seen in many of the figures,
produced clear patterns that could be carefully analyzed for $R = \frac{1}{n}$, resulting
in the observed larger regions of stability.

In this paper we analytically proved some results to support our claims. However, more
analytic work is needed around these singularities that occur in the characteristic
equation for $R = \frac{1}{n}$. Furthermore, other rational delays, like $R = \frac{2}{5}$,
 show similar increases in their regions of stability, but we have not investigated the
 details on how these rationally dependent delays produce larger regions of stability.
 We have produced a framework for future studies of the DDE (\ref{DDE2}) and have excellent
 MatLab programs for further geometric investigations.

 Understanding the stability properties of DDE (\ref{DDE2}) is very important for a number
 of applications with time delays. Our results show that selecting delays of $R = \frac{1}{n}$
 for $n$ small in a model could give the investigator stability that is easily lost with
 only a very small change in the delay. This ultra-sensitivity in the model can be explained
 by our results. This two delay problem is very complex, but our geometric analysis provides
 a valuable tool for future stability analysis of delay models.

\newpage
\bibliography{bibfile}

\begin{thebibliography}{10}

\bibitem{BEL}
J.~B\'elair.
\newblock Stability of a differential-delay equation with two time delays.
\newblock In F.~V. Atkinson, W.~F. Langford, and A.~B. Mingarelli, editors,
  {\em Oscillations, Bifurcations, and Chaos}, volume~8, pages 305--315. AMS,
  1987.

\bibitem{BEM}
J.~B\'elair and M.~Mackey.
\newblock A model for the regulation of mammalian platelet production.
\newblock {\em Ann. N. Y. Acad. Sci.}, 1:1--3, 1987.

\bibitem{BEMb}
J.~B\'elair and M.~Mackey.
\newblock Consumer memory and price fluctuations in commodity markets: {A}n
  integrodifferential model.
\newblock {\em J. Dyn. and Diff. Eqns.}, 3:299--325, 1989.

\bibitem{BMM}
J.~B\'elair, M.~C. Mackey, and J.~M. Mahaffy.
\newblock Age-structured and two delay models for erythropoiesis.
\newblock 1995.

\bibitem{BM01}
J.~B\'elair and J.~M. Mahaffy.
\newblock Variable maturation velocity and parameter sensitivity in a model of
  haematopoiesis.
\newblock {\em IMA J. Math. Appl. Med. \& Biol.}, 18:193--211, 2001.

\bibitem{BeCa}
Jacques B{\'e}lair and Sue~Ann Campbell.
\newblock Stability and bifurcations of equilibria in a multiple-delayed
  differential equation.
\newblock {\em SIAM J. Appl. Math.}, 54(5):1402--1424, 1994.

\bibitem{BCvdD}
Jacques B\'elair, Sue~Ann Campbell, and P.~van~den Driessche.
\newblock Frustration, stability, and delay-induced oscillations in a neural
  network model.
\newblock {\em SIAM Journal on Applied Mathematics}, 56(1):pp. 245--255, 1996.

\bibitem{BelC}
R.~Bellman and K.~L. Cooke.
\newblock {\em Differential-Difference Equations}.
\newblock Academic Press, New York, N.Y., 1963.
\newblock Lectures in Applied Mathematics, Vol. 17.

\bibitem{BOE}
F.~G. Boese.
\newblock The delay-independent stability behaviour of a first order
  differential-difference equation with two constant lags.
\newblock Preprint, 1993.

\bibitem{bohay}
F.~G. Boese.
\newblock A new representation of a stability result of {N}. {D}. {H}ayes.
\newblock {\em Z. Angew. Math. Mech.}, 73(2):117--120, 1993.

\bibitem{Boese94}
F.~G. Boese.
\newblock Stability in a special class of retarded difference-differential
  equations with interval-valued parameters.
\newblock {\em Journal of Mathematical Analysis and Applications}, 181(1):227
  -- 247, 1994.

\bibitem{bortz}
D.~M. Bortz.
\newblock Eigenvalues for two-lag linear delay differential equations.
\newblock (submitted) arXiv:1206.6364, 2012.

\bibitem{BRV}
R.~D. Braddock and P.~van~den Driessche.
\newblock A population model with two time delays.
\newblock In D.~G. Chapman and V.~F. Gallucci, editors, {\em Quantitative
  Population Dynamics}, volume~13 of {\em Statistical Ecology Series}.
  International Cooperative Publishing House, Fairland, MD, 1981.

\bibitem{Busk}
T.~C. Busken.
\newblock On the asymptotic stability of the zero solution for a linear
  differential equation with two delays, 2012.
\newblock Master's Thesis, San Diego State University.

\bibitem{CamBel}
Sue~Ann Campbell and Jacques B{\'e}lair.
\newblock Analytical and symbolically-assisted investigation of {H}opf
  bifurcations in delay-differential equations.
\newblock In {\em Proceedings of the {G}. {J}. {B}utler {W}orkshop in
  {M}athematical {B}iology ({W}aterloo, {ON}, 1993)}, volume~3, pages 137--154,
  1995.

\bibitem{CoY}
K.~L. Cooke and J.~A. Yorke.
\newblock Some equations modelling growth processes and gonorrhea epidemics.
\newblock {\em Math. Biosci.}, 16:75--101, 1973.

\bibitem{ELS}
L.~E. El'sgol'Ts and S.B. Norkin.
\newblock {\em Introduction to the Theory of Differential Equations with
  Deviating Arguments}.
\newblock Academic Press, New York, NY, 1977.

\bibitem{elsken}
Thomas Elsken.
\newblock The region of (in)stability of a 2-delay equation is connected.
\newblock {\em J. Math. Anal. Appl.}, 261(2):497--526, 2001.

\bibitem{GCC}
Cuneyt Guzelis and Leon~O. Chua.
\newblock Stability analysis of generalized cellular neural networks.
\newblock {\em International Journal of Circuit Theory and Applications},
  21(1):1--33, 1993.

\bibitem{HIT}
J.~Hale, E.~Infante, and P.~Tsen.
\newblock Stability in linear delay equations.
\newblock {\em J. Math. Anal. Appl.}, 105:533--555, 1985.

\bibitem{HAL}
J.~K. Hale.
\newblock {\em Theory of Functional Differential Equations}.
\newblock Springer-Verlag, New York, 1977.

\bibitem{Hal1}
J.~K. Hale.
\newblock {\em Nonlinear oscillations in equations with delays}.
\newblock American Math. Soc., Providence, R. I., 1979.
\newblock Lectures in Applied Mathematics, Vol. 17.

\bibitem{HAHw}
J.~K. Hale and W.~Huang.
\newblock Global geometry of the stable regions for two delay differential
  equations.
\newblock {\em J. Math. Anal. Appl.}, 178:344--362, 1993.

\bibitem{HT}
J.~K. Hale and S.~M. Tanaka.
\newblock Square and pulse waves with two delays.
\newblock {\em Journal of Dynamics and Differential Equations}, 12:1--30, 2000.
\newblock 10.1023/A:1009052718531.

\bibitem{HaS}
G.~Haller and G.~St{\'e}p{\'a}n.
\newblock Codimension two bifurcation in an approximate model for delayed robot
  control.
\newblock In R.~Seydel, F.~W. Schneider, Kupper T., and H.~Troger, editors,
  {\em Bifurcation and Chaos: {A}nalysis, Algorithms, Applications}, pages
  155--159. Birkhauser, Basel, 1991.

\bibitem{HAY}
N.~Hayes.
\newblock Roots of the transcendental equation associated with a certain
  differential difference equation.
\newblock {\em J. London Math. Soc.}, 25:226--232, 1950.

\bibitem{HOR}
T.~D. Howroyd and A.~M. Russell.
\newblock Cournot oligopoly models with time lags.
\newblock {\em J. Math. Econ.}, 13:97--103, 1984.

\bibitem{Inf}
E.~F. Infante, 1975.
\newblock Personal Communication.

\bibitem{LEV}
I.~S. Levitskaya.
\newblock Stability domain of a linear differential equation with two delays.
\newblock {\em Comput. Math. Appl.}, 51(1):153--159, 2006.

\bibitem{LRW}
Xiangao Li, Shigui Ruan, and Junjie Wei.
\newblock Stability and bifurcation in delay-differential equations with two
  delays.
\newblock {\em Journal of Mathematical Analysis and Applications}, 236(2):254
  -- 280, 1999.

\bibitem{McDd}
N.~MacDonald.
\newblock Cyclical neutropenia; models with two cell types and two time lags.
\newblock In A.~J. Valleron and P.~D.~M. Macdonald, editors, {\em
  Biomathematics and Cell Kinetics}, pages 287--295. North Holland, Amsterdam,
  1979.

\bibitem{McDc}
N.~MacDonald.
\newblock An activation-inhibition model of cyclic granulopoiesis in chronic
  granulocytic leukemia.
\newblock {\em Math. Biosci.}, 54:61--70, 1980.

\bibitem{MACe}
M.~C. Mackey.
\newblock Commodity price fluctuations: Price dependent delays and
  nonlinearities as explanatory factors.
\newblock {\em J. Econ. Theory}, 48:497--509, 1989.

\bibitem{MAHb}
J.~M. Mahaffy.
\newblock A test for stability of linear differential delay equations.
\newblock {\em Quart. Appl. Math.}, 40:193--202, 1982.

\bibitem{MAHcII}
J.~M. Mahaffy.
\newblock Cellular control models with linked positive and negative feedback
  and delays. {I}. {L}inear analysis and local stability.
\newblock {\em J. theor. Biol.}, 106:103--118, 1984.

\bibitem{MAH85}
J.~M. Mahaffy.
\newblock Stability of periodic solutions for a model of genetic repression
  with delays.
\newblock {\em J. Math. Biol.}, 22:137--144, 1985.

\bibitem{MBM98}
J.~M. Mahaffy, J.~B\'elair, and M.~C. Mackey.
\newblock Hematopoietic model with moving boundary condition and state
  dependent delay: Applications in erythropoiesis.
\newblock {\em J. Theor. Biol.}, 190:135--146, 1998.

\bibitem{MJV}
J.~M. Mahaffy, D.~A. Jorgensen, and R.~L. Vanderheyden.
\newblock Stability results for a model of repression with external control.
\newblock {\em J. Math. Biol.}, 30:669--691, 1992.

\bibitem{MaS}
J.~M. Mahaffy and E.~S. Savev.
\newblock Stability analysis for mathematical models of the {\it lac} operon.
\newblock {\em Quart. Appl. Math}, 57:37--53, 1999.

\bibitem{MZJa}
J.~M. Mahaffy, P.~J. Zak, and K.~M. Joiner.
\newblock A three parameter stability analysis for a linear differential
  equation with two delays.
\newblock Technical report, Department of Mathematical Sciences, San Diego
  State University, San Diego, CA, 1993.

\bibitem{MZJ}
J.~M. Mahaffy, P.~J. Zak, and K.~M. Joiner.
\newblock A geometric analysis of stability regions for a linear differential
  equation with two delays.
\newblock {\em Internat. J. Bifur. Chaos Appl. Sci. Engrg.}, 5(3):779--796,
  1995.

\bibitem{MiI}
M.~Mizuno and K.~Ikeda.
\newblock An unstable mode selection rule: {F}rustrated optical instability due
  to two competing boundary conditions.
\newblock {\em Physica D}, 36:327--342, 1989.

\bibitem{GopMo}
S.~Mohamad and K.~Gopalsamy.
\newblock Exponential stability of continuous-time and discrete-time cellular
  neural networks with delays.
\newblock {\em Applied Mathematics and Computation}, 135(1):17 -- 38, 2003.

\bibitem{MNB}
W.~W. Murdoch, R.~M. Nisbet, S.~P. Blythe, W.~S.~C. Gurney, and J.~D. Reeve.
\newblock An invulnerable age class and stability in delay-differential
  parasitoid-host models.
\newblock {\em American Naturalist}, 129:263--282, 1987.

\bibitem{NUS}
R.~D. Nussbaum.
\newblock A hopf global bifurcation theorem for retarded functional
  differential equations.
\newblock {\em Trans. Amer. Math. Soc.}, 238:139�164, 1978.

\bibitem{PIOT}
M.J. Piotrowska.
\newblock A remark on the ode with two discrete delays.
\newblock {\em Journal of Mathematical Analysis and Applications}, 329(1):664
  -- 676, 2007.

\bibitem{RM}
C.~Grotta Ragazzo and C.~P. Malta.
\newblock Singularity structure of the hopf bifurcation surface of a
  differential equation with two delays.
\newblock {\em Journal of Dynamics and Differential Equations}, 4:617--650,
  1992.
\newblock 10.1007/BF01048262.

\bibitem{RUC}
J.~Ruiz-Claeyssen.
\newblock Effects of delays on functional differential equations.
\newblock {\em J. Diff. Eq.}, 20:404--440, 1976.

\bibitem{saka}
Sadahisa Sakata.
\newblock Asymptotic stability for a linear system of differential-difference
  equations.
\newblock {\em Funkcial. Ekvac.}, 41(3):435--449, 1998.

\bibitem{YTJ}
Toshiaki Yoneyama and Jitsuro Sugie.
\newblock On the stability region of differential equations with two delays.
\newblock {\em Funkcial. Ekvac.}, 31(2):233--240, 1988.

\bibitem{ZAR}
E.~Zaron.
\newblock The delay differential equation: $x'(t) = -ax(t) + bx(t-\tau_1) +
  cx(t-\tau_2)$.
\newblock Technical report, Harvey Mudd College, Claremont, CA, 1987.

\end{thebibliography}

\end{document}